\documentclass[11pt, english]{smfart}
\usepackage{amsfonts, amssymb, amsmath, latexsym}
\usepackage{mathptmx}
\usepackage[english]{babel}
\usepackage{mysectionsen}
\usepackage[latin1]{inputenc}
\usepackage[T1]{fontenc}
\usepackage{textcomp}
\usepackage[dvips]{graphicx}

\input xy
\xyoption{all}

\newcommand{\G}{\mathrm{G}}
\newcommand{\PGL}{\mathrm{PGL}}

\newcommand{\V}{\mathrm{V}}
\newcommand{\W}{\mathrm{W}}

\newcommand{\Br}{\mathrm{B}}
\newcommand{\K}{\mathrm{K}}
\newcommand{\Pp}{\mathbb{P}}
\newcommand{\Qr}{\mathrm{Q}}
\newcommand{\an}{\mathrm{an}}
\newcommand{\rP}{\mathrm{P}}
\newcommand{\U}{\mathrm{U}}
\newcommand{\Sr}{\mathrm{S}}
\newcommand{\R}{\mathrm{R}}

\newcommand{\X}{\mathrm{X}}
\newcommand{\A}{\mathrm{A}}
\newcommand{\T}{\mathrm{T}}
\newcommand{\N}{\mathrm{N}}
\newcommand{\I}{\mathrm{I}}

\newcommand{\Y}{\mathrm{Y}}
\newcommand{\F}{\mathrm{F}}
\newcommand{\Z}{\mathrm{Z}}
\newcommand{\Hr}{\mathrm{H}}
\newcommand{\B}{\mathrm{B}}
\newcommand{\M}{\mathrm{M}}
\newcommand{\Lr}{\mathrm{L}}
\newcommand{\C}{\mathrm{C}}
\newcommand{\D}{\mathrm{D}}

\newcommand{\rad}{\mathrm{rad}}
\newcommand{\radu}{\mathrm{rad}^{\rm u}}

\newcommand{\ov}{\overline}
\newcommand{\inv}{^{-1}}

\newcommand{\Phired}{\Phi^{\rm red}}

\begin{document}

\thispagestyle{empty}

\begin{center} {\Large \bfseries{\textsc{Bruhat-Tits theory from Berkovich's point of view.\\
\vspace{0.2cm}
II. Satake compactifications of buildings}}}
\end{center}

\vspace{1cm}

\begin{center}
\textsc{Bertrand R\'emy, Amaury Thuillier and Annette Werner}
\end{center}

\vspace{1cm}

\begin{center}    July 2009
\end{center}

\vspace{5cm}

\hrule

\vspace{0.5cm}

{\small
\noindent
{\bf Abstract:}
In the paper \emph{Bruhat-Tits theory from Berkovich's point of view. I --- Realizations and compactifications of buildings}, we investigated various realizations of the Bruhat-Tits building $\mathcal{B}(\G,k)$ of a connected and reductive linear algebraic group $\G$ over a non-Archimedean field $k$ in the framework of V. Berkovich's non-Archimedean analytic geometry. We studied in detail the compactifications of the building which naturally arise from this point of view. In the present paper, we give a representation theoretic flavor to these compactifications, following Satake's original constructions for Riemannian symmetric spaces.

We first prove that Berkovich compactifications of a building coincide with the compactifications, previously introduced by the third named author and obtained by a gluing procedure. Then we show how to recover them from an absolutely irreducible linear representation of $\G$ by embedding $\mathcal{B}(\G,k)$ in the building of the general linear group of the representation space, compactified in a suitable way. Existence of such an embedding is a special case of Landvogt's general results on functoriality of buildings, but we also give another natural construction of an equivariant embedding, which relies decisively on Berkovich geometry.

\vspace{0,2cm}


\vspace{0,1cm}

\noindent
{\bf AMS classification (2000):}
20E42,
% Groups with a $BN$-pair; buildings [See also 51E24]
51E24,
% Buildings and the geometry of diagrams
14L15,
% Group schemes
14G22.
% Rigid analytic geometry
}

\vspace{0,5cm}

\hrule

\newpage

\tableofcontents

\newpage

\section*{Introduction}

\textbf{1.} Let $k$ be field a endowed with a complete non-Archimedean absolute value, which we assume to be non-trivial. Let $\G$ be a connected reductive linear algebraic group over $k$. Under some assumptions on $\G$ or on $k$, the Bruhat-Tits building $\mathcal{B}(\G,\K)$ of $\G(\K)$ exists for any non-Archimedean field $\K$ extending $k$ and behaves functorially with respect to $\K$; this is for example the case if $\G$ is quasi-split, or if $k$ is discretely valued with a perfect residue field (in particular, if $k$ is a local field); we refer to \cite[1.3.4]{RTW1} for a discussion. Starting from this functorial existence of the Bruhat-Tits building of $\G$ over any non-Archimedean extension of $k$ and elaborating on some results of Berkovich \cite[Chapter 5]{Ber1}, we explained in \cite{RTW1} how to realize canonically the building $\mathcal{B}(\G,k)$ of $\G(k)$ in some suitable $k$-analytic spaces. The fundamental construction gives a canonical map from the building to the analytification $\G^{\rm an}$ of the algebraic group $\G$, from which one easily deduce another map from $\mathcal{B}(\G,k)$ to $\X^{\rm an}$, where $\X$ stands for any generalized flag variety of $\G$, i.e., a connected component of the projective $k$-scheme ${\rm Par}(\G)$ parametrizing the parabolic subgroups of $\G$. Recall that, if such a connected component $\X$ contains a $k$-rational point $\rP \in {\rm Par}(\G)(k)$, then $\X$ is isomorphic to the quotient scheme $\G/\rP$. In more elementary words, this simply means that $\mathcal{B}(\G,k)$ has a natural description in terms of multiplicative seminorms (of homothety classes of multiplicative seminorms, respectively) on the coordinate ring of $\G$ (on the homogeneous coordinate ring of any connected component of ${\rm Par}(\G)$, respectively).

Since the algebraic scheme ${\rm Par}(\G)$ is projective, the topological space underlying the analytification ${\rm Par}_t(\G)^{\rm an}$ of any connected component ${\rm Par}_t(\G)$ of ${\rm Par}(\G)$ is compact (that is, Hausdorff and quasi-compact), hence can be used to compactify $\mathcal{B}(\G,k)$ by passing to the closure (in a suitable sense if $k$ is not locally compact). In this way, one associates with each connected component ${\rm Par}_t(\G)$ of ${\rm Par}(\G)$ a compactified building $\overline{\mathcal{B}}_t(\G,k)$, which is a $\G(k)$-topological space containing some factor of $\mathcal{B}(\G,k)$ as a dense open subset. There is no loss of generality in restricting to connected components of ${\rm Par}(\G)$ having a $k$-rational point, i.e., which are isomorphic to $\G/\rP$ for some parabolic subgroup $\rP$ of $\G$ (well-defined up to $\G(k)$-conjugacy). Strictly speaking, $\overline{\mathcal{B}}_t(\G,k)$ is a compactification of $\mathcal{B}(\G,k)$ only if $k$ is a local field and if the conjugacy class of parabolic subgroups corresponding to the component ${\rm Par}_t(\G)$ of ${\rm Par}(\G)$ is \emph{non-degenerate}, i.e., consists of parabolic subgroups which do not contain a full almost simple factor of $\G$; however, we still refer to this enlargement of $\mathcal{B}(\G,k)$ as a "compactification" even if these conditions are not fulfilled. The compactified building $\overline{\mathcal{B}}_t(\G,k)$ comes with a canonical stratification into locally closed subspaces indexed by a certain set of parabolic subgroups of $\G$. The stratum attached to a parabolic subgroup $\rP$ is isomorphic to the building of the semi-simplification $\rP/{\rm rad}(\rP)$ of $\rP$, or rather to some factors of it. We obtain in this way one compactified building for each $\G(k)$-conjugacy class of parabolic subgroups of $\G$.

\vspace{0.3cm}
\textbf{2.} Assuming that $k$ is a local field, the third named author had already defined a compactification of $\mathcal{B}(\G,k)$ for each conjugacy class of parabolic subgroup of $\G$, see \cite{Wer2}. Inspired by Satake's approach for Riemannian symmetric spaces, the construction in [loc.cit] starts with an absolutely irreducible (faithful) linear representation $\rho$ of $\G$ and consists of two steps: \begin{itemize} \item[(i)] the apartment $\A(\Sr,k)$ of a maximal split torus $\Sr$ of $\G$ in $\mathcal{B}(\G,k)$ is compactified, say into $\overline{\A}(\Sr,k)_{\rho}$, by using the same combinatorial analysis of the weights of $\rho$ as in \cite{Sa};
\item[(ii)] the compactified building $\overline{\mathcal{B}}(\G,k)_{\rho}$ is defined as the quotient of $\G(k) \times \overline{\A}(\Sr,k)_{\rho}$ by a suitable extension of the equivalence relation used by Bruhat and Tits to construct $\mathcal{B}(\G,k)$ as a quotient of $\G(k) \times \A(\Sr,k)$.
\end{itemize}

It is proved in [loc.cit] that the so-obtained compactified building only depends on the position of a highest weight of $\rho$ with respect to Weyl chambers, or equivalently on the conjugacy class of parabolic subgroups of $\G$ stabilizing the line spanned by a vector of highest weight. As suggested in [loc.cit], these compactifications turn out to coincide with Berkovich ones.

\vspace{0.1cm} Let us define the \emph{type} $t(\rho)$ of an absolutely irreducible linear representation $\rho : \G \rightarrow {\rm GL}_\V$ as follows. If $\G$ is split, then each Borel subgroup $\B$ of $\G$ stabilizes a unique line $\Lr_{\B}$ in $\V$, its \emph{highest weight line}. One easily shows that there exists a largest parabolic subgroup $\rP$ of $\G$ stabilizing the line $\Lr_{\B}$. Now, the type $t(\rho)$ of the representation $\rho$ is characterized by the following condition: for any finite extension $k'/k$ splitting $\G$, the connected component ${\rm Par}_{t(\rho)}(\G)$ of ${\rm Par}(\G)$ contains each $k'$-point occurring as the largest parabolic subgroup of $\G \otimes_k k'$ stabilizing a highest weight line in $\V \otimes_k k'$.  Finally, the \emph{cotype} of the representation $\rho$ is defined as the type of the contragredient representation $\check{\rho}$. We establish in Section 2, Theorem \ref{thm.comparison}, the following comparison.

\vspace{0.2cm}
\noindent \emph{\textbf{Theorem 1}} --- \emph{Let $\rho$ be an absolutely irreducible (faithful) linear representation of $\G$ in some finite-dimensional vector space over $k$. Then the compactifications  $\overline{\mathcal{B}}(\G,k)_{\rho}$ and $\overline{\mathcal{B}}_{t(\rho)}(\G,k)$ of the building $\mathcal{B}(\G,k)$ are canonically isomorphic.}

\vspace{0.3cm}
\textbf{3.} We still assume that $k$ is a local field but the results below hold more generally for a discretely valued non-Archimedean field with perfect residue field. Another way to compactify buildings by means of linear representations consists first in compactifying the building of the projective linear group ${\rm PGL}_{\V}$ of the representation space and then using a representation in order to embed $\mathcal{B}(\G,k)$ into this compactified building. Finally, a compactification of $\mathcal{B}({\rm PGL}_{\V},k)$ can be obtained by embedding this building in some projective space, hence this viewpoint is the closest one in spirit to the original approach for symmetric spaces. It is also a way to connect Bruhat-Tits theory to Berkovich's interpretation of the space of seminorms on a given $k$-vector space \cite{Ber3}.

More precisely, let $\rho: \G \to {\rm GL}_{\V}$ be an absolutely irreducible linear representation of $\G$ in a finite-dimensional $k$-vector space $\V$. We use such a map $\rho$ in two ways to obtain continuous $\G(k)$-equivariant maps from the building $\mathcal{B}(\G,k)$ to a compact space $\mathcal{X}(\V,k)$ naturally attached to the $k$-vector space $\V$. Denoting by $\mathcal{S}(\V,k)$ the "extended Goldman-Iwahori space" consisting of non-zero seminorms on $\V$ (the space of norms was studied in \cite{GoldmanIwahori}), then the space $\mathcal{X}(\V,k)$ is the quotient of $\mathcal{S}(\V,k)$ by homotheties. It is the non-Archimedean analogue of the quotient of the cone of positive (possibly degenerate) Hermitian matrices in the projective space associated with ${\rm End}(\V)$ \cite{Sa}.

In the real case, the latter space is classically the target space of a suitable Satake map. In our case, we identify $\mathcal{X}(\V,k)$ with the compactification $\overline{\mathcal{B}}_{\delta}({\rm PGL}_\V,k)$ corresponding to the type $\delta$ of parabolic subgroups stabilizing a \emph{hyperplane} of $\V$. One could also consider the compactified building $\overline{\mathcal{B}}_{\pi}({\rm PGL}_\V, k)$ associated with the type $\pi$ of parabolic subgroups stabilizing a \emph{line} of $\V$ (see \cite{Wer0}). Note that $\overline{\mathcal{B}}_{\delta}({\rm PGL}_\V,k) \cong \overline{\mathcal{B}}_{\pi}({\rm PGL}_{\V^{\vee}},k)$, where $\V^{\vee}$ is the dual of $\V$.

\vspace{0.1cm} A first way to obtain a map $\mathcal{B}(\G,k) \to \mathcal{X}(\V,k)$ is to make use of E. Landvogt's work on the functoriality of Bruhat-Tits buildings (with respect both to the group and to the field). Indeed, specializing the results of \cite{LandvogtCrelle} to $k$-homomorphisms arising from linear representations $\rho: \G \to {\rm GL}_{\V}$, we obtain a (possibly non-uniquely defined) map $\rho_* : \mathcal{B}(\G,k) \rightarrow \mathcal{B}({\rm PGL}_{\V},k)$ between buildings. We can then compose it with the compactification map $\vartheta_\pi : \mathcal{B}({\rm PGL}_{\V},k) \to \overline{\mathcal{B}}_{\pi}({\rm PGL}_\V,k)$ in order to obtain an analogue of a Satake map.

There is another way to embed the building $\mathcal{B}(\G,k)$ into $\mathcal{X}(\V,k)$, which turns out to be very natural and relies crucially on Berkovich geometry. There exists a natural $k$-morphism $\widetilde{\rho}$ from the scheme ${\rm Bor}(\G)$ of Borel subgroups of $\G$ to the projective space $\Pp(\V)$ satisfying the following condition: for any extension $\K/k$, the map $\widetilde{\rho}_{\K}$ sends a Borel subgroup $\B$ of $\G \otimes_k \K$ to the unique $\K$-point $\widetilde{\rho}(\B)$ of $\mathbb{P}(\V)$ it fixes. By passing to analytic spaces, we get a map $\widetilde{\rho}: {\rm Bor}(\G)^{\rm an} \to \Pp(\V)^{\rm an}$.
Using the concrete description of $\mathcal{X}(\V,k)$ and $\mathbb{P}(\V)^{\rm an}$, we have a natural retraction $\tau: \Pp(\V)^{\an} \rightarrow \mathcal{X}(\V,k)$, so that the composition $\underline\rho = \tau \circ \widetilde{\rho} \circ \vartheta_\varnothing$ sends the Bruhat-Tits building $\mathcal{B}(\G,k)$ into $\mathcal{X}(\V,k)$.
This is our second way to obtain a non-Archimedean analogue of a Satake map, and it is easily seen that this canonical map sends an apartment into an apartment.

\vspace{0.1cm}
These two embedding procedures lead to the previous families of compactifications (cf. Theorem \ref{thm.Berkovich.Satake} and Theorem \ref{th-Satake_via_embeddings}):

\vspace{0.2cm}
\noindent \emph{\textbf{Theorem 2}} ---
\emph{Assume that $k$ is a non-Archimedean local field and let $\rho: \G \to {\rm GL}_\V$ be an absolutely irreducible linear representation of $\G$ in a finite-dimensional vector space $\V$ over $k$. \begin{itemize} \item[(i)] The map $\underline\rho: \mathcal{B}(\G,k) \rightarrow \mathcal{X}(\V,k)$ induces a $\G(k)$-equivariant homeomorphism between $\overline{\mathcal{B}}_{t(\check{\rho})}(\G,k)$ and the closure of the image of $\underline\rho$ in $\mathcal{X}(\V,k)$.
\item[(ii)] Any Landvogt map $\rho_\ast: \mathcal{B}(\G,k) \rightarrow \mathcal{B}({\rm PGL}_\V,k)$ induces a $\G(k)$-equivariant homeomorphism between $\overline{\mathcal{B}}_{t(\rho)}(\G,k)$ and the closure of its image in $\overline{\mathcal{B}}_{\pi}({\rm PGL}_\V,k)$.
\end{itemize}}

\

\vspace{0.2cm} \textbf{Conventions.} Assumptions on the field $k$ are made explicit at the beginning of each section. Notations and conventions from \cite{RTW1} are recalled in section 1. 

Let us stress one particular working hypothesis: the results in \textbf{[loc.cit]} were obtained under a functoriality assumption for buildings with respect to non-Archimedean extension of the ground field (see \textbf{[loc.cit, 1.3.4]} for a precise formulation). This assumption, which is fulfiled in particular if $k$ is discretely valued with perfect residue field or if the group under consideration is split, is made throughout the present work.  

\vspace{0.2cm} \textbf{Structure of the paper.} In the first section, we briefly review the constructions of \cite{RTW1} and state the results from [loc.cit] to be used in this work. The second section is devoted to the identification of Berkovich compactifications with the compactifications introduced in \cite{Wer2}. The third section contains a concrete description of the Berkovich compactification of the building $\mathcal{X}(\V,k) = \mathcal{B}({\rm PGL}_\V,k)$ associated with the projective space $\mathbb{P}(\V)$ seen as a generalized flag variety. The last two sections deal with the recovery of Berkovich compactifications via embeddings into $\mathcal{X}(\V,k)$, in the spirit of Satake's original construction for Riemannian symmetric spaces. In Section 4, we construct a canonical $\G(k)$-map from $\mathcal{B}(\G,k)$ to $\mathcal{X}(\V,k)$ for each absolutely irreducible linear representation of $\G$ in $\V$, and we show that taking the closure leads to the Berkovich compactification of $\mathcal{B}(\G,k)$ of type $t(\check{\rho})$. In Section 5, we rely on Landvogt's functoriality results to produce such a map and derive the same conclusion.

\newpage

\section{Berkovich compactifications of buildings}

This section provides a brief summary of realizations and compactifications of Bruhat-Tits buildings in the framework of Berkovich's non-Archimedean analytic geometry. We refer to \cite{RTW1} for proofs, details and complements.

\vspace{0.1cm} In the following, we consider a non-Archimedean field $k$, i.e., a field endowed with a complete non-Archimedean absolute value which we assume to be non-trivial, and a semisimple and connected linear $k$-group $\G$. 

\vspace{0.2cm} \noindent \textbf{(1.1)} For each point $x$ of the Bruhat-Tits building $\mathcal{B}(\G,k)$, there exists a unique affinoid subgroup $\G_x$ of $\G^{\rm an}$ satisfying the following condition: for any non-Archimedean extension $\K/k$, the group $\G_x(\K)$ is the stabilizer of $x_{\K}$ in $\G(\K)$, where $x_{\K}$ denotes the image of $x$ under the natural injection $\mathcal{B}(\G,k) \hookrightarrow \mathcal{B}(\G,\K)$.
Seen as a set of multiplicative seminorms on the coordinate algebra $\mathcal{O}(\G)$ of $\G$, the subspace $\G_x$ contains a unique maximal point, denoted by $\vartheta(x)$. One can recover $\G_x$ from $\vartheta(x)$ as its holomorphic envelope:  $$\G_x = \{z \in \G^{\rm an} \ ; \ |f|(z) \leqslant |f|(\vartheta(x)) \ \textrm{ for all } f \in \mathcal{O}(\G)\}.$$

We have thus defined a map $$\vartheta : \mathcal{B}(\G,k) \rightarrow \G^{\rm an}$$ which is continuous, injective and $\G(k)$-equivariant with respect to the $\G(k)$-action by conjugation on $\G^{\rm an}$. By its very construction $\vartheta$ is compatible with non-Archimedean extensions of $k$.

\vspace{0.2cm} \noindent \textbf{(1.2)} We let ${\rm Par}(\G)$ denote the $k$-scheme of parabolic subgroups of $\G$; this is a smooth and projective scheme representing the functor $$\mathbf{Sch}/k \rightarrow \mathbf{Sets},  \ \ \Sr \mapsto \{\textrm{parabolic subgroups of } \G \times_k \Sr\}.$$
The connected components of ${\rm Par}(\G)$ are naturally in bijection with ${\rm Gal}(k^a|k)$-stable subsets of vertices in the Dynkin diagram of $\G \otimes_k k^a$. Such a subset $t$ is called a \emph{type} of parabolic subgroups of $\G$ and we denote by ${\rm Par}_t(\G)$ the corresponding connected component of ${\rm Par}(\G)$. For example, ${\rm Par}_{\varnothing}(\G)$ is the scheme of Borel subgroups of $\G$ whereas the trivial type corresponds to the maximal parabolic subgroup $\G$. Finally, a type $t$ is said to be \emph{k-rational} if ${\rm Par}_t(\G)(k) \neq \varnothing$, i.e., if there exists a parabolic subgroup of $\G$ of type $t$.

\vspace{0.1cm} With each parabolic subgroup $\rP$ of $\G$ is associated a morphism $\omega_{\rP} : \G \rightarrow {\rm Par}(\G)$, defined functor-theoretically by $g \mapsto g \rP g^{-1}$ and inducing an isomorphism from $\G/\rP$ to the (geometrically) connected component of ${\rm Par}(\G)$ containing the $k$-point $\rP$. Composing $\vartheta$ with the analytification of $\omega_{\rP}$, we obtain a continuous and $\G(k)$-equivariant map from $\mathcal{B}(\G,k)$ to ${\rm Par}(\G)^{\rm an}$ which depends only on the type $t$ of $\rP$. This map is denoted by $\vartheta_t$ and its image lies in the connected component ${\rm Par}_t(\G)^{\rm an}$ of ${\rm Par}(\G)^{\rm an}$. The map $\vartheta_t$ only depends on the type $t$, not on the choice of $\rP$ in ${\rm Par}_t(\G)(k)$. It is defined more generally for any type $t$ of parabolic subgroups, even non-$k$-rational ones; however, we restrict to $k$-rational types in this section.

\vspace{0.1cm} The topological space underlying ${\rm Par}(\G)^{\rm an}$ is compact, hence leads to compactifications of the building $\mathcal{B}(\G,k)$ by closing. From now on, we fix a $k$-rational type $t$ and describe the corresponding compactification of $\mathcal{B}(\G,k)$. If $\Sr$ is a maximal split torus of $\G$, we recall that $\A(\Sr,k)$ denotes the corresponding apartment in the building $\mathcal{B}(\G,k)$.

\begin{Def} \label{def.compactifications} For any maximal split torus $\Sr$ of $\G$, we let $\overline{\A}_t(\Sr,k)$ denote the closure of $\vartheta_t(\A(\Sr,k))$ in ${\rm Par}(\G)^{\rm an}$. We set $$\overline{\mathcal{B}}_t(\G,k) = \bigcup_{\Sr} \overline{\A}_t(\Sr,k) \subset {\rm Par}(\G)^{\rm an},$$ where the union is taken over the set of maximal split tori of $\G$. This is a $\G(k)$-invariant subset of ${\rm Par}(\G)^{\rm an}$, which we endow with the quotient topology induced by the natural $\G(k)$-equivariant map $$\G(k) \times \overline{\A}_t(\Sr,k) \rightarrow \overline{\mathcal{B}}_t(\G,k).$$
\end{Def}

(See \cite[Definition 3.30]{RTW1}.)

\vspace{0.1cm}
The type $t$ is said to be \emph{non-degenerate} if it restricts non-trivially to each almost simple factor of $\G$, i.e., if $t$, seen as a ${\rm Gal}(k^a|k)$-stable set of vertices in the Dynkin diagram $\D$ of $\G \otimes_k k^a$, does not contain any connected component of $\D$. In general, there exist two semisimple groups $\Hr'$, $\Hr''$ and a central isogeny $\G \rightarrow \Hr' \times \Hr''$ such that $t$ has non-degenerate restriction to $\Hr'$ and trivial restriction to $\Hr''$. In this situation, $\mathcal{B}(\G,k) \cong \mathcal{B}(\Hr',k) \times \mathcal{B}(\Hr'',k)$ and we let $\mathcal{B}_t(\G,k)$ denote the factor $\mathcal{B}(\Hr',k)$.

\begin{Prop} \label{prop.compactifications}
\begin{itemize}
\item[(i)] The map $\vartheta_t : \mathcal{B}(\G,k) \rightarrow {\rm Par}(\G)^{\rm an}$ factors through the canonical projection of $\mathcal{B}(\G,k)$ onto $\mathcal{B}_t(\G,k)$ and induces an injection of the latter building in ${\rm Par}(\G)^{\rm an}$.
\item[(ii)] If the field $k$ is locally compact, then $\overline{\mathcal{B}}_t(\G,k)$ is the closure of $\vartheta_t\left(\mathcal{B}(\G,k)\right)$ in ${\rm Par}(\G)^{\rm an}$, endowed with the induced topology.
\end{itemize}
\end{Prop}

(See \cite[Proposition 3.34]{RTW1}.)

\vspace{0.1cm}

If $k$ is not locally compact, the topological space $\overline{\mathcal{B}}_t(\G,k)$ is not compact. However, the map $\vartheta_t : \mathcal{B}_t(\G,k) \hookrightarrow \overline{\mathcal{B}}_t(\G,k)$ still induces a homeomorphism onto an open dense subset of $\overline{\mathcal{B}}_t(\G,k)$.

\vspace{0.2cm} \noindent \textbf{(1.3)} The topological space $\overline{\mathcal{B}}_t(\G,k)$ carries a canonical stratification whose strata are lower-dimensional buildings coming from semisimplications of suitable parabolic subgroups of $\G$.

We can attach to each parabolic subgroup $\Qr$ of $\G$ a closed and smooth subscheme ${\rm Osc}_t(\Qr)$ of ${\rm Par}_t(\G)$, homogeneous under $\Qr$ and representing the subfunctor $$\mathbf{Sch}/k \rightarrow \mathbf{Sets}, \ \ \Sr \mapsto \left\{\begin{array}{c} \textrm{parabolic subgroups of } \G \times_k \Sr \\ \textrm{of type $t$, osculatory with } \Qr \times_k \Sr \end{array}\right\}.$$ We recall that two parabolic subgroups of a reductive $\Sr$-group scheme are osculatory if, \'etale locally on $\Sr$, they contain a common Borel subgroup. Letting $\Qr_{\rm ss}$ denote the semisimple $k$-group $\Qr/{\rm rad}(\Qr)$, the morphism $\iota_{\Qr} : {\rm Osc}_t(\Qr) \rightarrow {\rm Par}_t(\Qr_{\rm ss})$ defined functor-theoretically by $\rP \mapsto (\rP \cap \Qr)/{\rm rad}(\Qr)$ is an isomorphism.

\vspace{0.1cm}
There exists a largest parabolic subgroup $\Qr'$ stabilizing ${\rm Osc}_t(\Qr)$. By construction, we have $\Qr \subset \Qr'$ and ${\rm Osc}_t(\Qr') = {\rm Osc}_t(\Qr)$, and we say that $\Qr$ is \emph{t-relevant} if $\Qr = \Qr'$. In general, $\Qr'$ is the smallest $t$-relevant parabolic subgroup of $\G$ containing $\Qr$.

\begin{Ex} a) It $t_{\rm min}$ denotes the type of minimal parabolic subgroups of $\G$, then each parabolic subgroup of $\G$ is $t_{\rm min}$-relevant. Indeed, for any two parabolic subgroups $\rP$ and $\Qr$ such that $\Qr \subsetneq \rP$, there exists a minimal parabolic subgroup contained in $\rP$ but not in $\Qr$; this implies ${\rm Osc}_{t_{\rm min}}(\Qr) \neq {\rm Osc}_{t_{\rm min}}(\rP)$, hence $\Qr$ is the largest parabolic subgroup stabilizing ${\rm Osc}_{t_{\rm min}}(\Qr)$.

b) Let $\V$ be a finite-dimensional $k$-vector space. We assume that $\G = {\rm PGL}_{\V}$ and that $\delta$ is the type of parabolic subgroups of ${\rm PGL}_{\V}$ stabilizing a hyperplane. In this case, ${\rm Par}_\delta(\G)$ is the projective space $\mathbb{P}(\V)$, i.e., the scheme of hyperplanes in $\V$. Each parabolic subgroup $\Qr$ of ${\rm PGL}_{\V}$ is the stabilizer of a well-defined flag $\V^\bullet$ of linear subspaces, and two parabolic subgroups are osculatory if and only if the corresponding flags admit a common refinement, i.e., are subflags of the same flag. It follows that ${\rm Osc}_{\delta}(\Qr)$ is the closed subscheme $\mathbb{P}(\V/\W)$ of $\mathbb{P}(\V)$, where $\W$ is the largest proper linear subspace of $\V$ occurring in the flag $\V^{\bullet}$, and therefore $\delta$-relevant parabolic subgroups of ${\rm PGL}_{\V}$ are precisely the stabilizers of flags $(\{0\} \subset \W \subset \V)$, where $\W$ is any linear subspace of $\V$.
\end{Ex}

\vspace{0.1cm} We can now describe the canonical stratification on the compactified building $\overline{\mathcal{B}}_t(\G,k)$.

\begin{Thm} \label{thm.stratification} For any parabolic subgroup $\Qr$ of $\G$, we use the map $\iota_{\Qr}^{-1} \circ \vartheta_t$ to embed $\mathcal{B}_t(\Qr_{\rm ss},k)$ into ${\rm Osc}_t(\Qr)^{\rm an} \subset {\rm Par}_t(\G)^{\rm an}$.
\begin{itemize}
\item[(i)] As a subset of ${\rm Par}_t(\G)^{\rm an}$, the building $\mathcal{B}(\Qr_{\rm ss},k)$ is contained in $\overline{\mathcal{B}}_t(\G,k)$.
\item[(ii)] We have the following stratification by locally closed subsets: $$\overline{\mathcal{B}}_t(\G,k) = \bigsqcup_{\textrm{t-relevant } \Qr's} \mathcal{B}_t(\Qr_{\rm ss},k),$$ where the union is indexed by the $t$-relevant parabolic subgroups of $\G$. The closure of the stratum $\mathcal{B}_t(\Qr_{\rm ss},k)$ is the union of all strata $\mathcal{B}_t(\rP_{\rm ss},k)$ with $\rP \subset \Qr$ and is canonically homeomorphic to the compactified building $\overline{\mathcal{B}}_t(\Qr_{\rm ss},k)$.
\end{itemize}
\end{Thm}

(See \cite[Theorem 4.1]{RTW1}.)

\vspace{0.1cm}
\begin{Ex} a) Suppose that $t = t_{\rm min}$ is the type of minimal parabolic subgroups of $\G$. This type is non-degenerate and each parabolic subgroup of $\G$ is $t_{\rm min}$-relevant, hence the boundary of $\overline{\mathcal{B}}_{t_{\rm min}}(\G,k)$ contains a copy of the building of $\Qr_{\rm ss}$ for each proper parabolic subgroup $\Qr$ of $\G$.

b) Let $\V$ be a finite-dimensional $k$-vector space. We assume that $\G = {\rm PGL}_{\V}$ and that $t = \delta$ is the type of parabolic subgroups of ${\rm PGL}_{\V}$ stabilizing a hyperplane. In this case, the boundary of $\overline{\mathcal{B}}_\delta({\rm PGL}_{\V},k)$ is the union of the buildings $\mathcal{B}({\rm PGL}(\V/\W), k)$, where $\W$ runs over the set of proper non-zero linear subspaces of $\V$.
\end{Ex}

\vspace{0.2cm} \noindent \textbf{(1.4)} We now look at the compactified apartment $\overline{\A}_t(\Sr,k)$ of a maximal split torus $\Sr$ of $\G$. The apartment $\A(\Sr,k)$ is an affine space under the vector space $\V(\Sr) = {\rm Hom}_{\mathbf{Ab}}(\X^*(\Sr),\mathbb{R})$, where $\X^*(\Sr) = {\rm Hom}_{k-\mathbf{Gr}}(\Sr, \mathbb{G}{\rm m}_k)$ is the group of characters of $\Sr$. Let $\Phi = \Phi(\G,\Sr) \subset \X^*(\Sr)$ denote the set of roots of $\G$ with respect to $\Sr$. With each parabolic subgroup $\rP$ of $\G$ containing $\Sr$ we associate its \emph{Weyl cone} $$\mathfrak{C}(\rP) = \{u \in \V(\Sr) \ ; \ \langle \alpha, u \rangle \geqslant 0 \ \textrm{ for all roots } \alpha \textrm{ of } \rP\},$$ which is a strictly convex rational polyhedral cone in $\V(\Sr)$. The collection of Weyl cones of parabolic subgroups of $\G$ containing $\Sr$ is a complete \emph{fan} on the vector space $\V(\Sr)$, i.e., a finite family of strictly convex rational polyhedral cones stable under intersection, in which any two cones intersect along a common face, and satisfying the additional condition that $\V(\Sr)$ is covered by the union of these cones.

\vspace{0.1cm} Relying on the $k$-rational type $t$, we can define a new complete fan on $\V(\Sr)$, which we denote by $\mathcal{F}_t$. The fan of Weyl cones will turn out to be $\mathcal{F}_{t_{\rm min}}$. First of all, if $\rP$ is a parabolic subgroup of type $t$ containing $\Sr$, we  define $\C_t(\rP)$ as the "combinatorial neighborhood" of $\mathfrak{C}(\rP)$ in $\V(\Sr)$, i.e., $$\C_t(\rP) = \bigcup_{\tiny \begin{array}{c} \Qr \ \textrm{parabolic} \\ \Sr \subset \Qr \subset \rP \end{array}} \mathfrak{C}(\Qr).$$ This is a convex rational polyhedral cone, and $\C_t(\rP)$ is strictly convex if and only if the type $t$ is non-degenerate. More precisely, the central isogeny $\G \rightarrow \Hr' \times \Hr''$ introduced after Definition \ref{def.compactifications} corresponds to a decomposition of $\Phi$ as the union $\Phi' \cup \Phi''$ of two closed and disjoint subsets, and the largest linear subspace of $\C_t(\rP)$ is the vanishing locus of $\Phi''$, namely $$\langle \Phi'' \rangle = \{ u \in \X^*(\Sr) \ ; \ \langle \alpha, u \rangle = 0 \ \textrm{ for all } \alpha \in \Phi''\}.$$ When $\rP$ runs over the set of parabolic subgroups of $\G$ of type $t$ and containing $\Sr$, one checks that the set $\mathcal{F}_t$, consisting of the cones $\C_t(\rP)$ together with their faces, induces a complete fan on the quotient space $\V(\Sr)/\langle \Phi'' \rangle$.

\vspace{0.1cm} Any strictly convex rational polyhedral cone $\C$ in $\V(\Sr)$ has a canonical compactification $\overline{\C}$, whose description is nicer if we switch to multiplicative notation for the real dual of $\X^*(\Sr)$. Hence, we set $\Lambda(\Sr) = {\rm Hom}_{\mathbf{Ab}}(\X^{*}(\Sr),\mathbb{R}_{>0})$ and use the isomorphism $\mathbb{R} \rightarrow \mathbb{R}_{>0}, \ x \mapsto e^x$ in order to identify $\V(\Sr)$ with $\Lambda(\Sr)$.

Let $\M$ denote the set of characters $\chi \in \X^*(\Sr)$ such that $\langle \chi, u \rangle \leqslant 1$ for any $u \in \C \subset \Lambda(\Sr)$. This is a finitely generated semigroup of $\X^*(\Sr)$ and the map $$\C \rightarrow {\rm Hom}_{\mathbf{Mon}}(\M, ]0,1]), \ \ u \mapsto (\chi \mapsto \langle \chi, u\rangle)$$ identifies $\C$ with the set ${\rm Hom}_{\mathbf{Mon}}(\M, ]0,1])$ of morphisms of unitary monoids, endowed with the coarsest topology making each evaluation map continuous. We define $\overline{\C}$ as the set ${\rm Hom}_{\mathbf{Mon}}(\M,[0,1])$ endowed with the analogous topology; this is a compact space in which $\C$ embeds as an open dense subspace. Each complete fan $\mathcal{F}$ of strictly convex rational polyhedral cones on $\Lambda(\Sr)$ gives rise to a compactification $\overline{\Lambda}(\Sr)^{\mathcal{F}}$ of this vector space, defined by gluing together the compactifications of the cones $\C \in \mathcal{F}$. More generally, one can compactify in this way any affine space under $\Lambda(\Sr)$.

\begin{Prop} \label{prop.fan.compactification} Let $\Sr$ be a maximal split torus of $\G$. The compactified apartment $\overline{\A}_t(\Sr,k)$ is canonically homeomorphic to the compactification of $\A(\Sr,k)/\langle \Phi'' \rangle$ associated with the complete fan $\mathcal{F}_t$.
\end{Prop}

(See \cite[Proposition 3.35]{RTW1}.)

\vspace{0.1cm} The connection between $t$-relevant parabolic subgroups on the one hand and cones belonging to $\mathcal{F}_t$ on the other hand is the following.

\begin{Prop} \label{prop.relevant.fan} For each parabolic subgroup $\Qr$ of $\G$ containing $\Sr$, there is a smallest cone $\C_t(\Qr)$ in $\mathcal{F}_t$ containing the Weyl cone $\mathfrak{C}(\Qr)$. The following two conditions are equivalent:
\begin{itemize}
\item[(i)] $\Qr$ is $t$-relevant;
\item[(ii)] $\Qr$ is the largest parabolic subgroup defining the cone $\C_t(\Qr)$.
\end{itemize}
In particular, the map $\Qr \mapsto \C_t(\Qr)$ gives a one-to-one correspondence between  $t$-relevant parabolic subgroups containing $\Sr$ and cones in the fan $\mathcal{F}_t$.
\end{Prop}

(See \cite[Remark 3.25]{RTW1}.)

\vspace{0.2cm} \noindent \textbf{(1.5)} For any parabolic subgroup $\Qr$ of $\G$ containing $\Sr$, the cone $\C_t(\Qr)$ admits the following root-theoretic description. Let $\rP$ be a parabolic subgroup of type $t$ osculatory with $\Qr$. We have $$\C_t(\rP) = \{z \in \Lambda(\Sr) \ ; \ \langle \alpha, z \rangle \leqslant 1 \textrm{ for all } \alpha \in \Phi({\rm rad}^{\rm u}(\rP^{\rm op}),\Sr)\},$$ and $\C_t(\Qr)$ is the face of $\C_t(\rP)$ cut out by the linear subspace $$\langle \C_t(\Qr)\rangle = \{z \in \Lambda(\Sr) \ ; \ \langle \alpha, z \rangle = 1 \ \textrm{ for all } \alpha \in \Phi(\Lr_{\Qr},\Sr) \cap \Phi({\rm rad}^{\rm u}(\rP^{\rm op}),\Sr)\},$$ where ${\rm rad}^{\rm u}(\cdot)$ stands for the unipotent radical and $\Lr_{\Qr}$ denotes the Levi subgroup of $\Qr$ associated with $\Sr$ (\cite[Lemma 3.15]{RTW1}).

\vspace{0.1cm}
One deduces the following root-theoretical characterization of $t$-relevancy. Let $\Sr$ be a maximal split torus of $\G$. We fix a minimal parabolic subgroup $\rP_0$ of $\G$ containing $\Sr$ and write $\Delta$ for the corresponding basis of $\Phi(\G,\Sr)$, which we identify with the set of vertices in the Dynkin diagram of $\G$. The map $$\left\{\begin{array}{c}\textrm{parabolic subgroups of } \G \\ \textrm{containing } \Sr \end{array}\right\} \rightarrow \{\textrm{subsets of } \Delta \}, \ \ \  \Qr \mapsto \Y_{\Qr} = \Delta \cap \Phi(\Lr_{\Qr},\Sr)$$ is a bijection.

\begin{Prop} \label{prop.relevant.root} Let $\Qr$ be a parabolic subgroup of $\G$. We denote by $\Y_t$ the subset of $\Delta$ associated with the parabolic subgroup of type $t$ containing $\rP_0$ and let $\widetilde{\Y_{\Qr}}$ denote the union of the connected components of $\Y_{\Qr}$ meeting $\Delta - \Y_t$.

\begin{itemize}
\item[(i)] The parabolic subgroup $\Qr$ is $t$-relevant if and only if for any root $\alpha \in \Delta$, we have $$(\alpha \in \Y_t \ \textrm{ and } \ \alpha \perp \widetilde{\Y_{\Qr}}) \Longrightarrow \alpha \in \Y_{\Qr}.$$
\item[(ii)] More generally, the smallest $t$-relevant parabolic subgroup of $\G$ containing $\Qr$ is associated with the subset of $\Delta$ obtained by adjoining to $\Y_{\Qr}$ all roots in $\Y_t$ which are orthogonal to each connected component of $\Y_{\Qr}$ meeting $\Delta - \Y_t$.
\item[(iii)] The linear subspace of $\Lambda(\Sr)$ spanned by the cone $\C_t(\Qr)$ is the vanishing locus of $\widetilde{\Y_{\Qr}}$: $$\langle \C_t(\Qr) \rangle = \{z \in \Lambda(\Sr) \ ; \ \langle \alpha, z \rangle = 1 \ \textrm{ for all } \alpha \in \widetilde{\Y_{\Qr}}\}.$$
\end{itemize}
\end{Prop}

(For assertions (i) and (ii), see \cite[Proposition 3.24]{RTW1} and \cite[Remark 3.25, 2]{RTW1}. Assertion (iii) follows from \cite[Proposition 3.22]{RTW1} and \cite[Remark 3.25, 2]{RTW1}.)

Here, orthogonality is understood with respect to a scalar product on $\X^{*}(\Sr) \otimes_{\mathbb{Z}} \mathbb{R}$ invariant under the Weyl group of $\Phi(\G,\Sr)$.

\vspace{0.1cm} \begin{Rk} \label{rk.cones} Given a maximal split torus $\Sr$ and a parabolic subgroup $\Qr$ containing $\Sr$, we have the following inclusions of cones $$\mathfrak{C}(\Qr) = \C_{\varnothing}(\Qr) \subset \C_t(\Qr) \subset \C_{t(\Qr)}(\Qr)$$ for any $k$-rational type $t$.
Up to a central isogeny, we can write $\Lr_{\Qr}$ as the product $\Lr' \times \Lr''$ of two reductive groups such that $t$ has non-degenerate restriction to $\Lr'$ and trivial restriction to $\Lr''$. This amounts to decomposing $\Phi(\Lr_{\Qr},\Sr)$ as the union of two disjoint closed subsets $\Phi(\Lr',\Sr)$ and $\Phi(\Lr'',\Sr)$, with $$\Phi(\Lr',\Sr) = \langle \widetilde{\Y_{\Qr}} \rangle \cap \Phi(\G,\Sr)$$ if we use the notation introduced in the preceding proposition. It follows from the latter that the cone $\C_t(\Qr)$ is the intersection of $\C_{t(\Qr)}(\Qr)$ with the linear subspace of $\Lambda(\Sr)$ cut out by all roots in $\Phi(\Lr',\Sr)$.
\end{Rk}

\vspace{0.2cm} \noindent \textbf{(1.6)} Finally, we describe the stabilizer of a point of $\overline{\mathcal{B}}_t(\G,k)$.

\begin{Thm} \label{thm.stabilizer} Let $x$ be a point in $\overline{\mathcal{B}}_t(\G,k)$ and let $\Qr$ denote the $t$-relevant parabolic subgroup of $\G$ corresponding to the stratum containing $x$.

\begin{enumerate}
\item There exists a largest smooth and connected closed subgroup $\R_t(\Qr)$ of $\G$ satisfying the following conditions: \begin{itemize} \item[$\bullet$] $\R_t(\Qr)$ is a normal subgroup of $\Qr$ and contains ${\rad}(\Qr)$; \item[$\bullet$] for any non-Archimedean extension $\K/k$, the subgroup $\R_t(\Qr)(\K)$ of $\G(\K)$ acts trivially on the stratum $\mathcal{B}(\Qr_{\rm ss},\K)$. \end{itemize}
\item The canonical projection $\Qr_{\rm ss} \rightarrow \Qr/\R_t(\Qr)$ identifies the buildings $\mathcal{B}_t(\Qr_{\rm ss},k)$ and $\mathcal{B}(\Qr/\R_t(\Qr),k)$.
\item There exists a unique geometrically reduced $k$-analytic subgroup ${\rm Stab}_{\G}^{\ t}(x)$ of $\G^{\rm an}$ such that, for any non-Archimedean extension $\K/k$, the group ${\rm Stab}_{\G}^t(x)(\K)$ is the subgroup of $\G(\K)$ fixing $x$ in $\overline{\mathcal{B}}_t(\G,\K)$.
\item We have $\R_t(\Qr)^{\rm an} \subset {\rm Stab}_{\G}^t(x)^{\rm an} \subset \Qr^{\rm an}$ and the canonical isomorphism $\Qr^{\rm an}/\R_t(\Qr)^{\rm an} \cong (\Qr/\R_t(\Qr))^{\rm an}$ identifies the quotient group ${\rm Stab}_{\G}^{\ t}(x)/\R_t(\Qr)^{\rm an}$ with the affinoid subgroup $(\Qr/\R_t(\Qr))_x$ of $(\Qr/\R_t(\Qr))^{\rm an}$ attached in (1.1) to the point $x$ of $\mathcal{B}_t(\Qr_{\rm ss},k) = \mathcal{B}(\Qr/\R_t(\Qr),k)$.
\end{enumerate}
\end{Thm}

(See \cite[Proposition 4.7 and Theorem 4.11]{RTW1}.)

\vspace{0.1cm} \begin{Rk} \label{rk.bounded} If $\Qr$ is a proper $t$-relevant parabolic subgroup of $\G$, then ${\rm rad}(\Qr)(k)$ is an unbounded subgroup of $\G(k)$. Since ${\rm rad}(\Qr) \subset \R_t(\Qr) \subset {\rm Stab}_{\G}^{\ t}(x)$ for any $x \in \mathcal{B}_t(\Qr_{\rm ss},k)$, it follows that any point lying in the boundary $\overline{\mathcal{B}}_t(\G,k) - \mathcal{B}_t(\G,k)$ has an unbounded stabilizer in $\G(k)$. If the type $t$ is non-degenerate, the converse assertion is true.
\end{Rk}

%\vspace{0.1cm} \begin{Ex}
%\end{Ex}

\vspace{0.1cm} We can give a more precise description of the subgroup ${\rm Stab}_{\G}^{\ t}(x)(k)$ of $\G(k)$ stabilizing a point $x$ of $\overline{\mathcal{B}}_t(\G,k)$. Let us fix some notation. We pick a maximal split torus $\Sr$ of $\G$ whose compactified apartment contains $x$ and set $\N = {\rm Norm}_{\G}(\Sr)$. Let $\Qr$ denote the $t$-relevant parabolic subgroup of $\G$ attached to the stratum containing $x$ and write $\Lr$ for the Levi factor of $\Qr$ with respect to $\Sr$. We set $\Lr'' = \R_t(\Qr) \cap \Lr$ and let $\Lr'$ denote the semisimple subgroup of $\Lr$ generated by the isotropic almost simple components of $\Lr$ on which $t$ is non-trivial. Both the product morphism $\Lr' \times \Lr'' \rightarrow \Lr$ and the morphism $\Lr' \rightarrow \Qr/\R_t(\Qr)$ are central isogenies. We introduce also the split tori $\Sr' = (\Lr' \cap \Sr)^{\circ}$ and $\Sr'' = (\Lr'' \cap \Sr)^{\circ}$.

Let $\N(k)_x$ denote the stabilizer of $x$ in the $\N(k)$-action on $\overline{\A}_t(\Sr,k)$. Finally, we fix a special point in $\A(\Sr,k)$ and we recall that, for each root $\alpha \in \Phi(\G,\Sr)$, Bruhat-Tits theory endows the group $\U_{\alpha}(k)$ with a decreasing filtration $\{\U_{\alpha}(k)_r\}_{r \in [-\infty,\infty]}$.

\begin{Thm} \label{thm.stabilizer.k-points}  Let $x$ be a point in $\overline{\mathcal{B}}_t(\Qr,k)$ and let $\Qr$ denote the $t$-relevant parabolic subgroup of $\G$ attached to the stratum containing $x$.

The group ${\rm Stab}_{\G}^{\ t}(x)(k)$ is Zariski dense in $\Qr$ and is generated by the following subgroups of $\G(k)$:
\begin{itemize} \item $\N(k)_x$; \item all $\U_{\alpha}(k)$ with $\alpha \in \Phi({\rm rad}^{\rm u}(\Qr),\Sr)$; \item all $\U_{\alpha}(k)$ with $\alpha \in \Phi(\Lr'', \Sr'')$; \item all $\U_{\alpha}(k)_{-\log \alpha(x)}$ with $\alpha \in \Phi(\Lr',\Sr')$
\end{itemize}
\end{Thm}

(See \cite[Theorem 4.14]{RTW1}.)

An easy consequence of this description of stabilizers is the following generalization of well-known properties of Bruhat-Tits buildings.

\begin{Thm} \label{thm.bruhat.decomposition} 1. Let $\Sr$ be a maximal split torus of $\G$ and set $\N = {\rm Norm}_{\G}(\Sr)$. The compactified building $\overline{\mathcal{B}}_t(\G,k)$ is the topological quotient of $\G(k) \times \overline{\A}_t(\Sr,k)$ by the following equivalence relation: $$(g,x) \sim (h,y) \Longleftrightarrow \left(\exists n \in \N(k), \ y = n \cdot x \ \textrm{ and } g^{-1}hn \in {\rm Stab}_{\G}^{\ t}(x)(k)\right).$$

\noindent 2. Let $x$ and $y$ be two points in $\overline{\mathcal{B}}_t(\G,k)$.
\begin{itemize}
\item[(i)] There exists a maximal split torus $\Sr$ in $\G$ such that $x$ and $y$ lie in $\overline{\A}_t(\Sr,k)$.
\item[(ii)] The group ${\rm Stab}_\G^t(x)(k)$ acts transitively on the compactified apartments containing $x$.
\item[(iii)] We have the following \emph{mixed Bruhat decomposition}: $$\G(k) = {\rm Stab}_{\G}^{\ t}(x)(k) \N(k) {\rm Stab}_{\G}^{\ t}(y)(k).$$
\end{itemize}
\end{Thm}

(See \cite[Corollary 4.15 and Theorem 4.20]{RTW1}.)

\vspace{0.2cm} \noindent \textbf{(1.7)} Many statements listed above are proved by using an explicit formula for the map $\vartheta_t$ when $\G$ is \emph{split}.

\vspace{0.1cm} Let $\rP$ be a parabolic subgroup of $\G$ of type $t$ and pick a maximal split torus $\Sr$ of $\G$ contained in $\rP$. The morphism $${\rm rad}^{\rm u}(\rP^{\rm op}) \rightarrow {\rm Par}(\G), \ g \mapsto g \rP g^{-1}$$ is an isomorphism onto an open subscheme of ${\rm Par}(\G)$ which we denote by $\Omega(\Sr,\rP)$. Let $\Phi(\G,\Sr)$ be the set of roots of $\G$ with respect to $\Sr$. Since $\G$ is split, the choice of a special point $o$ in $\A(\Sr,k)$ determines a $k^{\circ}$-Chevalley group $\mathcal{G}$ with generic fibre $\G$. Any
Chevalley basis in $\mathrm{Lie}(\mathcal{G})(k^{\circ})$ leads to an isomorphism of
$\radu(\rP^{\rm op})$ with the affine space
$$\prod_{\alpha \in \Psi} \U_{\alpha} \simeq \prod_{\alpha \in \Psi} \mathbb{A}^{1}_k,$$
where $\Psi = \Phi(\radu(\rP^{\rm op}),\Sr) = -\Phi(\radu(\rP),\Sr)$.

\begin{Prop} \label{prop.explicit}
We assume that the group $\G$ is split and we use the notation introduced above.

\begin{itemize}
\item[(i)] The map $\vartheta_t$ sends the point $o$ to the point of
$\Omega(\Sr,\rP)^{{\rm an}}$ corresponding to the multiplicative (semi)norm
$$k\left[(\X_{\alpha})_{\alpha \in \Psi}\right] \rightarrow  \mathbb{R}_{\geqslant 0}, \ \
\sum_{\nu \in \mathbb{N}^{\Psi}} a_{\nu} \X^{\nu} \mapsto \max_{\nu} |a_{\nu}|.$$
\item[(ii)] Using the point $o$ to identify the apartment $\A(\Sr,k)$ with the vector space $\Lambda(\Sr) = {\rm Hom}_{\mathbf{Ab}}(\X^*(\Sr),\mathbb{R}_{>0})$,
the map $\Lambda(\Sr) \rightarrow {\rm Par}(\G)^{{\rm an}}$ induced by $\vartheta_t$
associates with an element $u$ of $\Lambda(\Sr)$
the point of $\Omega(\Sr,\rP)^{{\rm an}}$ corresponding to the multiplicative
seminorm
$$k\left[(\X_{\alpha})_{\alpha \in \Psi}\right] \rightarrow  \mathbb{R}_{\geqslant 0}, \ \
\sum_{\nu \in \mathbb{N}^{\Psi}} a_{\nu} \X^{\nu} \mapsto \max_{\nu}
|a_{\nu}|\prod_{\alpha \in \Psi} \langle u, \alpha \rangle^{\nu(\alpha)}.$$
\end{itemize}
\end{Prop}

(See \cite[Proposition 2.18]{RTW1}.)

\section{Comparison with gluings}

\label{ss - comparison gluings}

We show in this section that the compactifications defined in~\cite{Wer2} occur among the Berkovich compactifications. Let $k$ be a non-Archimedean local field and let $\G$ be a connected semisimple $k$-group. We consider a faithful and geometrically irreducible linear representation $\rho: \G \rightarrow \mathrm{GL}_\V$ of $\G$. In~\cite{Wer2}, a compactification $\overline{\mathcal{B}} (\G,k)_\rho$ of the Bruhat-Tits building is constructed using the combinatorics of weights for $\rho$. It only depends on the Weyl chamber face position of the highest weight of the representation.

\vspace{0.2cm}
\noindent \textbf{(2.1)} We fix a maximal split torus $\Sr$ in $\G$ and denote by $\Phi = \Phi(\G, \Sr)$ the root system of $\G$ with respect to $\Sr$. We denote by $\W$ the Weyl group of $\Phi$ and choose a $\W$-invariant scalar product $( \cdot | \cdot )$ on the character group $\X^*(\Sr)$ of $\Sr$, which we use to embed $\X^*(\Sr)$ in the vector space $\Lambda(\Sr) = \mathrm{Hom}_{\mathbf{Ab}}(\X^*(\Sr),\mathbb{R}_{>0})$ via the map $$\X^*(\Sr) \rightarrow \Lambda(\Sr),\ \ \ \chi \mapsto e^{(\chi|\cdot)}.$$

Let $\Delta$ be a basis of $\Phi$. For every subset $\Y$ of $\Delta$, we denote as in \cite{Wer2} by $\rP_\Y^\Delta$ the standard parabolic subgroup associated with $\Y$; in particular, $\rP_{\varnothing}^{\Delta}$ is the minimal parabolic subgroup of $\G$ containing $\Sr$ and corresponding to $\Delta$. The weights with respect to the action of $\Sr$ on $\V$ are called the $k$-weights of $\rho$. If $\T$ is a maximal torus containing $\Sr$ and if $k'/k$ is a finite extension splitting $\T$, then we have a natural projection $$\X^*(\T \otimes_k k') \rightarrow \X^{*}(\Sr \otimes_k k') = \X^*(\Sr)$$ and there exists a basis $\Delta'$ of $\Phi(\G \otimes_k k', \T \otimes_k k')$ lifting $\Delta$. With the basis $\Delta'$ is associated a well-defined character of $\T \otimes_k k'$, the \emph{highest weight} $\lambda_0(\Delta')$, whose restriction to $\Sr$ does not depend on any choice made for $\T$, $k'$ and $\Delta'$. This character of $\Sr$, denoted $\lambda_0(\Delta)$, is called the highest $k$-weight of $\rho$ with respect to $\Delta$; it defines an element in $\Lambda(\Sr)$ lying in the Weyl cone $\mathfrak{C}(\rP^\Delta_\varnothing)$. Setting $$\Z = \{\alpha \in \Delta \ ; \ (\lambda_0(\Delta) | \alpha) = 0\},$$ the linear subspace $\{\alpha = 1 \ ; \ \alpha \in \Z\}$ cuts out the only face of $\mathfrak{C}(\rP^{\Delta}_{\varnothing})$ whose interior contains $\lambda_0 (\Delta)$.

\vspace{0.1cm}

The purpose of this paragraph is to prove the following theorem.

\vspace{0.1cm}

\begin{Thm}
\label{thm.comparison}
Let $\tau$ denote the type of the parabolic subgroup $\rP_{\Z}^{\Delta}$.
The compactified buildings $\overline{\mathcal{B}}(\G,k)_\rho$ and $\overline{\mathcal{B}}_\tau(\G,k)$ are canonically isomorphic, and $\tau$ is the only $k$-rational type satisfying this condition.
\end{Thm}

\vspace{0.1cm} \begin{Rk} \label{rk.type}
Up to conjugacy, it is clear that the parabolic subgroup $\rP_{\Z}^{\Delta}$ does not depend on the choice of $\Sr$ and $\Delta$. Therefore, the $k$-rational type $t(\rP_{\Z}^{\Delta})$ is canonically associated with the absolutely irreducible representation $\rho$. One the other hand, the theory of highest weights of irreducible linear representations of split reductive groups singles out naturally a well-defined type $t(\rho)$ of parabolic subgroups of $\G$, maybe non-$k$-rational: the connected component ${\rm Par}_{t(\rho)}(\G)$ of ${\rm Par}(\G)$ is characterized by the condition that, for any finite extension $k'/k$ splitting $\G$, this component contains all the maximal parabolic subgroups of $\G$ stabilizing a highest line in $\V \otimes_k k'$ (see paragraph 4.1). We conclude this article by establishing that $t(\rP_{\Z}^{\Delta})$ is the unique $k$-rational type defining the same compactification of $\mathcal{B}(\G,k)$ as the type $t(\rho)$ (cf. \cite[Appendix C]{RTW1}); equivalently, the compactification $\overline{\mathcal{B}}(\G,k)_{\rho}$ defined in \cite{Wer2} is canonically isomorphic to the Berkovich compactification $\overline{\mathcal{B}}_{t(\rho)}(\G,k)$ (see Proposition \ref{prop.types}).
\end{Rk}

\vspace{0.3cm}

Before proving this theorem, we can derive at once a comparison with the group-theoretic compactification~\cite{GuiRem}.

\vspace{0.1cm}

\begin{Cor}
\label{cor.comparisonGuiRem}
Let $t_{\rm min}$ be the type of a minimal parabolic subgroup of $\G$.
We denote by $\mathcal{V}_{\mathcal{B}(\G,k)}$ the set of vertices in the Bruhat-Tits building ${\mathcal{B}}(\G,k)$.
Then the closure of $\mathcal{V}_{\mathcal{B}(\G,k)}$ in the maximal Berkovich compactification $\overline{\mathcal{B}}_{t_{\rm min}}(\G,k)$ is $\G(k)$-equivariantly homeomorphic to the group-theoretic compactification of $\mathcal{V}_{{\mathcal{B}}(\G,k)}$.
\end{Cor}

\vspace{0.1cm}
\noindent
\emph{\textbf{Proof of corollary}}. By~\cite[Theorem 20]{GuiRem}, the group-theoretic compactification of $\mathcal{V}_{{\mathcal{B}}(\G,k)}$ is $\G(k)$-equivariantly homeomorphic to the closure of $\mathcal{V}_{\mathcal{B}(\G,k)}$ in the polyhedral compactification of ${\mathcal{B}}(\G,k)$ defined by E. Landvogt.
By~\cite{Wer2}, we know that the latter compactification is $\G(k)$-equivariantly homeomorphic to $\overline{\mathcal{B}}(\G,k)_\rho$ where $\rho$ is any weight lying in the interior of some Weyl chamber, i.e., such that $\Z = \varnothing$ with the notation above.
Our claim follows from Theorem \ref{thm.comparison}.
\hfill $\Box$

\vspace{0.3cm}

Recall that every $k$-weight of $\rho$ is of the form $\lambda_0(\Delta) - \sum_{\alpha \in \Delta} n_\alpha \alpha$ for certain non-negative integers $n_\alpha$. We denote by $[\lambda_0(\Delta) - \lambda] = \{\alpha \in \Delta \ ; \ n_{\alpha} >0\}$ the support of $\lambda_0(\Delta) - \lambda$.
In~\cite[Definition 1.1]{Wer2}, a subset $\Y \subset \Delta$ is called \emph{admissible}, if the set $\Y \sqcup \{\lambda_0(\Delta)\}$ is connected in the following sense: the graph with vertex set $\Y \cup \lambda_0(\Delta)$ and edges between all $\alpha$ and $\beta$ such that $(\alpha | \beta) \neq 0$  is connected.

\vspace{0.1cm}

The following lemma is well-known, at least in characteristic 0~\cite[12.16]{BoTi}.
It is a link between the abstract root-theoretic definition of admissibility, and its interpretation in terms of representations.

\begin{Lemma}
\label{lemma.admissibility}
A set $\Y \subset \Delta$ is admissible if and only if there exists a $k$-weight $\mu$ whose support $[\lambda_0( \Delta) - \mu]$ is equal to $\Y$.
\end{Lemma}

\vspace{0.1cm}
\noindent
\emph{\textbf{Proof}}.
For the sake of completeness, we show that this statement holds whatever the characteristic of $k$ is.
In order to be short, we freely use the notation of~\cite[\S 24.B]{Borel}, which sums up the basic results of representation theory of reductive groups over arbitrary fields.
In particular, given $\G$ as above, we denote by ${\rm E}^\lambda$ the unique Weyl $\G$-module of highest weight $\lambda$ and by ${\rm F}^\lambda$ its unique irreducible submodule (which in turn determines ${\rm E}^\lambda$); in characteristic 0, we have ${\rm F}^\lambda={\rm E}^\lambda$.
Note that in the setting of this section, the $\G$-module $\V$ is isomorphic to some ${\rm F}^\lambda$ and remains irreducible after extension of the ground field to the algebraic closure of $k$.

\vspace{0.1cm}

Let us first assume that $\Y$ is the support of some weight.
Since the irreducible module ${\rm F}^\lambda$ is a submodule of the Weyl $\G$-module ${\rm E}^\lambda$, we deduce that $\Y$ is the support of some weight for ${\rm E}^\lambda$.
Moreover the Weyl module ${\rm E}^\lambda$ has the same character formula as the irreducible module of highest weight $\lambda$ in characteristic 0, so the connectedness of the graph under consideration comes from the result in this case~\cite[12.16]{BoTi}.
Note that we use the classification of semisimple groups in order to find a group over a field of characteristic 0 having the same representations as $\G$.

\vspace{0.1cm}

Conversely, let us assume that the graph $\Y \sqcup \{\lambda_0(\Delta)\}$ is connected.
Recall that the set of weights is stable under the spherical Weyl group.
We investigate first the case when $\Y$ is connected. We write $\Y = \{ \beta_1, \beta_2, \ldots , \beta_m \}$ in such a way that
$\beta_1$ is connected to $ \lambda_0(\Delta)$ (i.e., $(\lambda_0(\Delta) \mid \beta_1) \neq 0$) and that for any $i \leqslant m$ there exists $j<i$ such that $\beta_i$ is connected to $\beta_j$ (i.e., $(\beta_i \mid \beta_j) \neq 0$).
Then it is easy to show by a finite induction on $l \leqslant m$, that the support of the weight $r_{\beta_l}r_{\beta_{l-1}}...r_{\beta_1}(\lambda_0(\Delta))$ is equal to $\{\beta_1,\beta_2, ... \ , \beta_l \}$.
Indeed, for $l=1$ this is clear since $r_{\beta_1}(\lambda_0(\Delta)) = \lambda_0(\Delta) - 2\frac{(\lambda_0(\Delta)\mid\beta_1)}{(\beta_1|\beta_1)} \beta_1$; and to pass from one step to the next one, we argue as follows.
First, we have:
$$r_{\beta_l}r_{\beta_{l-1}}...r_{\beta_1}(\lambda_0(\Delta)) = r_{\beta_l}\left(\lambda_0(\Delta) - \sum_{i=1}^{l-1} c_i \beta_i\right),$$
with each $c_i>0$ by induction hypothesis. This gives:
$$
r_{\beta_l}r_{\beta_{l-1}}...r_{\beta_1}(\lambda_0(\Delta)) = \lambda_0(\Delta) - \sum_{i=1}^{l-1} c_i \beta_i - 2\left(\frac{(\lambda_0(\Delta)\mid \beta_l)}{(\beta_l|\beta_l)} - \sum_{i=1}^{l-1} c_i \frac{(\beta_i \mid \beta_l)}{(\beta_l|\beta_l)}\right)\beta_l,
$$
which implies our claim by the numbering of the $\beta_i$'s and the fact that $\lambda_0(\Delta)$ is dominant.

In the general case, we use a numbering $\Y_1, \Y_2, ... \, \Y_s$ of the connected components of  $\Y$.
The previous argument shows that there is a weight, say $\mu$, with support equal to $\Y_1$.
Then we note that for each $\alpha \in \Y_1$ and each $\beta \in \Y_2$ we have $r_\beta(\alpha)=\alpha$.
This allows us to apply the previous argument, replacing $\lambda_0(\Delta)$ by $\mu$ and $\Y$ by $\Y_2$. Our claim follows by induction on the number of connected components of $\Y$.
\hfill $\Box$

\vspace{0.2cm}
\noindent \textbf{(2.2)} For every admissible subset $\Y \subset \Delta$ we set
$$\Y^{*} = \{\alpha \in \Delta \ ; \ (\alpha |\lambda_0(\Delta)) = 0 \ \textrm{and} \ (\alpha|\Y) = 0\}$$
and let $\C_{\Y}^{\Delta}$ denote the cone in $\Lambda(\Sr)$ defined by the following conditions
$$\left\{\begin{array}{ll} \alpha = 1, & \mbox{ for all } \alpha \in \Y \\
\lambda_0(\Delta)-\lambda \geqslant  1, & \mbox{ for all $k$-weights } \lambda \mbox{ such that } [\lambda_0(\Delta) - \lambda] \not\subset \Y.\end{array}\right.$$

Identifying the additive and multiplicative duals of $\X^*(\Sr)$ via the map $\mathbb{R} \rightarrow \mathbb{R}_{>0}, \ x \mapsto e^{x}$, the cone $\C_{\Y}^{\Delta} \subset \Lambda(\Sr)$ is the closure of the subset $\F_{\Y}^{\Delta}$ of $\V(\Sr) = {\rm Hom}_{\mathbf{Ab}}(\X^{*}(\Sr), \mathbb{R})$ defined in~\cite[section 2]{Wer2}. It is shown in [loc. cit.] that $\V(\Sr)$ is the disjoint union of the subsets $\F_{\Y}^{\Delta}$, where $\Y$ runs over the set of admissible subsets of $\Delta$.

\begin{Lemma}
\label{lemma.comparison.fans1}
Recall that $\Z = \varnothing^*$ and let $\tau$ denote the type of the parabolic subgroup $\rP_{\Z}^{\Delta}$.
\begin{itemize}
\item[(i)] A subset $\Y$ of $\Delta$ is admissible if and only if each of its connected components meets $\Delta - \Z$.
\item[(ii)] For any admissible subset $\Y$ of $\Delta$, we have $$\C_{\Y}^{\Delta} = \C_\tau(\rP_{\Y}^{\Delta}).$$
\item[(iii)] The correspondence $\Y \mapsto \rP_{\Y \cup \Y^{*}}^{\Delta}$ is a bijection between admissible subsets of $\Delta$ and $\tau$-relevant parabolic subgroups containing $\rP_{\varnothing}^{\Delta}$.
\end{itemize}
\end{Lemma}

\vspace{0.1cm}
\noindent
\emph{\textbf{Proof}}.
(i) This assertion is clear, since $\Y \cup \{\lambda_{0}(\Delta)\}$ is connected if and only if each connected component of $\Y$ contains a root $\alpha \in \Delta$ with $(\alpha | \lambda_0(\Delta)) \neq 0$, i.e., a root in $\Delta - \Z$.

\vspace{0.1cm} (ii) Let $\Y$ be an admissible subset of $\Delta$. It follows from (i) and from Proposition \ref{prop.relevant.root} (iii) that the linear space $\{\alpha =1; \ \alpha \in \Y\}$ cuts out a face of the cone $\C_\tau(\rP_{\varnothing}^{\Delta})$, namely the cone $\C_\tau(\rP_{\Y}^{\Delta})$. Since this subspace cuts out the face $\C_{\Y}^{\Delta}$ of $\C_{\varnothing}^{\Delta}$, it suffices to check that the cones $\C_\tau(\rP_{\varnothing}^{\Delta})$ and $\C_{\varnothing}^{\Delta}$ coincide.

Let $\Delta'$ be another basis of the root system $\Phi$. If $\lambda_0(\Delta') = \lambda_0(\Delta)$, then every $x$ in the Weyl cone $\mathfrak{C}(\rP_{\varnothing}^{\Delta'})$ satisfies $(\lambda_0(\Delta) - \lambda)(x) \geqslant 1$ for all $k$-weights $\lambda$, hence $\mathfrak{C}(\rP_{\varnothing}^{\Delta'})$ is contained in $\C_\varnothing^\Delta$. On the other hand, every point in the interior of $\C_\varnothing^\Delta$ is contained in  the Weyl cone $\mathfrak{C}(\rP_{\varnothing}^{\Delta'})$ for some basis $\Delta'$. By~\cite[Proposition 4.4 and Lemma 2.1]{Wer2}, this implies $\lambda_0(\Delta') = \lambda_0(\Delta)$. Hence $\C_\varnothing^\Delta$ is equal to the union of all Weyl cones
$\mathfrak{C}(\rP_{\varnothing}^{\Delta'})$ with $\lambda_0(\Delta) = \lambda_0(\Delta')$. By definition, the cone $\C_\tau(\rP_\varnothing^\Delta)$ is the union of all $\mathfrak{C}(\rP_{\varnothing}^{\Delta'})$ such that the minimal parabolic subgroup $\rP_\varnothing^{\Delta'}$ is contained in $\rP_{\Z}^\Delta$. Therefore, it remains to check that $\lambda_0(\Delta) = \lambda_0(\Delta')$, if and only if $\rP_\varnothing^{\Delta'}$ is contained in $\rP_{\Z}^\Delta$.

\vspace{0.1cm}
Let $n$ be an element of $\mathrm{Norm}_{\G}(\Sr)(k)$ satisfying $n \rP_\varnothing^\Delta n^{-1} = \rP_\varnothing^{\Delta'}$, and let $w$ be its image in the Weyl group $\W$ of $\Phi$. Then $w (\Delta) = \Delta'$, hence $w (\lambda_0(\Delta)) = \lambda_0(\Delta')$. Besides, we have $n \rP_{\Z}^\Delta n^{-1}  = \rP_{w(\Z)}^{\Delta'}$.

Assume that $\lambda_0(\Delta) = \lambda_0(\Delta')$. Then $w$ fixes $\lambda_0(\Delta)$, which implies that
$w(\Z) = \Z$ since the scalar product on $\X^*(\Sr)$ is $\W$-invariant. Besides, for every $\alpha \in \Delta - \Z$ there exists a $k$-weight $\lambda$ such that $[\lambda_0(\Delta) - \lambda] = \{\alpha\}$ for $\{ \alpha\}$ is an admissible subset of $\Delta$. Since $w(\lambda)$ is a weight and $w(\lambda_0(\Delta)) = \lambda_0(\Delta)$, we deduce that $w(\alpha)$ is a positive root for $\Delta$. Hence $\rP^{\Delta}_{\Z}$ contains $\rP_\varnothing^{w(\Delta)} = \rP_\varnothing^{\Delta'}$.

Now assume that $\rP_\varnothing^{\Delta'} = n \rP_\varnothing^\Delta n^{-1}$ is contained in $\rP_{\Z}^\Delta$. Then $n$ is contained
in $\rP_{\Z}^\Delta$, which implies that $w$ is in the Weyl group of the parabolic $\rP_{\Z}^\Delta$. Hence $w$ is a product of reflections corresponding to roots in $\Z$. Since roots in $\Z$ are perpendicular to $\lambda_0(\Delta)$, the corresponding reflections leave $\lambda_0(\Delta)$ invariant and therefore $\lambda_0(\Delta')= w(\lambda_0(\Delta)) = \lambda_0(\Delta)$.

\vspace{0.1cm} (iii) Let $\Y$ be an admissible subset of $\Delta$. By Proposition \ref{prop.relevant.root} (ii), the smallest $\tau$-relevant parabolic subgroup containing $\rP_{\Y}^{\Delta}$ is $\rP_{\Y'}^{\Delta}$, where $\Y'$ is obtained by adjoining to $\Y$ all roots in $\Z$ which are perpendicular to each connected component of $\Y$ meeting $\Delta - \Z$, hence to $\Y$ by (i). It follows that $\Y' = \Y \cup \Y^*$. Conversely, if $\rP_{\Z}^{\Delta}$ is a $\tau$-relevant parabolic subgroup, then $\C_\tau(\rP^{\Delta}_{\Z}) = \C^{\Delta}_{\Y}$ for some admissible subset $\Y$ and $\C_{\Y}^{\Delta} = \C(\rP_{\Y}^{\Delta})$ by (i). It follows from what we have just said that $\rP^{\Delta}_{\Y \cup \Y^*}$ is the smallest $\tau$-relevant parabolic subgroup containing $\rP_{\Y}^{\Delta}$, hence $\C_\tau(\rP^{\Delta}_{\Z}) = \C^{\Delta}_{\Y} = \C_\tau(\rP^{\Delta}_{\Y \cup \Y^*})$ and therefore $\Z = \Y \cup \Y^*$.
\hfill $\Box$

\vspace{0.3cm}
Thus, the fan consisting of all polyhedral cones $\C^\Delta_\Y$ coincides with the fan $\mathcal{F}_\tau$ defined in (1.4). Note that the type $\tau$ is non-degenerate since the representation $\rho$ is faithful. Relying on \cite[Proposition B.3]{RTW1}, it is not hard to check that the identity map of the apartment $\A(\Sr,k)$ extends to a homeomorphism $j$ between the compactification $\overline{\A}_\tau(\Sr,k) = \overline{\A(\Sr,k)}^{\mathcal{F}_\tau}$ introduced in Definition \ref{def.compactifications} and the compactification $\overline{\A(\Sr,k)}_{\rho}$ of $\A(\Sr,k)$ defined from a different viewpoint in~\cite[Sect. 2]{Wer2} (where it is simply denoted $\overline{\A}$). This homeomorphism is compatible with the action of the group $\mathrm{Norm}_{\G}(\Sr)(k)$ on each space since this action is in both cases the unique continuous extension of the standard action of $\mathrm{Norm}_{\G}(\Sr)(k)$ on $\A(\Sr,k)$.

\vspace{0.2cm} \noindent \textbf{(2.3)} Seen as a function $\Lambda(\Sr) \rightarrow \mathbb{R}_{>0}$, each root $\alpha \in \Phi$ has a continuous extension $\widetilde{\alpha} : \overline{\C} \rightarrow [0,\infty]$ for every cone $\C$ in the fan $\mathcal{F}_\tau$ over which either $\alpha \leqslant 1$ or $\alpha \geqslant 1$; this is obvious if we write $\C = {\rm Hom}_{\mathbf{Mon}}(\M,]0,1])$ and $\overline{\C} = {\rm Hom}_{\mathbf{Mon}}(\M,[0,1])$, where $\M$ is the saturated and finitely generated semigroup in $\X^*(\Sr)$ defined by $$\M = \{\alpha \in \X^*(\Sr) \ ; \ \alpha_{|\C} \leqslant 1\}.$$
If $\tau = t_{\rm min}$ is the type of a minimal parabolic subgroup, then $\mathcal{F}_\tau$ is the Weyl fan and every root $\alpha$ satisfies $\alpha_{|\C} \leqslant 1$ or $\alpha_{|\C} \geqslant 1$ for each cone $\C \in \mathcal{F}_{t_{\rm min}}$, hence extends continuously to the corresponding compactified vector space $\overline{\Lambda}(\Sr)^{\mathcal{F}_{t_{\rm min}}}$. Since we have either $\alpha < 1$, $\alpha > 1$ or $\alpha = 1$ on the interior $\F^{\circ}$ of each face $\F$ of $\C \in \mathcal{F}_{t_{\rm min}}$, the extension $\widetilde{\alpha}$ of $\alpha$ to $\overline{\C}$ satisfies
$$\left\{\begin{array}{ll} \widetilde{\alpha}_{|\C_{\F}} = 0 & \textrm{ if } \ \alpha_{|\F^{\circ}} < 1, \\ 0 < \widetilde{\alpha}_{|\C_{\F}} < \infty & \textrm{ if } \ \alpha_{|\F} = 1, \\ \widetilde{\alpha}_{|\C_{\F}} = \infty & \textrm{ if } \ \alpha_{|\F^{\circ}} > 1,\end{array} \right.$$ where $\C_{\F}$ is the stratum of $\overline{\C}$ corresponding to the face $\F$, namely the subset of $\overline{\C}$ defined by the conditions $$\left\{\begin{array}{ll} \varphi = 0, & \textrm{ for all }  \varphi \in \M \textrm{ such that } \varphi_{|\F} \neq 1, \\ \varphi > 0, & \textrm{ for all } \varphi \in \M  \textrm{ such that } \varphi_{|\F} = 1.\end{array} \right.$$
This situation is illustrated by Figure 1 below with $\G = {\rm SL}(3)$.

\vspace{0.2cm} In general, we can always extend each root $\alpha$ to a \emph{upper semicontinuous} function $\widetilde{\alpha}: \overline{\Lambda(\Sr)}^{\mathcal{F}_\tau} \rightarrow [0,\infty]$ by setting $$\widetilde{\alpha}(x) = \sup \{c \in \mathbb{R}_{>0} \ ; \ x \in \overline{\{\alpha \geqslant c\}}\}.$$ This function coincides with the continuous extension of $\alpha_{|\C}$ to $\overline{\C}$ for any cone $\C$ over which $\alpha \leqslant 1$ or $\alpha \geqslant 1$. In general, given a cone $\C$ and a face $\F$ of $\C$, the upper semicontinuous extension $\widetilde{\alpha}$ of $\alpha$ to $\overline{\C}$ satisfies
$$\left\{\begin{array}{ll} \widetilde{\alpha}_{|\C_{\F}} = 0 & \textrm{ if } \ \alpha_{|\F^{\circ}} < 1 \\ 0 < \widetilde{\alpha}_{|\C_{\F}} < \infty & \textrm{ if } \alpha_{|\F} = 1 \\ \widetilde{\alpha}_{|\C_{\F}} = \infty & \textrm{ if } \ \alpha_{|\F^{\circ}} > 1 \\ \widetilde{\alpha}_{|\C_{\F}} = \infty & \textrm{ otherwise}. \end{array} \right. $$ This follows easily from the existence of an affine function $\beta : \C \rightarrow ]0,1]$ such that $\beta_{|\F} = 1$.

This situation is illustrated by Figure 2 below, where $\G = {\rm SL}(3)$ and $\tau$ is a type of maximal proper parabolic subgroups.

\vspace{0.2cm}
With each point $x$ of $\overline{\A(\Sr,k)}_{\rho}$ is associated in~\cite{Wer2} a subgroup $\rP_x$ of $\G(k)$ defined as follows. Set $\N = \mathrm{Norm}_{\G}(\Sr)$ and recall that Bruhat-Tits theory provides us with a decreasing filtration $\{\U_{\alpha}(k)_r\}_{r \in [-\infty, \infty]}$ on each unipotent root group $\U_{\alpha}(k)$, with $\U_{\alpha}(k)_{-\log(\infty)} = \U_{\alpha}(k)_{-\infty} = \U_{\alpha}(k)$ and $\U_{\alpha}(k)_{-\log(0)} = \U_{\alpha}(k)_{\infty} = \{1\}$. Then $\rP_x$ is the subgroup of $\G(k)$ generated by $\N(k)_x = \{n \in \N(k); \ nx = x\}$ and $\U_{\alpha}(k)_{-\log \widetilde{\alpha}(x)}$ for all $\alpha \in \Phi$.

\vspace{0.1cm}
Let $\Qr$ be a $\tau$-relevant parabolic subgroup of $\G$ containing $\Sr$ and denote by $\Lr$ the Levi subgroup of $\Qr$ associated with $\Sr$. We consider the following decomposition of $\Phi$ in mutually disjoint closed subsets:
$$\Phi = \big(-\Phi(\radu(\Qr),\Sr)\big) \cup \Phi(\radu(\Qr),\Sr) \cup \Phi(\Lr',\Sr') \cup \Phi(\Lr'',\Sr''),$$ where $\Lr'$ and $\Lr''$ are the normal and connected reductive subgroups of $\Lr$ such that the natural morphisms $\Lr' \times \Lr'' \rightarrow \Lr$ and $\Lr' \rightarrow \Qr/\R_\tau(\Qr)$ are central isogenies, and where $\Sr'$ and $\Sr''$ are the connected components of $\Sr \cap \Lr'$ and $\Sr \cap \Lr''$ respectively (see the discussion before Theorem \ref{thm.stabilizer.k-points}). Equivalently, the subset $\Phi(\Lr',\Sr')$ of $\Phi(\Lr,\Sr)$ is the union of root systems $\Phi(\Hr,\Sr)$, where $\Hr$ runs over the set of quasi-simple components of $\Lr$ on which the restriction of $\tau$ is non-trivial, and $\Phi(\Lr'',\Sr'') = \Phi(\Lr,\Sr) - \Phi(\Lr',\Sr')$.

\vspace{0.1cm}
\begin{Lemma}
\label{lemma.extension.roots}
Let $x$ be a point in the stratum $\Sigma = \A(\Sr,k)/\langle \C_\tau(\Qr) \rangle$ of $\overline{\A(\Sr,k)}^{\mathcal{F}_\tau}$.
\begin{itemize}
\item[(i)] For any root $\alpha$ in $\Phi$, we have:
$$\left\{\begin{array}{ll} \widetilde{\alpha}(x) = 0 \textrm{ and } \widetilde{-\alpha}(x) = \infty & \textrm{ if }\ \alpha \in \Phi({\rm rad}^{\rm u}(\Qr^{\rm op}),\Sr); \\
\widetilde{\alpha}(x) = \infty \textrm{ and } \widetilde{-\alpha}(x) = 0 & \textrm{ if } \ \alpha \in - \Phi({\rm rad}^{\rm u}(\Qr^{\rm op}),\Sr); \\
\widetilde{\alpha}(x) = \widetilde{-\alpha}(x) = \infty & \textrm{ if } \ \alpha \in \Phi(\Lr'',\Sr''); \\
0 < \widetilde{\alpha}(x) < \infty & \textrm{ if } \ \alpha \in \Phi(\Lr',\Sr').
\end{array} \right.$$
\item[(ii)] $\rP_x = \mathrm{Stab}_{\G}^{\ t}(x)(k)$.
\end{itemize}
\end{Lemma}

\vspace{0.1cm}
\noindent
\emph{\textbf{Proof}}. (i) This assertion follows from the identities $$\Phi({\rm rad}^{\rm u}(\Qr^{\rm op}),\Sr) = \{\alpha \in \Phi \ ; \ \alpha < 1 \textrm{ on the interior of } \C_\tau(\Qr)\},$$ $$\Phi(\Lr',\Sr') = \{\alpha \in \Phi \ ; \ \alpha = 1 \textrm{ on } \C_\tau(\Qr)\}$$ and $$\Phi(\Lr'',\Sr'') = \{\alpha \in \Phi \ ; \ \alpha \textrm{ takes values } <1 \textrm{ and } >1 \textrm{ on } \C_\tau(\Qr)\}$$ (see Remark \ref{rk.cones}).

%Consider first of all a root $\alpha \in \Phi(\Qr^{\rm op},\Sr)$ and let $\rP$ denote a parabolic subgroup of type $\tau$ osculatory with $\Qr$ such that $\alpha$ belongs to $\Phi(\rP^{\rm op},\Sr)$. We have $$\C_\tau(\rP) = {\rm Hom}_{\mathbf{Mon}}(\Phi({\rm rad}^{\rm u}(\rP^{\rm op}),\Sr), ]0,1])$$ and $\C_\tau(\Qr)$ is the face cut out by the linear subspace $$\{\alpha = 1 \ ; \ \alpha \in \Phi({\rm rad}^{\rm u}(\rP^{\rm op}), \Sr) \cap \Phi(\Lr, \Sr)\} = \{alpha = 1 \ ; \ \alpha \in \Phi(\Lr',\Sr')\}$$ of $\Lambda(\Sr)$. It follows that $\Sigma \cap \overline{\C_\tau(\rP)}$ is the subset of $\overline{\C_\tau(\rP)} = {\rm Hom}_{\mathbf{Mon}}(\Phi({\rm rad}^{\rm u}(\rP^{\rm op}),\Sr), [0,1])$ defined by the conditions $$\left\{\begin{array}{ll} u(\alpha) > 0 & \textrm{ if } \alpha \in \Phi({\rm rad}^{\rm u}(\rP^{\rm op}), \Sr) \cap \Phi(\Lr, \Sr) \\ u(\alpha) = 0 &  \textrm{ if } \alpha \notin \Phi({\rm rad}^{\rm u}(\rP^{\rm op}), \Sr) \cap \Phi(\Lr, \Sr)\end{array}\right.$$

\vspace{0.1cm} (ii) This assertion follows immediately from (i) and from the explicit description of $\mathrm{Stab}_{\G}^{\ t}(x)(k)$ in Theorem \ref{thm.stabilizer.k-points} since both $\rP_x$ and $\mathrm{Stab}_{\G}^{\ t}(x)(k)$ are the subgroups of $\G(k)$ generated by $\N(k)_x$ and all $\U_{\alpha}(k)_{-\operatorname{log}\widetilde{\alpha}(x)}$, $\alpha \in \Phi$.
\hfill $\Box$

\vspace{0.1cm} The compactification $\overline{\mathcal{B}}(\G,k)_{\rho}$ defined in~\cite{Wer2} is the topological quotient of $\G(k) \times \overline{\A(\Sr,k)}_{\rho}$ by the following equivalence relation: $$(g,x) \sim (h,y) \Longleftrightarrow \big(\exists n \in \N(k), \ y=nx \ \textrm{ and} \ g^{-1}hn \in \rP_x \big).$$ It follows immediately from assertion (ii) in the previous lemma and from the first assertion of Theorem \ref{thm.bruhat.decomposition} that the canonical homeomorphism $$\G(k) \times \overline{\A(\Sr,k)}^{\mathcal{F}_\tau} \xymatrix{{} \ar@{->}[r]^{\sim} & {}} \G(k) \times \overline{\A(\Sr,k)}_{\rho}$$ induces a $\G(k)$-homeomorphism between the compactified buildings $\overline{\mathcal{B}}_\tau(\G,k)$ and $\overline{\mathcal{B}}(\G,k)_{\rho}$.

\vspace{0.1cm} Uniqueness of the $k$-rational type $\tau$ such that the compactifications $\overline{\mathcal{B}}(\G,k)_{\rho}$ and $\overline{\mathcal{B}}_\tau(\G,k)$ are isomorphic is easily checked. For any $k$-rational type $\tau'$ satisfying this condition, the compactifications $\overline{\mathcal{B}}_\tau(\G,k)$ and $\overline{\mathcal{B}}_{\tau'}(\G,k)$ are $\G(k)$-equivariantly homeomorphic. This homeomorphism identifies $0$-dimensional strata; taking stabilizers in $\G(k)$, we obtain two parabolic subgroups $\rP$ and $\rP'$ of types $\tau$ and $\tau'$ respectively, which satisfy $\rP(k) = \rP'(k)$, hence $\rP = \rP'$ by Zariski density of rational points in parabolics and, finally, $\tau'=\tau$.

%\vspace{3cm}
\begin{figure}[here]
\centering
\input{Maximal.pstex_t}
\caption{Compactified apartment in $\overline{\mathcal{B}}_{\varnothing}({\rm SL}(3),k)$}
\label{Maximal}
\end{figure}

\

\vspace{1cm}

\

\begin{figure}[here]
\centering
\input{Triangular2.pstex_t}
\caption{Compactified apartment in $\overline{\mathcal{B}}_{\tau}({\rm SL}(3),k)$, with $\tau \neq \varnothing$}
\end{figure}

\newpage
\section{Seminorm compactification for general linear groups}

We assume in this section that the non-Archimedean field $k$ is discretely valued. In the following, we study a particular compactification of the building $\mathcal{B}({\rm PGL}_{\V},k)$ of ${\rm PGL}_{\V}$, where $\V$ is a finite-dimensional $k$-vector space. From Berkovich's point of view, this is the compactification $\overline{\mathcal{B}}_{\delta}({\rm PGL_\V},k)$ associated with the flag variety ${\rm Par}_{\delta}({\rm PGL}_{\V}) = \mathbb{P}(\V)$, classifying flags of type $((0) \subset \Hr \subset \V)$, where $\Hr$ is a hyperplane of $\V$. One can give another description of this compactification as the projectivization of the cone of non-zero seminorms on $\V$, thereby extending Goldman-Iwahori's construction of the building $\mathcal{B}({\rm PGL}_{\V},k)$. This compactification of $\mathcal{B}({\rm PGL}_{\V},k)$ should be seen as the non-Archimedean analogue of the projectivization of the cone of positive semidefinite hermitian matrices for a finite-dimensional complex vector space, the latter being the ambient space for Satake compactifications of symmetric spaces.

Starting with some reminder of Berkovich's note \cite{Ber3} and of the third named author's paper \cite{Wer1}, we give an elementary description of the compactified building $\overline{\mathcal{B}}_{\delta}({\rm PGL}_\V,k)$ and make everything explicit: convergence of seminorms, strata, stabilizers. An important feature of this compactification is the existence of a canonical retraction $\tau : \mathbb{P}(\V)^{\rm an} \rightarrow \overline{\mathcal{B}}_{\delta}({\rm PGL_\V},k)$.

%We start with some reminders of Berkovich's note \cite{Ber3} and give an elementary description of the compactification $\overline{\mathcal{B}}_{\delta}(\PGL (\V), k)$ coming from the flag variety $\Pp(\V) = {\rm Par}_{\delta}({\rm PGL}_{\V})$ in terms of (homothety classes
%of) seminorms on the $k$-vector space $\V$. An fundamental fact is the existence of a canonical retraction $\tau: \Pp(\V)^{\an} \hookrightarrow \overline{\mathcal{B}}_{\delta}({\rm PGL}_{\V},k)$. Given an absolutely irreducible projective representation $\rho : \G \rightarrow \PGL (\V)$, we then define a continuous and $\G(k)$-equivariant map $\underline{\rho}: \mathcal{B}(\G,k) \rightarrow \overline{\mathcal{B}}_{\delta}({\rm PGL}_{\V},k)$, compatible with scalar extension and whose image is contained in $\mathcal{B}({\rm PGL}_{\V},k)$. We finally prove the main result of this section: the map $\underline{\rho}$ extends to the compactification of type $t(\check{\rho})$ and induces a homeomorphism between $\overline{\mathcal{B}}_{t(\check{\rho})}(\G,k)$ and the closure of
%$\underline{\rho}(\mathcal{B}(\G,k))$ into $\overline{\mathcal{B}}_{\delta}({\rm PGL}_{\V},k)$, where $t(\check{\rho})$ is the type of the contragredient representation $\check{\rho}$.

\vspace{0.2cm}
\noindent \textbf{(3.1)} Let $\Sr^{\bullet}\V$ be the symmetric algebra of the $k$-vector space
$\V$. This is a graded $k$-algebra of finite type whose spectrum
(whose homogeneous spectrum, respectively) is the affine space $\mathbb{A}(\V)$
(the projective space $\Pp(\V)$, respectively): $$\mathbb{A}(\V) =
\mathrm{Spec}(\Sr^{\bullet}\V) \ \textrm{ and } \  \Pp(\V) =
\mathrm{Proj}(\Sr^{\bullet}\V).$$

The underlying set of the $k$-analytic space $\mathbb{A}(\V)^{\an}$ consists of all multiplicative seminorms on $\Sr^{\bullet}\V$. The underlying set of the $k$-analytic space $\Pp(\V)^{\an}$ is the quotient of
$\mathbb{A}(\V)^{\an} - \{0\}$ by homothety: two non-zero seminorms $x,y$ are equivalent if there
exists a positive real number $\lambda$ such that $|f|(y)=\lambda^n
|f|(x)$ for any natural integer $n$ and any element $f \in \Sr^n\V$.

Let $\mathcal{S}(\V,k)$ be the set of all seminorms on the vector
space $\V$ and let $\mathcal{X}(\V,k)$ be the quotient of $\mathcal{S}(\V,k) - \{0\}$
by homothety: two non-zero
seminorms $x$ and $y$ on $\V$ are equivalent if there exists a positive real number
$\lambda \in \mathbb{R}_{>0}$ such that $|f|(y)=\lambda |f|(x)$ for any $f \in \V$. Since each (multiplicative)
seminorm on $\Sr^{\bullet}\V$ induces a seminorm on $\V = \Sr^1\V$ by
restriction, we have a natural map $\tau: \mathbb{A}(\V)^{\rm an}
\rightarrow \mathcal{S}(\V,k)$ such that $\tau(x)=0$ if and only if
$x=0$. This map is obviously compatible with the above
equivalence relations and therefore descends to a map $\tau: \Pp(\V)^{\an}
\rightarrow \mathcal{X}(\V,k)$.

A seminorm $x$ on the $k$-vector space $\V$ is \emph{diagonalizable}
if there exists a basis $(e_0, \ldots, e_d)$ of $\V$ such that for every $v = \sum_{0 \leqslant i \leqslant d} a_i e_i$ in $\V$, $$|v|(x) = \max_i |a_i||e_i|(x).$$

\begin{Prop}
\label{prop.diagonalize}
Any non-zero seminorm on the $k$-vector space $\V$ is diagonalizable.
\end{Prop}

\vspace{0.1cm}
\noindent \emph{\textbf{Proof}}. As the absolute value of $k$ is assumed to be
discrete, this fact is established by F. Bruhat and J. Tits in
\cite[Proposition 1.5 (i)]{BT2}. It was initially proved by A. Weil in the locally compact case. \hfill $\Box$

\vspace{0.1cm}
Diagonalizability of seminorms on $\V$ allows us to define a canonical
section $j$ for both maps $\tau$. Given a point $x$ in $\mathcal{S}(\V,k) - \{0\}$, pick a
diagonalizing basis $(e_0, \ldots, e_d)$ of $\V$ and consider the
multiplicative seminorm defined on $\Sr^{\bullet}\V$ by
$$\sum_{\nu \in \mathbb{N}^d} \lambda_\nu e^\nu \mapsto \max_\nu
|\lambda_\nu|\prod_{i=0}^{d}|e_i|(x)^{\nu_i}.$$
For any
multiplicative seminorm $z$ on $\Sr^{\bullet}\V$ inducing $x$ on $\V$, we have:
$$|e^\nu|(z) = \prod_{0 \leqslant i \leqslant d}|e_i|(z)^{\nu_i} = \prod_{0
\leqslant i \leqslant d}|e_i|(x)^{\nu_i},$$ hence $$\left|\sum_{\nu} \lambda_\nu
e^\nu \right|(z) \leqslant \max_{\nu} |\lambda_\nu||e^{\nu}|(z) =
\max_{\nu}|\lambda_\nu|\prod_{i=0}^{d}|e_i|(x)^{\nu_i}.$$
Thus, the seminorm which we have just defined on $\Sr^{\bullet}\V$ is maximal among multiplicative seminorms on $\Sr^{\bullet}(\V)$ inducing $x$ on $\V$ and therefore it does
not depend on the basis we picked; it will be denoted by $j(x)$. We also set $j(0) = 0$. The map
$j: \mathcal{S}(\V,k) \rightarrow \mathbb{A}(\V)^{\rm an}$ so obtained
is obviously a section of $\tau$ such that $j(x) = 0$ if and only if
$x=0$. Moreover, this map is compatible with above equivalence relations, hence
descends to a map $j: \mathcal{X}(\V,k) \rightarrow \Pp(\V)^{\an}$ which is
a section of $\tau$.

\vspace{0.1cm}
\begin{Prop}
\label{prop.retraction}
\begin{itemize}
\item[(i)] For any points $x \in \mathcal{S}(\V,k)$ and $z \in
\mathbb{A}(\V)^{\rm an}$ with $\tau(z)=x$, we have $$z \leqslant j(x).$$
\item[(ii)] If we equip the sets $\mathcal{S}(\V,k)$ and $\mathcal{X}(\V,k)$ with the
natural actions of the groups ${\rm GL}_\V$ and ${\rm PGL}_\V$
respectively, then the maps $j: \mathcal{S}(\V, k) \rightarrow
\mathbb{A}(\V)^{\rm an}$ and $\tau: \mathbb{A}(\V)^{\rm an} \rightarrow \mathcal{S}(\V,k)$
($j: \mathcal{X}(\V,k) \rightarrow \Pp(\V)^{\rm an}$ and $\tau: \Pp(\V)^{\rm an}
\rightarrow \mathcal{X}(\V,k)$, respectively) are equivariant.
\end{itemize}
\end{Prop}

\vspace{0.1cm}

\noindent \emph{\textbf{Proof}} (i) We checked this inequality in the discussion
above while defining $j$.

(ii) It is enough to prove that the maps $j:
\mathcal{S}(\V,k) \rightarrow \mathbb{A}(\V)^{\an}$ and $\tau:
\mathbb{A}(\V)^{\an} \rightarrow \mathcal{S}(\V,k)$ are
$\mathrm{GL}_\V(k)$-equivariant.
This is trivially true for $\tau$ since this map sends a seminorm on
$\Sr^{\bullet}\V$ to its restriction to $\V = \Sr^1 \V$.
For any elements $x \in \mathcal{S}(\V,k) - \{0\}$ and $g \in
\mathrm{GL}_\V(k)$, the point $z = g^{-1}j(gx)$ of
$\mathbb{A}(\V)^{\an}$ satisfies $\tau(z)=g^{-1}\tau j
(gx)=g^{-1}gx=x$, hence $g^{-1}j(gx) \leqslant j(x)$ according to (i).
Substituting $gx$ to $x$ and $g$ to $g^{-1}$ in this inequality, we
obtain $gj(x) = gj(g^{-1}gx) \leqslant j(gx)$ and therefore
$j(gx)=gj(x)$. \hfill $\Box$

\vspace{0.1cm}

In the special case of the semisimple group
$\PGL_\V$ and of the flag variety $\Pp(\V) =
\mathrm{Par}_{\delta}(\PGL_\V)$, where $\delta$ is the type of parabolic subgroups stabilizing a hyperplane in $\V$, this elementary picture provides us with
an alternative description of the general construction of \cite[2.4]{RTW1}, recalled in section 1. We thus recover the classical
realization of the building $\mathcal{B}(\PGL_\V,k)$ as the
space of norms on $\V$ up to homothety (\cite{GoldmanIwahori}, \cite{BT2}) and the
construction of a compactification in terms of seminorms
\cite{Wer1}.

\vspace{0.1cm}
\begin{Prop} \label{prop.comparison.PGL}
There exists one and only one map $\iota: \overline{\mathcal{B}}_{\delta}({\rm PGL}_\V,k)
\rightarrow \mathcal{X}(\V,k)$ such that the diagram
$$\xymatrix{\overline{\mathcal{B}}_{\delta}({\rm PGL}_\V,k) \ar@{->}[d]_{\iota}
\ar@{->}[r]^{\vartheta_{\delta} \hspace{0.5cm}} &
{\rm Par}_{\delta}({\rm PGL}_{\V},k)^{\rm an} \ar@{=}[d] \\ \mathcal{X}(\V,k)
\ar@{->}[r]_{j} & \Pp(\V)^{{\rm an}}}$$ is commutative.
This map has the following properties:
\begin{itemize}
\item[(i)] it is bijective and ${\rm PGL}_\V$-equivariant;
\item[(ii)] it identifies $\mathcal{B}({\rm PGL}_{\V},k)$ with the subset of $\mathcal{X}(\V,k)$
consisting of all homothety classes of norms on $\V$; more generally,
given a subspace $\W$ of $\V$, $\iota$ identifies the stratum
$\mathcal{B}(\V/\W, k)$ of $\overline{\mathcal{B}}_{\delta}({\rm PGL}_\V,k)$ with the subset
of $\mathcal{X}(\V,k)$ consisting of all homothety classes of seminorms on $\V$ with
kernel $\W$;
\item[(iii)] for any maximal split torus $\T$ in ${\rm PGL}_{\V}$, the map
$\iota$ identifies the compactified apartment
$\overline{\A}_{\delta}(\T,k)$ in $\overline{\mathcal{B}}_{\delta}({\rm PGL}_\V,k)$ with
the set of homothety classes of $\T$-diagonalizable seminorms on
$\V$ (i.e., seminorms which are diagonalizable in a basis of $\V$
consisting of eigenvectors for the maximal split torus in ${\rm GL}_{\V}$
lifting $\T$).
\end{itemize}
\end{Prop}

\vspace{0.1cm}
\noindent \emph{\textbf{Proof}}. If it exists, such a map $\iota$
is unique since $j$ is injective.

The existence of $\iota$ follows easily from the explicit description of the map $\vartheta_{\delta}$ recalled in (1.7). Pick a maximal split
torus $\T$ in $\PGL_\V$ and a basis $(e_0, \ldots, e_d)$ of $\V$
consisting of eigenvectors for the maximal split torus in $\mathrm{GL}_\V$
lifting $\T$. Using Proposition \ref{prop.explicit}, one sees that the map
$\vartheta_{\delta}$ realizes a bijection between the compactified
apartment $\overline{\A}_{\delta}(\T,k)$ and the subset of
$\Pp(\V)^{\an}$ consisting of homothety classes of all
multiplicative seminorms $x$ on $\Sr^{\bullet}\V$ satisfying the
following condition: there exist non-negative real numbers $c_0,
\ldots, c_d$, not all equal to zero, such that $|\sum_{\nu}
\lambda_\nu e^\nu|(x) = \max_\nu |\lambda_\nu|\prod_{0\leqslant i\leqslant
d}c_i^{\nu_i}$. The subset
$\vartheta_{\delta}\left(\overline{\A}_{\delta}(\T,k)\right)$ of $\Pp(\V)^{\an}$ is
therefore the image under $j$ of the subset
$\mathcal{X}_\T(\V,k)$ of $\mathcal{X}(\V,k)$ consisting of homothety classes of all
$\T$-diagonalizable seminorms on $\V$ (i.e., diagonalizable by the split maximal torus of ${\rm GL}_{\V}$ lifting $\T$).
Since $\overline{\mathcal{B}}_{\delta}({\rm PGL}_\V,k)$ is the union of all compactified
apartments associated with maximal tori in $\PGL_\V$, the image of
the map $\vartheta_{\delta}$ is therefore contained in the image of
$j$. This observation establishes the existence of the application
$\iota$; it also proves (iii).

The map $\iota$ is injective, because so is $\vartheta_{\delta}$. Surjectivity
follows from the fact that $\mathcal{X}(\V,k)$ is the union of the subsets
$\mathcal{X}_{\T}(\V,k)$, where $\T$ runs over the set of maximal split tori in $\PGL_\V$. To
see that the map $\iota$ is $\PGL_\V(k)$-equivariant, it suffices to observe that $\iota$ is the composition $\tau
\vartheta_{\delta}$ of two equivariant maps. Indeed, since $\tau j =
\mathrm{id}_{\mathcal{X}(\V,k)},$ $$j\tau
\vartheta_{\delta} = j \tau j \iota = j \iota$$
and thus $\tau \vartheta_{\delta}=\iota$.

We now check (ii). Let $\W$ be a linear subspace of $\V$ and
consider a seminorm $x$ on $\V$. The point $j(x)$ in $\Pp(\V)^{\an}$
belongs to the subspace $\Pp(\V/\W)^{\an}$ of $\Pp(\V)^{\an}$
if and only if the seminorm $j(x): \Sr^{\bullet}\V \rightarrow \mathbb{R}_{\geqslant
0}$ factors through the canonical homomorphism $\Sr^{\bullet}\V
\rightarrow \Sr^{\bullet}(\V/\W)$. By multiplicativity, this is the
case if and only if
$x$ vanishes identically on $\W$. Since the stratum
$\mathcal{B}({\rm PGL}_{\V/\W}, k)$ of $\overline{\mathcal{B}}_{\delta}({\rm PGL}_\V,k)$ is the
preimage under $\vartheta_{\delta}$ of the space
$$\Pp(\V/\W)^{\an} - \bigcup_{\W \subsetneq \W' \subsetneq \V}
\Pp(\V/\W')^{\an},$$ we conclude that $\iota$ identifies this stratum with the
subspace of $\mathcal{X}(\V,k)$ consisting of homothety classes of seminorms
on $\V$ with kernel $\W$; in particular, this map is a bijection
between $\mathcal{B}({\rm PGL}_{\V},k)$ and the set of homothety classes of \emph{norms}
on $\V$. \hfill $\Box$

%\vspace{0.1cm} \begin{Rk} \label{rk.apartment.pgl} Let $(e_i)_{1 \leqslant i \leqslant d}$ be a basis of $\V$ and denote by $\T$ the corresponding torus of ${\rm GL}_{\V}$. The map $\iota : \overline{\mathcal{B}}_{\delta}(\V,k) \rightarrow \mathcal{X}(\V,k)$ identifies the compactified apartment $\overline{\A}_{\delta}(\T,k)$ with the set $\mathcal{N}_{\T}(\V,k)$ of non-zero seminorms on $\V$ diagonalized by this basis.
%\end{Rk}

\vspace{0.1cm}
We can introduce a natural topology on $\mathcal{X}(\V,k)$: equip the
set $\mathcal{S}(\V,k)$ with the coarsest topology such that each
evaluation map $(x \mapsto |v|(x)$, $v \in \V)$ is continuous and
consider the quotient topology on $\mathcal{X}(\V,k)$. The map $\tau:
\Pp(\V)^{\an} \rightarrow \mathcal{X}(\V,k)$ is obviously continuous. If the field $k$ is locally
compact, then the map $j: \mathcal{X}(\V,k) \rightarrow
\Pp(\V)^{\an}$ is continuous (see point (ii) below).

\vspace{0.1cm}
\begin{Prop} \label{prop.compactification.PGL} The set $\mathcal{X}(\V,k)$ is equipped with the topology which we
have just defined.
\begin{itemize}
\item[(i)] The map $\iota: \overline{\mathcal{B}}_{\delta}({\rm PGL}_\V,k) \rightarrow
\mathcal{X}(\V,k)$ is continuous and, for any maximal split torus $\T$ in
${\rm PGL}_{\V}$, it induces a homeomorphism between the compactified
apartment $\overline{\A}_{\delta}(\T,k)$ and the subspace
$\mathcal{X}_{\T}(\V,k)$ of $\mathcal{X}(\V,k)$ consisting of homothety classes of
$\T$-diagonalizable
seminorms on $\V$.
\item[(ii)] If $k$ is locally compact, the map $\iota$ is a
homeomorphism and the map $j: \mathcal{X}(\V,k) \rightarrow \Pp(\V)^{\rm an}$
is a homeomorphism onto its image.
\end{itemize}
\end{Prop}

\vspace{0.1cm}
\noindent \emph{\textbf{Proof}}. (i) Continuity of $\iota$ is obvious
if we write this map as the composition $\tau
\vartheta_{\delta}$.
Given a maximal split torus $\T$ in $\PGL_\V$, the map $\iota$ induces a
continuous bijection between the compact space
$\overline{\A}_{\delta}(\T)$ and its image in $\mathcal{X}(\V,k)$; this map
is a homeomorphism since the topological space $\mathcal{X}(\V,k)$ is
Hausdorff.

(ii) If the field $k$ is locally compact, the topological space
$\overline{\mathcal{B}}_{\delta}({\rm PGL}_\V,k)$ is compact and the continuous bijection $\iota$
onto the Hausdorff topological space $\mathcal{X}(\V,k)$ is a homeomorphism.
The map $\vartheta_{\delta}$ is a homeomorphism onto its image; writing
the map $j$ as the composition $\vartheta_{\delta} \iota^{-1}$, we see that the
same is true for $j$. \hfill $\Box$

\vspace{0.1cm}
The topology which we consider on $\mathcal{X}(\V,k)$ is relevant only if the field $k$ is locally compact.
In general, we have to modify it and endow $\mathcal{X}(\V,k)$ with the topology deduced from $\overline{\mathcal{B}}_{\delta}({\rm PGL}_\V,k)$ via the bijection $\iota$.
Equivalently, pick a maximal split torus $\T$ in $\PGL_\V$, endow $\mathcal{X}_{\T}(\V,k)$ with the coarsest topology such that all evaluations $(x \mapsto |v|(x)$, $v \in \V)$ are continuous and equip $\mathcal{X}(\V,k)$ with the quotient topology deduced from the surjective map
$$\G(k) \times \mathcal{X}_{\T}(\V,k) \rightarrow \mathcal{X}(\V,k), \ (g,x) \mapsto g \cdot x.$$

\vspace{0.1cm} The above identification between $\overline{\mathcal{B}}_{\delta}({\rm PGL}_\V,k)$ and $\mathcal{X}(\V,k)$ allows us to describe the subgroup of $\PGL_\V$ fixing a given point $x$ of $\mathcal{X}(\V,k)$. Let $\W$ be the kernel of $x$ and let $\rP$ be the parabolic subgroup of $\PGL_{\V}$ stabilizing $\W$. The subgroup of $\PGL_\V (k) (k)$ fixing $x$ is contained in $\rP(k)$; this is the extension of the maximal bounded subgroup of $\PGL_{\V/\W}(k)$ fixing the norm (induced by) $x$ on $\V/\W$ by the subgroup of $\rP(k)$ acting trivially on $\W$.

More explicitly, if $(e_0,  \ldots, e_d)$ is a basis of $\V$ diagonalizing $x$ and chosen so that $\W = {\rm Span}(e_m, \ldots, e_d)$, then ${\rm P}(k)$ is the subgroup of lower triangular block matrices $$\left(\begin{array}{ll} {\rm GL}(m,k) & \hspace{1cm} 0 \\ \ast & {\rm GL}(d+1-m, k) \end{array} \right)$$ modulo homothety. Moreover, if the basis can be chosen so that $x$ satisfies $|e_i|(x) = 1$ for any $i \in \{0, \ldots, m-1\}$, i.e., if $x$ is a vertex of $\mathcal{X}(\V/\W,k)$, then its stabilizer in ${\rm PGL}_\V (k)$ is a conjugate of the subgroup of matrices $$\left(\begin{array}{ll} k^{\times} \cdot {\rm GL}(m, k^{\circ}) & \hspace{1cm} 0 \\ \ast &  {\rm GL}(d+1-m, k)\end{array}\right)$$ modulo homothety.

\vspace{0.2cm} \noindent \textbf{(3.2)} Assuming that the field $k$ is locally compact, we complete our description of $\mathcal{X}(\V,k) \cong \overline{\mathcal{B}}_{\delta}({\rm PGL}_\V,k)$ in terms of seminorms. We fix a basis $(e_0, \ldots, e_d)$ of $\V$ and denote by $\T$ and $\widetilde{\T}$ the corresponding split maximal tori in ${\rm PGL}_{\V}$ and ${\rm GL}_{\V}$ respectively. We also denote by $o$ the norm on $\V$ defined by $$\left|\sum_{i=0}^{d} a_i e_i \right|(o) = \max_{0\leqslant i \leqslant d} |a_i|$$ and set $$\K(o) = \{g \in {\rm GL}_\V(k) \ ; \ g \cdot o = o\}.$$

\begin{Prop} \label{prop.representatives} A complete set of representatives for the action of  ${\rm GL}_{\V}(k)$ on $\mathcal{S}(\V,k) - \{0\}$ consists of all non-zero $\widetilde{\T}$-diagonalizable seminorms $x$ on $\V$ satisfying $0 \leqslant |e_d|(x) \leqslant \ldots \leqslant |e_1|(x) \leqslant |e_0|(x) \leqslant q$, where $q > 1$ generates the group $|k^{\times}|$.

(ii) The set $\overline{\mathcal{C}}$ of non-zero $\widetilde{\T}$-diagonalizable seminorms $x$ on $\V$ satisfying $0 \leqslant |e_d|(x) \leqslant \ldots \leqslant |e_1|(x) \leqslant |e_0|(x)$ is a fundamental domain for the $\K(o)$-action on $\mathcal{S}(\V,k) - \{0\}$.
\end{Prop}

\noindent \textbf{\emph{Proof}}. (i) Since each seminorm on $\V$ is diagonalizable by some maximal split torus, it follows from conjugacy of maximal split tori that each orbit of ${\rm GL}_{\V}(k)$ in $\mathcal{S}(\V,k) - \{0\}$ meets the set $\mathcal{S}_{\widetilde{\T}}(\V,k)$ of non-zero $\widetilde{\T}$-diagonalizable seminorms.

Let $\varpi$ be a generator of the maximal ideal of $k^{\circ}$, i.e., $|\varpi| = q^{-1} < 1$ generates $|k^{\times}|$, and pick $\nu \in \mathbb{N}^{d+1}$. By definition of the ${\rm GL}_{\V}$-action on $\mathcal{S}(\V,k) - \{0\}$, ${\rm diag}(\varpi^\nu)\cdot o$ is the $\widetilde{\T}$-diagonalizable seminorm on $\V$ such that $$|e_i|({\rm diag}(\varpi^\nu)\cdot o) = |{\rm diag}(\varpi^{-\nu}) \cdot e_i|(o) = |\varpi^{-\nu_i}e_i|(o) = q^{\nu_i}.$$
Accordingly, for any permutation $w \in \mathfrak{S}_{d+1}$ the permutation matrix $n(w)$ maps a $\widetilde{\T}$-diagonalizable seminorm $x$ to the $\widetilde{\T}$-diagonalizable seminorm $n(w)\cdot x$ satisfying $$|e_i|(n(w) \cdot x) = |n(w)^{-1}\cdot e_i|(x) = |e_{w^{-1}(i)}|(x).$$
Combining these two observations, one checks immediately that each ${\rm GL}_{\V} (k)$-orbit in $\mathcal{S}(\V,k) - \{0\}$ meets the subset of $\mathcal{S}_{\widetilde{\T}}(\V,k)$ consisting of seminorms $x$ such that $$0 \leqslant |e_d|(x) \leqslant \ldots \leqslant |e_1|(x) \leqslant |e_0|(x) \leqslant q.$$

\vspace{0.1cm}
\noindent (ii) As in (i), one easily shows that any $\K(o)$-orbit meets $\overline{\mathcal{C}}$.

For any point $x$ in $\mathcal{X}_{\T}(\V,k)$, we can extend $x$ to a seminorm on the exterior algebra $\Lambda^\bullet \V$ as follows: defining as usual $e_{\I}$ as the product $e_{i_1} \wedge \ldots \wedge e_{i_m}$ for any subset $\I = \{i_1, \ldots, i_m\}$ of $\{0, \ldots, d\}$ with $i_1 < \ldots < i_m$, we set $$|e_{\I}|(x) = \prod_{i \in \I} |e_{i}|(x) \ \ \textrm{ and } \ \ \left|\sum_{\I} a_{\I} e_{\I}\right|(x) = \max_{\I} |a_{\I}| \cdot |e_{\I}|(x).$$  Pick $x$ in $\mathcal{S}_{\widetilde{\T}}(\V,k)$ and assume that we have $g\cdot x \in \mathcal{S}_{\widetilde{\T}}(\V,k)$ for some $g \in \K(o)$. If we use the basis $(e_0, \ldots, e_d)$ to identify $\V$ with $k^{d+1}$, then $\K(o)$ is the subgroup ${\rm GL}(d+1,k^{\circ})$ of ${\rm GL}(d+1,k)$. For each $m \in \{1, \ldots, d\}$, this observation implies immediately $$\max_{\I} |e_{\I}|(g \cdot x) = \max_{\I} |{\small \Lambda}^m g^{-1} \cdot e_{\I}|(x) = \max_{\I} |e_{\I}|(x),$$ where the maximum is taken over all subsets $\I \subset \{0, \ldots, d\}$ of cardinality $m$. If we assume that both $x$ and $g \cdot x$ belong to $\overline{\mathcal{C}}$, it follows recursively that $|e_i|(g \cdot x)  = |e_i|(x)$ for any $i \in \{0, \ldots, d\}$, hence $g\cdot x = x$. Therefore, each $\K(o)$-orbit contains a unique point lying in $\overline{\mathcal{C}}$.
\hfill $\Box$

\vspace{0.2cm} \noindent \emph{Convergence of seminorms up to homothety}. We examine now the convergence of sequences in $\mathcal{X}(\V,k)$, from which one can recover that $\mathcal{X}(\V,k)$ is a compactification of the Bruhat-Tits building $\mathcal{B}({\rm PGL}_{\V},k)$.

\vspace{0.1cm} Let $(z_n)$ be a sequence of $\widetilde{\T}$-diagonalizable seminorms. We say that this sequence is \emph{normalized from below} if $|e_i|(z_n) \geqslant 1$ for all $i \in \{0, \ldots, d\}$ and all $n \geqslant 0$ such that $|e_i|(z_n) \neq 0$. Furthermore, we say that $(z_n)$ is \emph{distinguished} if there exists a non-empty subset $\I$ of $\{0,\ldots, d\}$ such that: \begin{itemize} \item[(a)] for any $i,j \in \I$, the sequence $\left(\frac{|e_i|(z_n)}{|e_j|(z_n)}\right)_n$ converges to a positive real number; \item[(b)] for any $i \in \I$ and $j \in \{0,\ldots, d\} - \I$, the sequence $\left(\frac{|e_j|(z_n)}{|e_i|(z_n)}\right)_n$ converges to $0$.\end{itemize}

In this situation, we set $\frac{|e_i|(z_{\infty})}{|e_j|(z_{\infty})} = \lim_n \left(\frac{|e_i|(z_n)}{|e_j|(z_n)}\right)$ for any $i,j \in \I$ and we say that $\I$ is the \emph{index set at infinity} of the sequence $(z_n)$.

\vspace{0.1cm}
The following proposition describes the convergence of sequences in $\mathcal{X}_{\T}(\V,k)$. We recall that $\overline{\mathcal{C}}$ denotes the subset of $\mathcal{S}_{\widetilde{\T}}(\V,k)$ consisting of seminorms $x$ satisfying $0 \leqslant |e_d|(x) \leqslant \ldots \leqslant |e_1|(x) \leqslant |e_0|(x)$.

\begin{Prop} Let $(x_n)$ be a sequence of points in $\mathcal{X}_{\T}(\V,k)$.
\begin{itemize}
\item[(i)] Up to going over to a subsequence, there exists a sequence $(z_n)$ in $\mathcal{S}_{\widetilde{\T}}(\V,k)$ lifting $(x_n)$ and an element $w$ of $\mathfrak{S}_{d+1}$ such that the sequence $(n(w)z_n)$ is normalized from below, distinguished and contained in $\overline{\mathcal{C}}$.
\item[(ii)] Assume that $(x_n)$ comes from a sequence $(z_n)$ of points in $\overline{\mathcal{C}}$ normalized from below and distinguished, with index set at infinity $\I$. We have $\lim (x_n) = x_{\infty}$, where $x_{\infty}$ is the homothety class of the $\widetilde{\T}$-diagonalizable seminorms $z^i_{\infty}$ defined by picking an element $i$ of $\I$ and setting $$|e_j|(z_\infty^i) = \left\{ \begin{array}{ll} \frac{|e_j|(z_{\infty})}{|e_i|(z_{\infty})} & \textrm{ if } j \in \I; \\ 0 & \textrm{ if } j \notin \I. \end{array} \right. $$
\item[(iii)] The topological space $\mathcal{X}(\V,k)$ is metrizable and compact. It contains the Bruhat-Tits building of ${\rm PGL}_{\V}(k)$ as a dense open subset.
\end{itemize}
\end{Prop}

\vspace{0.1cm} \noindent \textbf{\emph{Proof}}. (i) Let $(z_n)$ be any sequence in $\mathcal{S}_{\widetilde{\T}}(\V,k)$ lifting $(x_n)$. The seminorm $z_n$ is non-zero, so the real number $\mu_n$, defined as the minimum of the finite set $\{|e_i|(z_n) \ ; \ 0 \leqslant i \leqslant d \ \textrm{ and } |e_i|(z_n) \neq 0\}$, is positive.
For $\lambda_n = \mu_n^{-1}$, the sequence $\{\lambda_n \cdot z_n \}_{n \geqslant 0}$ is normalized from below.
It is therefore enough to show that any sequence in $\mathcal{S}_{\widetilde{\T}}(\V,k)$ which is normalized from below admits a distinguished subsequence, up to multiplication by a permutation matrix $n(w)$.
For simplicity, let us denote again by $(z_n)$ such a sequence.

For each $n \geqslant 0$, there exists $i_n \in \{0, 1, \ldots, d \}$ such that $|e_{i_n}|(z_n) = \max_{0 \leqslant i \leqslant d} \{ |e_i|(z_n)\}$.
The sequence $(i_n)_{n}$ takes its values in a finite set, so up to extracting, we may assume that it is constant.
By iterating the same argument, we find $w \in \mathfrak{S}_{d+1}$ such that:
$$|e_{w(0)}|(z_n) \geqslant |e_{w(1)}|(z_n) \geqslant \ldots \geqslant |e_{w(d)}|(z_n)$$ for any $n \geqslant 0$, that is such that the sequence $(n(w^{-1})\cdot z_n)$ lies in $\overline{\mathcal{C}}$.

Note that since $(z_n)$ is normalized from below, we have $|e_{w(0)}|(z_n) \geqslant 1$ for each $n \geqslant 0$.
For each $i \in \{0, 1, \ldots, d\}$, let us set $\beta_i = \limsup_n \frac{|e_{w(i)}|(z_n)}{|e_{w(0)}|(z_n)}$; we have:
$1=\beta_0 \geqslant \beta_1 \geqslant ... \geqslant \beta_d \geqslant 0$.
Up to extracting, we may assume that $ \lim_n \left(\frac{|e_{w(i)}|(z_n)}{|e_{w(0)}|(z_n)}\right) = \beta_i$ for each $i$. Define $\I$ as the subset of $\{0, \ldots, d\}$ consisting of indices $i$ such that $\beta_i > 0$; note that $\I$ contains $0$ by assumption, hence is non-empty. For any $i,j \in \I$, the sequence $$\frac{|e_{w(i)}|(z_n)}{|e_{w(j)}|(z_n)} = \frac{|e_{w(i)}|(z_n)}{|e_{w(0)}|(z_n)} \cdot \frac{|e_{w(0)}|(z_n)}{|e_{w(j)}|(z_n)}$$ converges to the positive real number $\frac{\beta_i}{\beta_j}$, whereas for any $i \in \I$ and $j \in \{0, \ldots, d\} - \I$ the sequence $$\frac{|e_{w(j)}|(z_n)}{|e_{w(i)}|(z_n)} = \frac{|e_{w(j)}|(z_n)}{|e_{w(0)}|(z_n)} \cdot \frac{|e_{w(0)}|(z_n)}{|e_{w(i)}|(z_n)}$$ converges to $\frac{\beta_j}{\beta_i} = 0$. Thus, the sequence $(n(w^{-1})\cdot z_n)$ is distinguished.

\vspace{0.1cm}
(ii) Let $(z_n)$ be a sequence in $\mathcal{S}_{\widetilde{\T}}(\V,k)$ lifting $(x_n)$, which we assume to be normalized from below and distinguished. Let $\I$ denote its index set at infinity. Since $$\frac{|e_\ell|(z_\infty^i)}{|e_{\ell}|(z^j_{\infty})} = \frac{|e_\ell|(z_{\infty})}{|e_i|(z_{\infty})} \cdot \frac{|e_j|(z_{\infty})}{|e_\ell|(z_{\infty})} = \frac{|e_j|(z_{\infty})}{|e_i|(z_{\infty})}$$ for any $i, j, \ell \in \I$, the $\widetilde{\T}$-diagonalizable seminorms $z_{\infty}^i$ and $z_{\infty}^j$ define the same homothety class in $\mathcal{X}(\V,k)$. Given $i \in \I$, the seminorm $y^i_n = |e_i|(z_n)^{-1} \cdot z_n$ satisfies $$\lim_{n} |e_{\ell}|(y^i_n) = \lim_n \frac{|e_{\ell}|(z_n)}{|e_i|(z_n)} = \left\{ \begin{array}{ll} \frac{|e_\ell|(z_{\infty})}{|e_i|(z_{\infty})} & \textrm{ if } \ell \in \I \\ 0 & \textrm{ if } \ell \in \{0, \ldots, d\} - \I\end{array}\right.$$ since the sequence $(z_n)$ is distinguished and thus the sequence $(y_n^i)$ converges to the seminorm $z^i_{\infty}$ in $\mathcal{S}(\V,k) - \{0\}$.

\vspace{0.1cm}
(iii) Let $k_0$ denote a dense and countable subfield of $k$ and let $\V_0$ be a $k_0$-vector subspace of $\V$ such that $\V = \V_0 \otimes_{k_0} k$; this is a dense and countable subset of $\V$. Each non-zero seminorm on $\V$ is completely determined by its restriction to $\V_0$, hence the map $$\mathcal{S}(\V,k) \rightarrow \mathbb{R}^{\V_0}, \ x \mapsto (v \mapsto |v|(x))$$ is a continuous injection. Since $\mathcal{S}(\V,k)$ is locally compact, this injection is a homeomorphism of $\mathcal{S}(\V,k)$ onto its image. This map induces a homeomorphism of $\mathcal{X}(\V,k)$ onto a subspace of $\mathbb{R}^{\V_0}/\mathbb{R}_{>0}$ and, since the latter topological space is metrizable, so is $\mathcal{X}(\V,k)$.

It follows from (ii) that the image $\overline{\mathcal{Q}}$ of $\overline{\mathcal{C}}$ in $\mathcal{X}(\V,k)$ is compact. The map $\pi : \K(o) \times \overline{\mathcal{Q}} \rightarrow \mathcal{X}(\V,k)$ induced by the ${\rm GL}_{\V} (k)$-action is continuous, and it is surjective by Proposition \ref{prop.representatives}, (ii). Since $\K(o) \simeq {\rm GL}_{d+1}(k^{\circ})$, the source is compact; as the target is Hausdorff, compactness of $\mathcal{X}(\V,k)$ follows.

Identifying the Bruhat-Tits building $\mathcal{B}({\rm PGL}_{\V},k)$ with the subspace of $\mathcal{X}(\V,k)$ consisting of classes of norms on $\V$, the complementary subspace $\mathcal{X}(\V,k) - \mathcal{B}({\rm PGL}_{\V},k) = \K(o) \cdot \left[\overline{\mathcal{Q}} \cap (\mathcal{X}(\V,k) - \mathcal{B}({\rm PGL}_{\V},k))\right]$ is closed and therefore $\mathcal{B}({\rm PGL}_{\V},k)$ is open in $\mathcal{X}(\V,k)$. Density is obvious. \hfill $\Box$

%\vspace{0.1cm} More generally, if $(x_n)$ is an arbitrary sequence of points in $\mathcal{X}(\V,k)$, then Proposition \ref{prop.representatives} (i) furnishes a sequence $(g_n)$ in ${\rm PGL}_{\V}$ such that $g_n \cdot x_n$ belongs to $\mathcal{X}_{\T}(\V,k)$ for each $n \geqslant 0$. Up to extracting a subsequence, we may assume that $(g_n \cdot x_n)$ converges to a point $z$ in $\mathcal{X}_{\T}(\V,k)$ by Proposition \ref{prop.convergence} (i).

\vspace{0.2cm} \noindent \emph{Orbit structure}. We have already observed in Proposition \ref{prop.comparison.PGL} that the canonical identification $\mathcal{X}(\V,k) \cong \overline{\mathcal{B}}_{\delta}({\rm PGL}_\V,k)$ transforms the natural stratification of $\overline{\mathcal{B}}_{\delta}({\rm PGL}_\V,k)$ into the stratification of $\mathcal{X}(\V,k)$ by kernels: with each point $x$ of $\mathcal{X}(\V,k)$ is associated the non-zero linear subspace $\V(x) = \{v \in \V \ ; \ |v|(x) = 0\}$ and two points $x, y \in \mathcal{X}(\V,k)$ belong to the same stratum if $\V(x) = \V(y)$. The set of strata is indexed by the set of non-zero linear subspaces of $\V$ and the stratum associated with a linear subspace $\W$ is canonically isomorphic to the building $\mathcal{B}({\rm PGL}_{\V/\W},k)$.

Given any point $x$ of $\mathcal{X}(\V,k)$, its stabilizer in ${\rm PGL}_{\V} (k)$ is the extension of a maximal compact subgroup of ${\rm PGL}_{\V/\V(x)}(k)$ by ${\rm PGL}_{\V(x)}(k)$, and its Zariski closure is the parabolic subgroup fixing $\V(x)$.

All these assertions can be easily proved starting from the definition of $\mathcal{X}(\V,k)$, without knowing the structure of the Berkovich compactification $\overline{\mathcal{B}}_{\delta}({\rm PGL}_\V,k)$. One can also show that the unique closed orbit for the ${\rm PGL}_{\V}$-action on $\mathcal{X}(\V,k)$ consists of the homothety classes of seminorms of the form $|.| \circ \varphi$, where $\varphi$ is a non-zero linear form on $\V$; this orbit is ${\rm PGL}_{\V}(k)$-equivariantly homeomorphic to $\mathbb{P}(\V)(k)$, i.e., to the set of hyperplanes in $\V$.

\vspace{0.2cm} \noindent \textbf{(3.3)} We end this section on the compactified
building $\overline{\mathcal{B}}_{\delta}({\rm PGL}_\V,k)$ with a couple of technical results to be used in the next paragraph.

Recall that, for any Banach $k$-algebra $\A$ and any non-Archimedean
extension $\K/k$, the formula $$||f|| = \inf \left\{\max_{i \in \I}
|\lambda_i| \cdot ||f_i||; \ \lambda_i \in \K, f_i \in \A \ \textrm{and}
\ f = \sum_{i \in \I} f_i \otimes \lambda_i \right\}$$ defines a
seminorm on the $\K$-algebra $\A \otimes_k \K$ and that $\A
\widehat{\otimes}_k \K$ is the Banach $\K$-algebra one gets by
completion \cite[2.1.7 and 3.4.3]{BGR}. The following definition
is due to Berkovich \cite[Sect. 5.2]{Ber1}.

\begin{Def} \label{def.peaked} Let $\X$ be a $k$-analytic space. A point $x$ in $\X$ is
\emph{peaked} if, for any non-Archimedean extension $\K/k$, the norm
on the Banach $\K$-algebra $\mathcal{H}(x) \widehat{\otimes}_k \K$
is multiplicative.

Let $x$ be a peaked point of $\X$. For any non-Archimedean extension
$\K/k$, the norm on $\mathcal{H}(x) \widehat{\otimes}_{k} \K$
defines a point in $\mathcal{M}(\mathcal{H}(x) \widehat{\otimes}_{k}
\K)$ and $\sigma_{\K}(x)$ denotes its image under the canonical map
$\mathcal{M}(\mathcal{H}(x) \widehat{\otimes}_k \K) \rightarrow \X
\widehat{\otimes}_k \K$.
\end{Def}

\begin{Rk} \label{rk.peaked} For a point
$x$ in a $k$-analytic space $\X$, being peaked or not depends only on the completed residue field $\mathcal{H}(x)$.

%2. The norm of a $k$-affinoid algebra $\A$ is multiplicative if and only if the Shilov boundary of $\X = \mathcal{M}(\A)$ consists of a unique point. It follows that the norm of $\A$ is \emph{universally multiplicative}, i.e., the norm of $\A \widehat{\otimes}_k \K$ is multiplicative for any non-Archimedean extension $\K/k$, if and only if the Shilov boundary of $\X = \mathcal{M}(\A)$ consists of a unique point, which is peaked.

\end{Rk}

\vspace{0.1cm}
\begin{Lemma} \label{lemma.peaked} For any point $x$ in $\Pp(\V)^{\rm an}$, there exists a
point $y$ in $\mathbb{A}(\V)^{\an}$ lifting $x$ and such that
$\mathcal{H}(x)=\mathcal{H}(y)$. In particular, each peaked point $x$
in $\Pp(\V)^{\rm an}$ can be lifted to a peaked point in
$\mathbb{A}(\V)^{\rm an}$.
\end{Lemma}

\vspace{0.1cm}
\noindent \emph{\textbf{Proof}}. This is obvious since the canonical map
$\mathbb{A}(\V)(\K) - \{0\} \rightarrow \Pp(\V)(\K)$ is surjective for any
field extension $\K/k$. \hfill $\Box$

\vspace{0.1cm}
\begin{Prop} \label{prop.peaked} Let $x$ be a peaked point of
$\Pp(\V)^{\rm an}$. For any discretely valued non-Archimedean field $\K$
extending $k$, the canonical injection of $\mathcal{X}(\V,k)$ into
$\mathcal{X}(\V,\K)$ maps the point $\tau(x)$ to the point
$\tau(\sigma_{\K}(x))$.
\end{Prop}

\vspace{0.1cm}
\noindent \emph{\textbf{Proof}}. Consider a peaked point $y$ in
$\mathbb{A}(\V)^{\an}$ lifting $x$ and denote by $\tau(y)_{\K}$ the
image of $\tau(y)$ under the canonical injection $\mathcal{X}(\V,k)
\rightarrow \mathcal{X}(\V,\K)$. We want to show: $\tau(y)_{\K} = \tau(\sigma_{\K}(y))$.

The point $\sigma_{\K}(y)$ in
$\mathbb{A}(\V \otimes_k \K)^{\an}$ is the multiplicative seminorm
on $\Sr^{\bullet}(\V \otimes_k \K) = (\Sr^{\bullet}\V) \otimes_k \K$
defined by
$$|f|(\sigma_{\K}(y)) = \inf \left\{\max_{i \in \I} |\lambda_i||f_i|(y)
; \ \lambda_i \in \K, f_i \in \Sr^{\bullet}\V \ \textrm{and} \ f =
\sum_{i \in \I} f_i \otimes \lambda_i\right\}.$$ Hence
$$|f|(\tau(\sigma_{\K}(y))) = \inf \left\{\max_{i \in \I}
|\lambda_i||f_i|(y); \ \lambda_i \in \K, f_i \in \Sr^{1}\V = \V \
\textrm{and} \ f = \sum_{i \in \I} f_i \otimes \lambda_i\right\}$$
for any $f \in \V \otimes_k \K$.

Pick a basis $(e_0, \ldots, e_d)$ of $\V$ diagonalizing $\tau(y)$.
Given $f = \sum_{i \in \I} f_i \otimes \lambda_i$ in $\V \otimes_k
\K$, we can write $f_i = \sum_{0 \leqslant j \leqslant d} a_{ij}e_j$, and \begin{eqnarray*}
\max_{i \in \I} |\lambda_i||f_i|(y) & = & \max_{i \in \I}
|\lambda_i|\max_{0 \leqslant j \leqslant d}|a_{ij}||e_j|(y) \\ & = & \max_{0
\leqslant j \leqslant d} \max_{i \in \I} |\lambda_i||a_{ij}||e_j|(y) \\ & \geqslant
& \max_{0 \leqslant j \leqslant d} \left|\sum_{i \in \I} \lambda_i
a_{ij}\right||e_j|(y).
\end{eqnarray*}
We conclude that $$\max_{i \in \I} |\lambda_i||f_i|(y) \geqslant
|f|(\tau(y)_{\K}),$$ hence $|f|(\tau(y)_{\K}) \leqslant
|f|(\tau\sigma_{\K}(y))$.

The converse inequality is obvious: for any $f = \sum_{0 \leqslant i \leqslant
d}a_i e_i$ in $\V \otimes_k \K$, \begin{eqnarray*}|f|(\tau
\sigma_{\K}(y)) = |f|(\sigma_\K (y)) & \leqslant & \max_{0 \leqslant i \leqslant d}
|a_i||e_i|(\sigma_{\K}(y)) \\ & \leqslant & \max_{0 \leqslant i \leqslant d}
|a_i||e_i|(y) = |f|(\tau(y)_{\K})| \end{eqnarray*} and we finally get
$$\tau(y)_{\K} = \tau \sigma_{\K}(y).$$ \hfill $\Box$

\section{Satake compactifications via Berkovich theory}

In \cite{Sa}, Satake considers a Riemannian symmetric space $\Sr = \G/\K$ of non-compact type. Using a faithful representation $\rho$ of the real Lie group $\G$ in ${\rm PSL}(n,\mathbb{C})$, he embeds $\Sr$ in the symmetric space $\Hr$ associated with ${\rm PSL}(n,\mathbb{C})$, which can be identified with the space of all positive definite hermitian $n \times n$-matrices of determinant $1$. Observing that $\Hr$ has a natural compactification $\overline{\Hr}$, namely the projectivization of the cone of all positive semidefinite hermitian $n \times n$-matrices, Satake defines the compactification of $\Sr$ associated with $\rho$ as the closure of $\Sr$ in $\overline{\Hr}$.

In this section and the next one, we present an analogous construction for Bruhat-Tits buildings from two different viewpoints. Let $\G$ be a semisimple connected group over a discretely valued non-Archimedean field $k$. A faithful and absolutely irreducible linear representation $\rho : \G \rightarrow {\rm GL}_{\V}$ of $\G$ in some finite dimensional $k$-vector space $\V$ can be used to embed the building  of $\G$ in the building of ${\rm SL}_{\V}$, hence in any compactification of the latter, and we get a compactification of $\mathcal{B}(\G,k)$ by taking the closure. The Berkovich compactification of $\mathcal{B}({\rm SL}_{\V},k)$ corresponding to parabolics stabilizing a hyperplane has an elementary description as the space of seminorms up to scaling on $\V$ and will be the non-Archimedean analogue of the projective cone of semidefinite hermitian matrices.

\vspace{0.1cm} The difference between this section and the next one lies in the construction of the map from $\mathcal{B}(\G,k)$ to $\mathcal{B}({\rm SL}_\V,k)$. Whereas functoriality of buildings is a delicate question in general, it is quite remarkable that Berkovich theory allow us to attach very easily and in a completely canonical way a map $\underline{\rho} : \mathcal{B}(\G,k) \rightarrow \mathcal{B}({\rm PGL}_{\V},k)$ to each absolutely irreducible linear representation $\rho : \G \rightarrow {\rm GL}_{\V}$. General results of E. Landvogt on functoriality of buildings will be used in the next section.

\vspace{0.4cm} \noindent \textbf{(4.1)} {\itshape The map
$\underline{\rho}: \mathcal{B}(\G,k) \rightarrow \mathcal{X}(\V,k)$.} Let
$\G$ be a semisimple connected $k$-group and consider a projective
representation $\rho: \G \rightarrow \PGL_\V$, which we assume to
be \emph{absolutely irreducible}. We start by showing that the
morphism $\rho$ naturally leads to a continuous and
$\G(k)$-equivariant map $\underline{\rho}: \mathcal{B}(\G,k) \rightarrow
\mathcal{X}(\V,k)$, whose formation commutes with scalar extension
and whose image lies in the building $\mathcal{B}({\rm PGL}_{\V},k)$.

\vspace{0.1cm}
The two main ingredients in the definition of $\underline{\rho}$ are
the retraction $\tau: \Pp(\V)^{\an} \rightarrow \mathcal{X}(\V,k)$,
defined in 3.1, and the following well-known fact.

\vspace{0.1cm}
\begin{Prop} \label{prop.type.rep}
\begin{itemize}
\item[(i)] For any field extension $\K/k$ and any Borel subgroup $\B$ of
$\G \otimes_k \K$, there exists one and only one $\K$-point of
$\Pp(\V)$ invariant under $\B$.
\item[(ii)] There exists a unique $k$-morphism $\widetilde{\rho}:
{\rm Bor}(\G) \rightarrow \Pp(\V)$ such that: for any field
extension $\K/k$, the map $\widetilde{\rho}_{\K} : {\rm Bor}(\G)(\K) \rightarrow \Pp(\V)(\K)$ sends a Borel subgroup $\B$ to the unique $\K$-point of $\Pp(\V)$
invariant under~$\B$.
\end{itemize}
\end{Prop}

\vspace{0.1cm}
\noindent \emph{\textbf{Proof}}. We use the following two results:
\begin{enumerate}
\item If the field $k$ is algebraically closed, then for
each Borel subgroup $\Br \in \mathrm{Bor}(\G)(k)$ there exists one and only one
point in $\Pp(\V)(k)$ invariant under $\Br(k)$ \cite[Expos\'e 15, Proposition 6 and Corollaire 1]{Bible}.
\item If the group $\G$ is split over $k$, then for each
Borel subgroup $\Br \in \mathrm{Bor}(\G)(k)$ there exists at least one point in
$\Pp(\V)(k)$ invariant under $\Br(k)$ \cite[Proposition 15.2]{Borel}.
\end{enumerate}

(i) Let $\K/k$ be a field extension, pick an algebraic closure
$\K^a$ of $\K$ and consider the separable closure $\K^s$ of $\K$ in
$\K^a$. Given a Borel subgroup $\Br$ in $\G \otimes_k \K$, assertion 2 provides a
$\K^s$-point of $\Pp(\V)$, say $x$, invariant under the group $\Br(\K^s)$.
Since the $\K$-scheme $\Br$ is smooth, the subset $\Br(\K^s)$ is
dense in $\Br$, hence $x$ is invariant under
$\Br(\K^a)$ and assertion 1 provides uniqueness of this point.

For any $\gamma \in \mathrm{Gal}(\K^s/\K)$, the point
$\gamma \cdot x$ in $\Pp(\V)(\K^s)$ is invariant under the group $\gamma
\Br(\K^s) = \Br(\K^s)$; uniqueness implies $\gamma \cdot x = x$ and
therefore this point belongs to the subset $\Pp(\V)(\K)$ of $\Pp(\V)(\K^s)$. We have thus established existence and uniqueness of a
$\Br(\K^a)$-invariant point in $\Pp(\V)(\K)$. We still have to
check that this point is fixed by $\Br$, i.e., that  its image
in $\Pp(\V)(\Sr)$ is invariant under the group $\Br(\Sr)$ for any $\K$-scheme
$\Sr$.

\vspace{0.1cm} \emph{First step} --- {\itshape The functor
  $\K$-$\mathbf{Sch} \rightarrow \mathbf{Sets}, \ \Sr \mapsto
\mathrm{Stab}_{\G(\Sr)}(x)$ is representable by a closed
subgroup, say $\Pi$, of $\G$.}

As a direct verification shows, the functor $\K$-$\mathbf{Sch}
\rightarrow \textbf{Sets}, \ \Sr \mapsto \mathrm{Stab}_{\PGL_\V(\Sr)}(x)$, is represented by a closed and smooth subgroup $\rP_0$ of
$\PGL_\V$. Let $\Pi$ denote the $\K$-scheme $\rP_0 \times_{\PGL_\V} \G$.
The second projection $\Pi \rightarrow \G$ is
a closed immersion and $\Pi$ represents the functor
$\mathrm{Stab}_{\G}(x)$ since $$\Pi(\Sr) = \{(g,g') \in \G(\Sr)
\times \rP_0(\Sr) \ ; \ \rho(g)=g'\} = \mathrm{Stab}_{\G(\Sr)}(x)$$
for any $\K$-scheme $\Sr$.

\vspace{0.1cm} \emph{Second step} --- \emph{The subgroup $\Br$ of
$\G$ is contained in $\Pi$.}

Since $\Br$ is a reduced closed subscheme of $\G$, the
inclusion $\Br(\K^a) \subset \Pi(\K^a)$ implies the inclusion
$\Br \subset \Pi$ as subgroups of $\G$ and we have thus established
that the $\K$-point $x$ of $\Pp(\V)$ is invariant under $\Br$.
Note also that $\Pi$ (which may not be smooth) is a generalized parabolic subgroup of $\G$ since it
contains a Borel subgroup.

\vspace{0.2cm} (ii) Pick a finite Galois extension $k'/k$ splitting
$\G$ together with a Borel subgroup $\Br$ of $\G \otimes_k k'$, and
let $x$ be the only $k'$-point of $\Pp(\V)$ invariant under $\Br$.
By (i), the map $$\G(\Sr) \rightarrow
\Pp(\V)(\Sr), \ g \mapsto g \cdot x$$ factors through the canonical
projection $\G(\Sr) \rightarrow \G(\Sr)/\Br(\Sr)$ for any $k'$-scheme $\Sr$.
Thanks to the functorial identification $\G(\Sr)/\Br(\Sr)
\tilde{\rightarrow} \mathrm{Bor}(\G)(\Sr), \ g \Br(\Sr) \mapsto
g(\B \otimes_{k'} \Sr)g^{-1}$ \cite[Expos\'e XXVI, Corollaire 5.2]{SGA3} we thus get a
morphism of functors $\widetilde{\rho}: \mathrm{Bor}(\G \otimes_k
k') \rightarrow \Pp(\V \otimes_k k')$ and define therefore a
$k'$-morphism $\widetilde{\rho}: \mathrm{Bor}(\G \otimes_k k')
\rightarrow \Pp(\V \otimes_k k')$ such that, for any $k'$-scheme
$\Sr$ and any $\Br' \in \mathrm{Bor}(\G)(\Sr)$,
$\widetilde{\rho}(\Br')=g \cdot x$ if $\Br' = g\Br g^{-1}$, $g \in
\G(\Sr)$. In particular, for any field extension $\K/k$, the map
$\widetilde{\rho}$ associates with a Borel subgroup $\Br' \in
\mathrm{Bor}(\G)(\K)$ the only $\K$-point of $\Pp(\V)$ invariant
under $\Br'$.

By definition, the $k'$-morphism
$$\widetilde{\rho}: \mathrm{Bor}(\G) \otimes_k k' = \mathrm{Bor}(\G
\otimes_k k') \rightarrow \Pp(\V
\otimes_k k') = \Pp(\V) \otimes_k k'$$ commutes with the natural
action of $\mathrm{Gal}(k'|k)$ and thus $\widetilde{\rho}$ descends
to a $k$-morphism
$$\widetilde{\rho}: \mathrm{Bor}(\G) \rightarrow \Pp(\V)$$
satisfying the required condition.\hfill $\Box$

\vspace{0.2cm}
\begin{Prop} \label{prop.type.rep.2} There exists a largest type $t$ of parabolic subgroups of $\G$ such that the morphism $\widetilde{\rho} : {\rm Bor}(\G) \rightarrow \mathbb{P}(\V)$ factors through the canonical projection ${\rm Bor}(\G) \rightarrow {\rm Par}_t(\G)$. The so-obtained morphism ${\rm Par}_t(\G) \rightarrow \mathbb{P}(\V)$ induces a homeomorphism between ${\rm Par}_t(\G)^{\rm an}$ and a closed subspace of $\mathbb{P}(\V)^{\rm an}$.
\end{Prop}

\vspace{0.1cm} \noindent \textbf{\emph{Proof}.} Assume temporarily that the group $\G$ is split, pick a Borel
subgroup $\Br$ of $\G$ and let $x = \widetilde{\rho}(\B)$ be the only $k$-point of
$\Pp(\V)$ invariant under $\Br$. If we denote by $\Pi$ the stabilizer of $x$ in $\G$, then the underlying reduced scheme $\Pi^{\rm red}$ is the largest parabolic subgroup of $\G$ stabilizing $x$. Indeed, since we have proved above that $\Pi$ is a closed subgroup containing $\B$, the reduced scheme $(\Pi \otimes_k k^{a})^{\rm red}$ is a smooth closed subgroup of $\G \otimes_k k^a$ containing $\B \otimes_k k^{a}$, hence a parabolic subgroup of $\G \otimes_k k^a$. As $\G$ is split, there exists a unique parabolic subgroup $\rP$ of $\G$ containing $\B$ such that $(\Pi \otimes_k k^a)^{\rm red} = \rP \otimes_k k^a$. This identity implies $\rP = \Pi^{\rm red}$, hence $\Pi^{\rm red}$ is a parabolic subgroup of $\G$ stabilizing $x$. Since each parabolic subgroup $\Qr$ of $\G$ is smooth, $\Qr$ is a subgroup of $\Pi$ if it stabilizes $x$, and therefore $\Pi^{\rm red}$ contains any parabolic subgroup of $\G$ stabilizing $x$. Note also that the type of $\Pi^{\rm red}$ does not depend on the choice of $\B$ by $\G(k)$-conjugacy of Borel subgroups and equivariance of the map $\widetilde{\rho}$.

\vspace{0.1cm} The morphism $\widetilde{\rho}: \G/\B \rightarrow
\Pp(\V)$ induces a map $$\G/\Pi \hookrightarrow \mathbb{P}(\V)$$ which is a monomorphism in the category of $k$-schemes. Since the image of $\widetilde{\rho}$ is a closed subset of $\mathbb{P}(\V)$ by properness of ${\rm Bor}(\G)$, this map is a closed immersion. Moreover, we have an exact sequence of $k$-groups $$\xymatrix{e \ar@{->}[r] & \Pi/\Pi^{\rm red} \ar@{->}[r] & \G/\Pi^{\rm red} \ar@{->}[r]^p & \G/\Pi \ar@{->}[r] & e }$$ and $\Pi/\Pi^{\rm red}$ is a finite and connected $k$-group scheme \cite[Expos\'e VIA, 5.6]{SGA3}, hence the morphism $p$ is universally injective, i.e., induces an injection between $\K$-points for any extension $\K$ of $k$. Let $t$ denote the rational type of $\G$ defined by $\Pi^{\rm red}$. Composing $p$ with the morphism $\G/\Pi \rightarrow \mathbb{P}(\V)$ induced by $\rho$, we see that $\widetilde{\rho}$ factors through the canonical projection of ${\rm Bor}(\G)$ onto ${\rm Par}_t(\G)$. The induced morphism $f : {\rm Par}_t(\G) \rightarrow \mathbb{P}(\V)$ is  universally injective. At the analytic level, the associated map $f^{\rm an}$ is a continuous injection, hence a homeomorphism onto a closed subset of $\mathbb{P}(\V)^{\rm an}$ since ${\rm Par}_t(\G)^{\rm an}$ and $\mathbb{P}(\V)^{\rm an}$ are compact.

\vspace{0.1cm} In general, we pick a finite Galois extension $k'/k$ splitting $\G$ and set $\Gamma = {\rm Gal}(k'|k)$. For any $\gamma \in \Gamma$, there exists a unique $k'$-rational type $t'_{\gamma}$ such that the morphism $^\gamma \widetilde{\rho}_k' = \widetilde{\rho} \otimes_k \gamma$ factors through ${\rm Par}_{t'_{\gamma}}(\G \otimes_k k')$. The family $\{t'_{\gamma}\}_{\gamma \in \Gamma}$ is a Galois orbit, hence defines a type $t$ of parabolic subgroups of $\G$, and the morphism $\widetilde{\rho}$ factors through the canonical projection of ${\rm Bor}(\G)$ onto ${\rm Par}_t(\G)$ by Galois descent. \hfill $\Box$

\vspace{0.1cm} The above construction associates a well-defined rational type of parabolic subgroups of $\G$ with the representation $\rho$.

\begin{Def} \label{def.type.rep} Let $\rho$ be an absolutely irreducible projective representation $\G \rightarrow {\rm PGL}_{\V}$. Its \emph{cotype} $t(\check{\rho})$ is the largest rational type $t$ of $\G$ such that the canonical morphism $\widetilde{\rho} : {\rm Bor}(\G) \rightarrow \mathbb{P}(\V)$ factors through the projection of ${\rm Bor}(\G)$ onto ${\rm Par}_t(\G)$.
\end{Def}

\vspace{0.1cm} \begin{Rk} \label{rk.highest.weight} This definition is obviously related to the theory of the highest weight: if $\B$ is a Borel subgroup of $\G$, then the $k$-point $\widetilde{\rho}(\B)$ of $\mathbb{P}(\V)$ is a hyperplane of $\V$ invariant under $\B$, hence a line in $\V^{\vee}$ invariant under $\B$ in the contragredient representation $\check{\rho}$. The corresponding character of $\B$ is the highest weight of $\check{\rho}$ with respect to $\B$. This observation is the reason why we introduced the \emph{cotype} of the representation $\rho$; the \emph{type} of $\rho$ should be defined as the cotype of the contragredient representation, i.e., the type of the largest parabolic subgroup stabilizing a highest weight line in $\V$.
\end{Rk}

\vspace{0.2cm}
Composing the maps $$\xymatrix{\mathcal{B}(\G,k) \ar@{->}[r]^{\vartheta_{\varnothing}} &
\mathrm{Bor}(\G)^{\an} \ar@{->}[r]^{\widetilde{\rho}} &
\Pp(\V)^{\an} \ar@{->}[r]^{\tau} & \mathcal{X}(\V,k),}$$ we obtain a natural map $$\underline{\rho}:
\mathcal{B}(\G,k) \rightarrow \mathcal{X}(\V,k),$$ canonically associated with the homomorphism $\rho: \G
\rightarrow \PGL_\V$. Since all these maps are continuous and equivariant, so is $\underline{\rho}$.

\vspace{0.2cm} \noindent \textbf{(4.2)} The main properties of $\underline{\rho}$ are easily established. We first consider compatibility
with scalar extension.

\begin{Prop} \label{prop.rep.map.scalar} For any discretely valued non-Archimedean field $\K$ extending $k$,
the natural diagram $$\xymatrix{\mathcal{B}(\G,\K)
\ar@{->}[r]^{\underline{\rho_{\K}}} & \mathcal{X}(\V,\K) \\
\mathcal{B}(\G,k) \ar@{->}[u] \ar@{->}[r]_{\underline{\rho}} &
\mathcal{X}(\V,k) \ar@{->}[u]}$$ is commutative.
\end{Prop}

\vspace{0.1cm}
The proof of this proposition relies on the following lemma. We recall that, if $x$ is a peaked point of a $k$-analytic space $\X$ and if $\K/k$ is a non-Archimedean extension, then $\sigma_{\K}(x)$ denotes the canonical lift of $x$ to $\X \widehat{\otimes}_k \K$ (see Definition \ref{def.peaked}).

\begin{Lemma}
\label{lemma.rep.peaked}
For any rational type $t$ of $\G$ and any point $x$ in $\overline{\mathcal{B}}_t (\G,k)$, the point $\vartheta_t (x)$ of
$\mathrm{Par}_t (\G)^{\rm an}$ is peaked.
Moreover, given a non-Archimedean extension $\K/k$, the point $\sigma_{\K}(\vartheta_t(x))$
of $\mathrm{Par}_t (\G)^{\rm an} \widehat{\otimes}_k \K$ is the image of $x_{\K}$ under the map
$$\vartheta_t: \overline{\mathcal{B}}_t (\G, \K) \rightarrow {\rm Par}_t(\G \otimes_k \K)^{\rm an}
= {\rm Par}_t (\G)^{\rm an} \widehat{\otimes}_k \K.$$
\end{Lemma}

\vspace{0.1cm}
\noindent \emph{\textbf{Proof}}. Let us first consider a finite Galois extension $k'/k$ splitting $\G$ and consider a point $x'$ in
$\overline{\mathcal{B}}_t(\G,k')$. By Proposition 1.13, the point $\vartheta_t(x')$ is contained in some big cell $\Omega$ of ${\rm Par}_t(\G \otimes_k k')$. Choosing an isomorphism $\mathbb{G}_{{\rm a}, k'} \tilde{\rightarrow} \U_{\alpha}$ for each root $\alpha$ of $\G \otimes_k k'$
with respect to a maximal split torus $\T$ containing $\Sr \otimes_k k'$ leads to an isomorphism
$\mathbb{A}^n_{k'} \tilde{\rightarrow} \Omega \otimes_k k'$.
Then the point $\vartheta_t(x')$ corresponds to a seminorm
on the algebra $k'[\xi_1, \ldots, \xi_n]$ of the form
$$\sum_{\nu} a_{\nu} \xi^{\nu} \mapsto \max_{\nu} |a_{\nu}| \prod_{i=1}^{n}c_i^{\nu_i},$$ where $c_1, \ldots, c_n$ are
non-negative real numbers, not all equal to zero (with the convention $0^0 = 1$). Such a seminorm defines a peaked point in $\mathbb{A}^{n, \mathrm{an}}_{k'}$
\cite[Sect. 5.2]{Ber1} and the point $\vartheta_t(x')$ is therefore peaked.

In general, pick a point $x$ in $\overline{\mathcal{B}}_t(\G,k)$ and let $x_{k'}$ denote its image in $\mathcal{B}(\G,k')$, where $k'/k$ is a finite Galois extension splitting $\G$. We consider the completed residue field $\mathcal{H}(\vartheta_t(x))$ of $\vartheta_t(x)$. The point $\vartheta_t(x_{k'})$ induces a norm on the $k'$-Banach algebra $\mathcal{H}(\vartheta_t(x)) \otimes_k k'$ with respect to which the descent datum is an isometry (note that $\mathcal{H}(\vartheta_t(x)) \otimes_k k'$ is finite extension of $k'$). Since the point $\vartheta_t(x_{k'})$ is peaked, this norm is universally multiplicative. By \cite[Lemma A.10]{RTW1}, it follows that the norm induced on $\mathcal{H}(\vartheta_t(x))$ is also universally multiplicative, hence the point $\vartheta_t(x)$ is peaked.

\vspace{0.1cm} In order to prove the second assertion, consider a point $x$ in $\overline{\mathcal{B}}_t(\G,k)$ and let $\K/k$ be a non-Archimedean extension. Since the point $\vartheta_t(x)$ is peaked, the Banach norm on the $\K$-Banach algebra $\mathcal{H}(\vartheta_t(x)) \widehat{\otimes}_k\K$ coming from the absolute value of $\mathcal{H}(\vartheta_t(x))$ is multiplicative. On the other hand, the point $\vartheta_t(x_{\K})$ also defines a multiplicative norm on this $\K$-Banach algebra. Two such norms necessarily coincide, hence $\sigma_{\K}(\vartheta_t(x)) = \vartheta_t(x_{\K})$. \hfill $\Box$

\vspace{0.1cm} \noindent \emph{Proof of Proposition 4.5.} Let $\K$
be a discretely valued non-Archimedean field extending $k$. Denoting
by $t$ the cotype of the representation $\rho$, the morphism
$\widetilde{\rho}: \mathrm{Bor}(\G) \rightarrow \Pp(\V)$ factors
through the canonical projection $\mathrm{Bor}(\G) \rightarrow
\mathrm{Par}_t (\G)$ and leads to a homeomorphism between
$\mathrm{Par}_t(\G)^{\rm an}$ and a closed subset of $\Pp(\V)^{\rm an}$ (Proposition \ref{prop.type.rep.2}).
Pick a point $x$ in $\mathcal{B}(\G,k)$. The point
$\widetilde{\rho}(\vartheta_t (x))$ of $\Pp(\V)^{\an}$ is peaked since
$\mathcal{H}(\widetilde{\rho}\vartheta_t (x)) =
\mathcal{H}(\vartheta_t (x))$ and
$\vartheta_t (x)$ is a peaked point of $\mathrm{Par}_t (\G)^{\an}$ (Lemma 4.6).
Moreover, we have the identities
$$\sigma_\K\widetilde{\rho}\vartheta_t (x) =
\widetilde{\rho_\K}\sigma_\K \vartheta_t (x) = \widetilde{\rho_\K}
\vartheta_t(x).$$ The conclusion finally follows from Proposition
\ref{prop.peaked}: the points $\underline{\rho}(x) = \tau
\widetilde{\rho}\vartheta_t (x)$ and $\underline{\rho_\K}(x) = \tau
\widetilde{\rho_\K} \vartheta_t(x) = \tau \sigma_\K
\widetilde{\rho}\vartheta_t (x)$ coincide in $\mathcal{X}(\V,\K)$.
\hfill $\Box$

\vspace{0.1cm}
\begin{Prop} \label{prop.rep.open.stratum} The image of the map $\underline{\rho}: \mathcal{B}(\G,k)
\rightarrow \mathcal{X}(\V,k)$ is contained in the open stratum
$\mathcal{B}({\rm PGL}_{\V},k)$ of $\mathcal{X}(\V,k)$.
\end{Prop}

\vspace{0.1cm}
\noindent \emph{\textbf{Proof}}. Assume that there exists a point $x$ in
$\mathcal{B}(\G, k)$ whose image under the map $\underline{\rho}$ is not
contained the open stratum $\mathcal{B}({\rm PGL}_{\V},k)$ of $\mathcal{X}(\V,k)$.
Under this hypothesis, the point
$\underline{\rho}(x) = \tau \widetilde{\rho}
\vartheta_{\varnothing}(x)$ lies in $\mathcal{X}(\V,k) \cap
\Pp(\V/\W)^{\an}$ for some non trivial linear subspace $\W$ in
$\V$, hence $\widetilde{\rho} \vartheta_{\varnothing} (x) \in
\Pp(\V/\W)^{\an}$. Now consider the following diagram
$$\xymatrix{\mathrm{Bor}(\G)^{\an} \ar@{->}[r]^{\widetilde{\rho}}
\ar@{->}[d] & \Pp(\V)^{\an} \ar@{->}[d] \\ \mathrm{Bor}(\G)
\ar@{->}[r]_{\widetilde{\rho}} & \Pp(\V)}$$ in which the vertical
arrows are the maps sending a point $z$ of $\X^{\an}$, seen as a
multiplicative seminorm on the algebra $\mathcal{O}_{\X}(\U)$ of
some open affine subset $\U$ of $\X$, to the point of the scheme
$\X$ defined by the prime ideal $\mathrm{ker}(z) \in
\mathrm{Spec}(\mathcal{O}_{\X}(\U))$ (where $\X = \mathrm{Bor}(\G),$ or $\X=\Pp(\V)$). The point $x$ ($\widetilde{\rho}(x)$, respectively) is mapped to
the generic point of $\mathrm{Bor}(\G)$ (to the generic point
of $\Pp(\V/\W)$, respectively). Since the diagram above is commutative, it
follows that the morphism $\widetilde{\rho}$ maps the generic point
of $\mathrm{Bor}(\G)$ to the generic point of $\Pp(\V/\W)$, hence
maps $\mathrm{Bor}(\G)$ into the strict linear subspace $\Pp(\V/\W)$
of $\Pp(\V)$. Hence it would follow that
$\rho$ maps $\G$ into the nontrivial parabolic subgroup of
$\PGL_\V$ stabilizing $\Pp(\V/\W)$, thereby contradicting the irreducibility of $\rho$.\hfill $\Box$

\vspace{0.8cm} \noindent \textbf{(4.3)} We now state and prove the
main result of this section.

\begin{Thm} \label{thm.Berkovich.Satake} Let $k$ be a discretely valued non-Archimedean field and
 $\G$  a semisimple connected $k$-group. We consider a finite-dimensional
$k$-vector space $\V$ and an absolutely irreducible projective
representation $\rho: \G \rightarrow {\rm PGL}_{\V}$.
\begin{itemize}
\item[(i)] The map $\underline{\rho}: \mathcal{B}(\G,k) \rightarrow
\mathcal{X}(\V,k)$ extends continuously to the compactification $\mathcal{B}(\G,k)
\hookrightarrow \overline{\mathcal{B}}_{t(\check{\rho})} (\G,k)$.
\item[(ii)] The induced map is an injection of
$\overline{\mathcal{B}}_{t(\check{\rho})}(\G,k)$ into $\mathcal{X}(\V,k)$.
\item[(iii)] If the field $k$ is locally compact, the map $\underline{\rho}$
extends to a homeomorphism between $\overline{\mathcal{B}}_{t(\check{\rho})} (\G,k)$ and the
closure of $\underline{\rho}(\mathcal{B}(\G,k))$ in $\mathcal{X}(\V,k)$.
\end{itemize}
\end{Thm}

\vspace{0.1cm}
\noindent \emph{\textbf{Proof}}. Set $t = t(\check{\rho})$. %We first reduce to the case of a split group. Let $k'/k$ be a finite Galois extension
%splitting $\G$ and consider the commutative diagram
%$$\xymatrix{\mathcal{B}(\G,k') \ar@{->}[r]^{\underline{\rho_{k'}}} &
%\mathcal{X}(\V,k') \\ \mathcal{B}(\G,k) \ar@{->}[r]_{\underline{\rho}}
%\ar@{->}[u] & \mathcal{X}(\V,k). \ar@{->}[u]}$$ The vertical map on
%the right is injective and maps $\mathcal{X}(\V,k)$ homeomorphically
%onto its image. Denoting by $t'$ the cotype of the absolutely irreducible representation $\rho_{k'}$, the vertical map on the left extends to a topological embedding of $\overline{\mathcal{B}}_t(\G,k)$ into
%$\overline{\mathcal{B}}_{t'}(\G,k')$. It is clear that the theorem holds for $\G$
%if it holds for $\G \otimes_k k'$, hence we assume $\G$ to be
%split in what follows.

\vspace{0.1cm} (i) The morphism $\widetilde{\rho}: \mathrm{Bor}(\G)
\rightarrow \Pp(\V)$ factors through the canonical projection $\pi_t
: \mathrm{Bor}(\G) \rightarrow \mathrm{Par}_t (\G)$ and leads to a
homeomorphism between ${\rm Par}_t(\G)^{\rm an}$ and a closed subset
of $\mathbb{P}(\V)^{\rm an}$ (Proposition \ref{prop.type.rep.2}). The diagram
$$\xymatrix{\mathcal{B}(\G,k) \ar@{->}[r]^{\vartheta_{\varnothing}}
\ar@{->}[rd]_{\vartheta_t} & \mathrm{Bor}(\G)^{\an}
\ar@{->}[r]^{\widetilde{\rho}} \ar@{->}[d]^{\pi_t} & \Pp(\V)^{\an}
\ar@{->}[r]^{\tau} & \mathcal{X}(\V,k) \\ & \mathrm{Par}_t(\G)^{\an}
\ar@{^{(}->}[ru]_{\widetilde{\rho}} & & }$$ is commutative (use \cite[section 4.2]{RTW1} for the left-hand side triangle) and hence allows us to write
the map $\underline{\rho}$ as the composition $\tau \widetilde{\rho}
\vartheta_{t}$. For any maximal split torus $\Sr$ in $\G$, the restriction
of $\underline{\rho}$ to the apartment $\A(\Sr,k)$ extends continuously
to its closure $\overline{\A}_t(\Sr,k)$ in
$\mathrm{Par}_t(\G)^{\an}$. Since the image of
$\overline{\mathcal{B}}_t(\G,k)$ into $\mathrm{Par}_t(\G)^{\an}$ is the
union of these closures when $\Sr$ runs over all maximal split tori of
$\G$, the maps $\underline{\rho}$ extends to
$\overline{\mathcal{B}}_t(\G,k)$. This extension is continuous, for it is
$\G(k)$-equivariant and its restriction to $\overline{\A}_t(\Sr,k)$ is continuous.

\vspace{0.1cm} (ii) Let us now prove that the map
$\overline{\mathcal{B}}_t(\G,k) \rightarrow \mathcal{X}(\V,k)$ extending
$\underline{\rho}$, for which we keep the notation
$\underline{\rho}$, is injective. The fact that compatibility of $\underline{\rho}$ with scalar extension is proved only for discretely valued non-Archimedean extensions of $k$ in Proposition \ref{prop.rep.map.scalar} is a slight difficulty.

Given two points $x, y \in \overline{\mathcal{B}}_t(\G,k)$ with $\underline{\rho}(x) = \underline{\rho}(y)$, we will show that $\G_x(k^{\rm a}) = \G_{y}(k^{\rm a})$, where $\G_x = {\rm Stab}_{\G}^{\ t}(x)$ and $\G_y = {\rm Stab}_{\G}^{\ t}(y)$. Since the field $k$ is discretely valued, it follows from its description as a disjoint union of buildings (cf. Theorem \ref{thm.stratification}) that the
compactified building $\overline{\mathcal{B}}_t(\G,k)$ carries a (poly-)simplicial
decomposition and, by application of Bruhat-Tits theory to each stratum, the fixed-point set of ${\rm Stab}^t_{\G}(x)(k)$ is
precisely the facet of $\overline{\mathcal{B}}_t(\G,k)$ whose interior contains the point $x$.
Now, since two distinct points of $\overline{\mathcal{B}}_t(\G,k)$ belong to disjoint facets of
$\overline{\mathcal{B}}_t (\G,k')$ for a large enough finite extension $k'/k$, the equality $\G_x(k^{\rm a}) = \G_y(k^{\rm a})$ implies $x=y$.

\vspace{0.1cm}
We pick a point $x$ in $\overline{\mathcal{B}}_t(\G,k)$ and set $\G_{\underline{\rho}(x)} = \rho^{-1}\left({\rm Stab}_{{\rm PGL}_{\V}}^\delta (\underline{\rho}(x))\right)^{\rm red}$. This is an analytic subgroup of ${\rm G}^{\rm an}$, and $$\G_{\underline{\rho}(x)} (\K) = \{g \in \G(\K) \ ; \ \rho(g)\underline{\rho}(x) = \underline{\rho}(x)\}$$ for any non-Archimedean extension $\K/k$. Given any finite extension $k'/k$, it follows from Proposition \ref{prop.rep.map.scalar} that $\G_{\underline{\rho}(x)}(k')$ contains $\G_x(k')$. We have therefore $\G_x(k^{\rm a}) \subset \G_{\underline{\rho}(x)}(k^{\rm a})$, and we will now prove that equality holds. Notice that, if the point $x$ is rational (i.e., if it becomes a vertex over some finite extension of $k$), then the inclusion $\G_x(k^{\rm a}) \subset \G_{\underline{\rho}(x)}(k^{\rm a})$ implies $\G_x \subset \G_{\underline{\rho}(x)}$ by density.

\vspace{0.1cm} \emph{Notation} --- The point $x$ belongs to a stratum $\Sr$. Let $\rP = {\rm Stab}_{\G}(\Sr)$ denote the corresponding $t$-relevant parabolic subgroup of $\G$ and let $\R = \R_t(\rP)$ denote the largest connected, smooth and normal subgroup of $\G$ acting trivially on ${\rm Osc}_t(\rP)$. Similarly, the point $\underline{\rho}(x)$ belongs to a stratum $\Sigma$ of $\mathcal{X}(\V,k)$; we set $\Pi = {\rm Stab}_{{\rm PGL}_{\V}} (\Sigma)$ and we let $\R_{\delta}(\Pi)$ denote the largest connected, smooth and normal subgroup of $\Pi$ acting trivially on $\Sigma$. Up to replacing $k$ by a finite extension, we may assume that the reduced subschemes $\rP' = \rho^{-1}(\Pi)^{\rm red}$ and $\R'' = \rho^{-1}(\R_{\delta}(\Pi))^{\rm red}$ are smooth subgroups of $\G$. Note that $\R''$ is connected and invariant in $\rP'$.

\vspace{0.1cm} \emph{First step} --- The group $\G_x(k)$ is Zariski-dense in $\rP$ (Theorem \ref{thm.stabilizer.k-points}) and $\rho$ maps $\G_{\underline{\rho}(x)}(k)$ into $\Pi(k)$. Since $\rP$ is reduced, the inclusion $\G_x(k) \subset \G_{\underline{\rho}(x)}(k)$ implies that $\rho$ maps $\rP$ into $\Pi$ and therefore $\rP' = \rho^{-1}(\Pi)^{\rm red}$ is a parabolic subgroup of $\G$ containing $\rP$.

This parabolic subgroup $\rP'$ defines a stratum $\Sr'$ in $\overline{\mathcal{B}}_t(\G,k)$, the only one it stabilizes. We have $\Sr \subset \overline{\Sr'}$ since $\rP \subset \rP'$, and $\Sr = \Sr'$ if and only if $\rP = \rP'$ for $\rP$ is $t$-relevant. In order to establish the last identity, we let $\R' = \R_t(\rP')$ denote the largest smooth connected and normal subgroup of $\rP'$ acting trivially on ${\rm Osc}_{t}(\rP')$.

\vspace{0.1cm} \emph{Second step} --- We now prove that the parabolic subgroups $\rP$ and $\rP'$ coincide.

Since $\rP' = \rho^{-1}(\Pi)^{\rm red}$, the morphism $\widetilde{\rho}$ maps the closed subscheme ${\rm Osc}_t(\rP')$ of ${\rm Par}_t(\G)$ to the closed subscheme ${\rm Osc}_{\delta}(\Pi)$ of $\mathbb{P}(\V)$. By construction, $\widetilde{\rho}$ is universally injective (i.e., purely inseparable), hence the induced map ${\rm Osc}_t(\rP')(\K) \rightarrow {\rm Osc}_{\delta}(\Pi)(\K)$ is injective for any extension $\K$ of $k$. It follows that any element $g$ of ${\rm R}''(\K)$ acts trivially on ${\rm Osc}_t(\rP')(\K)$, which implies that the action of $\R''$ on the reduced scheme ${\rm Osc}_t(\rP')$ is itself trivial. As the subgroup $\R''$ is smooth, connected and normal in $\rP'$, we deduce that $\R''$ is contained in $\R'$ by maximality of the latter. On the other hand,  ${\R'}^{\rm an}$ is trivially contained in $\R^{\rm an}$, hence in $\G_{\underline{\rho(x)}}$, since any element acting trivially on $\Sr'$ fixes $\overline{\Sr'}$ pointwise.

\vspace{0.1cm} We consider now the quotient group $\Hr =  \rP'/\R'$, which is semisimple and satisfies $\Sr' = \mathcal{B}(\Hr,k)$. Thanks to the inclusion $\R'' \subset \R'$, this group is also a quotient of $\rP'/\R''$. Since $\rP' = \rho^{-1}(\Pi)^{\rm red}$ and $\R'' = \rho^{-1}(\R_{\delta}(\Pi))^{\rm red}$, we get a canonical morphism $$p : \xymatrix{\rP'/\R''  \ar@{^{(}->}[r] & \rho^{-1}(\Pi)/\rho^{-1}(\R_{\delta}(\Pi)) \ar@{^{(}->}[r] & \Gamma = \Pi/\R_{\delta}(\Pi)}$$ which is finite. By construction, we have ${\R''}^{\rm an} \subset \G_{\underline{\rho(x)}} \subset {\rP'}^{\rm an}$ and $\G_{\underline{\rho}(x)}/{\R''}^{\rm an} = p^{-1}(\Gamma_{\underline{\rho}(x)})^{\rm red}$, hence $\G_{\underline{\rho}(x)}/{\R''}^{\rm an}$ is bounded in $(\rP'/{\R''})^{\rm an}$ for $p$ is finite. It follows that $\G_{\underline{\rho(x)}}/{\R'}^{\rm an}$ is a bounded in $\Hr^{\rm an}$.

Since $\G_x(k^{\rm a}) \subset \G_{\underline{\rho}(x)}(k^{\rm a})$, the discussion above shows that the stabilizer $(\G_x/{\R'})(k^{\rm a})$ of $x$ in ${\rm H}(k^{\rm a})$ is bounded. By Remark \ref{rk.bounded}, this amounts to saying that $x$ belongs to the open stratum of $\overline{\Sr} = \overline{\mathcal{B}}_t(\Hr,k)$, hence $\Sr' = \Sr$ and $\rP' = \rP$.

\vspace{0.1cm} \emph{Third step} ---  We have just proved that the subgroup $\G_{\underline{\rho}(x)}(k^{\rm a})$ of $\G(k^{\rm a})$ is
contained in the parabolic subgroup $\rP$ and has bounded
image in the quotient group $\Hr = \rP/\R$. The inclusion $\G_x(k^{\rm a}) \subset \G_{\underline{\rho}(x)}(k^{\rm a})$ implies
$\G_x(k^{\rm a}) = \G_{\underline{\rho}(x)}(k^{\rm a})$ since $(\G_x/\R^{\rm an})(k^{\rm a}) = \G_x(k^{\rm a})/\R(k^{\rm a})$ is a maximal bounded subgroup of $\Hr(k^{\rm a})$.

\vspace{0.1cm} (iii) If the field $k$ is locally compact, the continuous extension of $\underline{\rho} : \mathcal{B}(\G,k) \rightarrow \mathcal{X}(\V,k)$ to $\overline{\mathcal{B}}_{t(\check{\rho})}(\G,k)$ is continuous injection between two locally compact spaces, hence is a homeomorphism on its image. \hfill $\Box$

\vspace{0.2cm} \noindent \textbf{(4.3)} We end this section by establishing a natural and expected  property of $\underline{\rho}$.

\begin{Prop} \label{prop.Berkovich.toral} For any maximal split torus $\Sr$ of $\G$, there exists a maximal split torus $\T$ of ${\rm PGL}_{\V}$ containing $\rho(\Sr)$ and such that $\underline{\rho}$ maps $\A(\Sr,k)$ into $\A(\T,k)$.
\end{Prop}

\vspace{0.1cm} \noindent \textbf{\emph{Proof}}. For any finite extension $k'/k$, we normalize the metrics so that the canonical embeddings $\mathcal{B}(\G,k) \hookrightarrow \mathcal{B}(\G,k')$ and $\mathcal{B}({\rm PGL}_{\V},k) \hookrightarrow \mathcal{B}({\rm PGL}_\V,k')$ are isometric.

\vspace{0.1cm}
Given a maximal split torus $\Sr$ of $\G$, our first goal is to find an apartment $\A'$ of $\mathcal{B}({\rm PGL}_{\V},k)$ containing the image of $\A(\Sr,k)$.

Let $\T$ be a maximal split torus of ${\rm PGL}_{\V}$ containing $\rho(\Sr)$ and let $x$ be a point in $\A(\Sr,k)$. For any $s \in \Sr(k)$, we have $\underline{\rho}(s \cdot x) = \rho(s) \cdot \underline{\rho}(x)$ and $\rho(s) \cdot \A(\T,k) = \A(\T,k)$, hence $${\rm dist}\left(\underline{\rho}(s \cdot x), \A(\T,k)\right) = {\rm dist}\left(\underline{\rho}(x), \A(\T,k)\right).$$ More generally, we have $${\rm dist}\left(\underline{\rho}(s \cdot x), \A(\T,k')\right) = {\rm dist}\left(\underline{\rho}(x),\A(\T,k')\right)$$ for any finite extension $k'/k$ and any $s \in \Sr(k')$. Since the points of $\A(\Sr,k)$ belonging to the orbit of $x$ under $\Sr(k')$ for some finite extension $k'/k$ are dense (in $\A(\Sr,k)$), it follows that ${\rm dist}\left(\underline{\rho}(z),\A(\T,k)\right)$ is independent of $z \in \A(\Sr,k)$. Now, the existence of a maximal split torus $\T'$ of ${\rm PGL}_{\V}$ such that $\rho(\Sr) \subset \T'$ and $\underline{\rho}(x) \in \A(\T',k)$, hence such that $\underline{\rho}(\A(\Sr,k)) \subset \A(\T',k)$, follows immediately from the next two facts:
\begin{enumerate}
\item the set of distances of $\underline{\rho}(x)$ to apartments in $\mathcal{B}({\rm PGL}_{\V},k)$ is discrete;
\item given a maximal split torus $\T$ of ${\rm PGL}_{\V}$ such that $\rho(\Sr) \subset \T$ and $\underline{\rho}(x) \notin \A(\T,k)$, there exists a maximal split torus $\T'$ of ${\rm PGL}_{\V}$ satisfying $\rho(\Sr) \subset \T'$ and $${\rm dist}\left(\underline{\rho}(x),\A(\T',k)\right) < {\rm dist}\left(\underline{\rho}(x),\A(\T,k)\right).$$
\end{enumerate}

\vspace{0.1cm} The first assertion follows easily from the (poly-)simplicial structure on $\mathcal{B}({\rm PGL}_\V,k)$, hence from the fact that the field $k$ is discretely valued. Let us then prove the second assertion.

\vspace{0.1cm} For any point $z \in \A(\Sr,k)$, let $p(z)$ denote the unique point of $\A(\T,k)$ satisfying $${\rm dist}\left(\underline{\rho}(z),p(z)\right) = {\rm dist}\left(\underline{\rho}(z),\A(\T,k)\right)$$ and observe that the image of the map $$p : \A(\Sr,k) \rightarrow \A(\T,k), \ \ z \mapsto p(z)$$ is an affine subspace under the image of $\Lambda(\Sr)$ in $\Lambda(\T)$.

We now use the (poly-)simplicial structure on $\A(\T,k)$. Suppose that there exists a point $z \in \A(\Sr,k)$ such that p($z)$ belongs to the interior of an alcove $c$. Any path in $\mathcal{B}({\rm PGL}_{\V},k)$ from $p(z)$ to a point lying outside $\A(\T,k)$ contains an initial segment $[p(z),z']$ with $z' \in \partial c$ and $[p(z),z'[ \subset c$. Applied to the geodesic path $[p(z),\underline{\rho}(z)]$, this observation leads to a contradiction if ${\rm dist}\left(\underline{\rho}(z),\A(\T,k)\right) > 0$, since then $${\rm dist}\left(\underline{\rho}(z),\A(\T,k)\right) = {\rm dist}\left(\underline{\rho}(z),p(z)\right) < {\rm dist}\left(\underline{\rho}(z),z'\right) \leqslant {\rm dist}\left(\underline{\rho}(z),\A(\T,k)\right).$$ Therefore, since $\underline{\rho}(x) \notin \A(\T,k)$, the affine subspace $p(\A(\Sr,k))$ of $\A(\T,k)$ is contained in some root hyperplane $\Hr_{\alpha, r} = \{\alpha = r\}$, where $\alpha \in \X^*(\T)$ is a root whose restriction to $\Sr$ is trivial: $\alpha_{|\Sr} = 1$, and $r \in |k^{\times}|$. By folding $\A(\T,k)$ along $\Hr_{\alpha,r}$, we will obtain a new apartment of $\mathcal{B}({\rm PGL}_{\V}, k)$ which is closer to $\underline{\rho}(\A(\Sr,k))$.

Let $x_0=\underline{\rho}(x), x_1, \ldots, x_n = p(x)$ denote the successive vertices of the simplicial decomposition of $[\underline{\rho}(x),p(x)]$ induced by the (poly-)simplicial structure of $\mathcal{B}({\rm PGL}_{\V},k)$. There exists an element $u$ of $\U_{\alpha}(k)_r$ satisfying the following two conditions:
\begin{itemize}
\item[(a)] $\A(\T,k) \cap u \cdot \A(\T,k)$ is the half-apartment $\{\alpha \leqslant r\}$;
\item[(b)] $u \cdot \A(\T,k) = \A(u \T u^{-1},k)$ contains $[x_{n-1},x_n]$.
\end{itemize}
Since $\alpha_{|\Sr} = 1$, we have $s u s^{-1} = u$ for any $s \in \Sr(k)$ and thus $\rho(\Sr(k'))$ stabilizes the apartment $\A(u \T u^{-1},k')$ for any finite extension $k'/k$. Setting $\N = {\rm Norm}_{{\rm PGL}_{\V}}(u\T u^{-1})$, the stabilizer of $\A(u \T u^{-1},k')$ in ${\rm PGL}_{\V}(k')$ is the group $\N(k')$, hence $\rho(\Sr(k')) \subset \N(k')$ for any finite extension $k'/k$ and thus $\rho(\Sr) \subset \N$ since both $\Sr$ and $\N$ are reduced $k$-groups. By connectedness, it follows that $\Sr$ is contained in $\N^{\circ} = u \T u^{-1}$.

We have $${\rm dist}\left(\underline{\rho}(x), \A(u \T u^{-1},k)\right) \leqslant {\rm dist}(\underline{\rho}(x),x_{n-1}) = {\rm dist}(\underline{\rho}(x),x_n) - {\rm dist}(x_{n-1}, x_n) < {\rm dist}\left(\underline{\rho}(x), \A(\T,k)\right)$$ since $x_{n-1} \neq x_n$. This concludes the proof of assertion 2 above.

\vspace{0.1cm} We have just proved that there exists a maximal split torus $\T'$ of ${\rm PGL}_\V$ such that $\rho(\Sr) \subset \T'$ and $\underline{\rho}(\A(\Sr,k)) \subset \A(\T',k)$. Thanks to compatibility of $\underline{\rho}$ with finite field extensions, the inclusion $\underline{\rho}(\A(\Sr,k)) \subset \A(\T',k)$ holds more generally after any such extension. As before, it follows that $\rho(\Sr)$ is contained in $\T' = {\rm Norm}_{{\rm PGL}_{\V}}(\T')^{\circ}$ and this completes the proof.
\hfill $\Box$

\vspace{0.1cm} \begin{Rk} Given two semisimple connected $k$-groups $\G$, $\Hr$ and a homomorphism $f : \G \rightarrow \Hr$, the above proof applies more generally to any continuous and $\G(k)$-equivariant map $\mathcal{B}(\G,k) \rightarrow \mathcal{B}(\Hr,k)$ which is compatible with finite extensions of $k$: the apartment of any maximal split torus $\Sr$ of $\G$ is mapped to the apartment of a maximal split torus of $\Hr$ containing $f(\Sr)$.
\end{Rk}

\vspace{0.1cm} Functoriality of buildings with respect to group homomorphisms has been studied by Landvogt in \cite{LandvogtCrelle}. Given a complete discretely valued field $k$ with perfect residue field and two semisimple connected $k$-groups $\G$ and $\Hr$, Landvogt proved that each homomorphism $f : \G \rightarrow \Hr$ gives rise to a non-empty set of $\G(k)$-equivariant and continuous maps $f_* : \mathcal{B}(\G,k) \rightarrow \mathcal{B}(\Hr,k)$. By construction, each such map is \emph{toral}, i.e., maps the apartment of a maximal split torus $\Sr$ of $\G$ to the apartment of a maximal split torus of $\Hr$ containing $f(\Sr)$. In the special case where $\Hr = {\rm PGL}_{\V}$ and $f$ is an absolutely irreducible representation, the map $\underline{f}$ introduced in this section is an instance of Landvogt's maps.

The canonical nature of the map $\underline{f}$ raises two obvious questions: is the set of Landvogt's maps reduced to an element when $f$ is an absolutely irreducible representation? If no, is there a way to single out $\underline{f}$ without using Berkovich geometry?

\section{Satake compactifications via Landvogt's functoriality}

In this last section, we present another approach to Satake compactifications using Landvogt's results on functoriality of Bruhat-Tits buildings. As before, $\G$ is a connected, semisimple group over a non-Archimedean local field $k$. We fix a faithful, absolutely irreducible representation $\rho: \G \rightarrow {\rm GL}_{\V}$ for some finite-dimensional $k$-vector space $\V$.
Using results from \cite{LandvogtCrelle}, the representation $\rho$  defines a continuous, $\G(k)$-equivariant embedding
$\rho_\ast: \mathcal{B}(\G,k) \rightarrow \mathcal{B}({\rm SL}_{\V},k)$.

As in the previous section, we want to use one fixed compactification of $\mathcal{B}({\rm SL}_{\V},k)$ on the right-hand side and take the closure of the image of $\mathcal{B}(\G,k)$ to retrieve $\overline{\mathcal{B}}(\G,k)_{\rho}$. For functoriality reasons, the natural candidate for this compactification of $\mathcal{B}({\rm SL}_{\V},k)$ is $\ov{\mathcal{B}}({\rm SL}_{\V},k)_{\rm id}$ for the identical representation ${\rm id} : {\rm SL}_{\V} \rightarrow {\rm GL}_{\V}$. According to Theorem \ref{thm.comparison}, $\ov{\mathcal{B}}({\rm SL}_{\V},k)_{\rm id} = \ov{\mathcal{B}}_{\pi}({\rm SL}_{\V},k)$, where $\pi$ is the type of parabolics stabilizing a line in $\V$. This space was studied in \cite{Wer0} and is canonically isomorphic to $\overline{\mathcal{B}}_{\delta}({\rm SL}_{\V^{\vee}},k)$, where $\V^{\vee}$ denotes the dual vector space. It can be identified with the union of all Bruhat-Tits buildings $\mathcal{B}({\rm SL}_{\V'},k)$, where $\V'$ runs through the linear subspaces of $\V$. Its points can be described as seminorms on $\V^{\vee}$ up to scaling and vertices correspond bijectively to the homothety classes of free $k^{\circ}$-submodules (of arbitrary rank) in $\V$.

\vspace{0.1cm} In the following, we let $\tau$ denote the unique $k$-rational type such that $\overline{\mathcal{B}}(\G,k)_\rho\cong \overline{\mathcal{B}}_\tau(\G,k)$, whose existence was established in section 2. It will eventually turn out that we can replace $\tau$ by the (non necessarily $k$-rational) type $t(\rho)$ naturally associated with $\rho$.

\vspace{0.2cm}
\noindent \textbf{(5.1)} We recall some results of \cite{LandvogtCrelle}, applied to the representation $\rho: \G \rightarrow {\rm GL}_{\V}$. Since $\G$ is semisimple, it is equal to its derived group. Hence $\rho$ comes from a representation $\rho: \G \rightarrow {\rm SL}_{\V}$, for which we use the same notation.

Let $\Sr$ be a maximal split torus in $\G$ with normalizer $\N$, and let $\A(\Sr,k)$ denote the corresponding apartment in $\mathcal{B}(\G,k)$.  Choose a special vertex $o$ in $\A(\Sr,k)$. By \cite{LandvogtCrelle}, there exists a maximal split torus $\T$ in ${\rm SL}_{\V}$ containing $\rho(\Sr)$, and there exists a point $o'$ in the apartment $\A(\T,k)$ of $\T$ such that the following properties hold:

\begin{enumerate}
\item There is a unique affine map $i: \A(\Sr,k) \rightarrow \A(\T,k)$ such that $i(o) = o'$. Its linear part is induced by $\rho: \Sr \rightarrow \T$.
\item The map $i$ satisfies $\rho(\rP_x) \subset \rP'_{i(x)}$ for all $x \in \A(\Sr,k)$, where $\rP_x$ denotes the stabilizer of the point $x$ with respect to the $\G(k)$-action on $\mathcal{B}(\G,k)$, and $\rP'_{i(x)}$ denotes the stabilizer of the point $i(x)$ with respect to the ${\rm SL}_\V(k)$-action on $\mathcal{B}({\rm SL}_{\V},k)$.
\item The map $\rho_\ast: \A(\Sr,k) \rightarrow \A(\T,k) \rightarrow \mathcal{B}({\rm SL}_{\V},k)$ defined by composing $i$ with the natural embedding of the apartment $\A(\T,k)$ in the building $\mathcal{B}({\rm SL}_{\V},k)$ is $\N(k)$-equivariant, i.e., for all $x \in \A(\Sr,k)$ and $n \in \N(k)$ we have $\rho_\ast (nx) = \rho(n) \rho_\ast(x)$.
\end{enumerate}

These properties imply that $\rho_\ast: \A(\Sr,k) \rightarrow \mathcal{B}({\rm SL}_\V, k)$ can be continued to a map
$\rho_\ast: \mathcal{B}(\G,k) \rightarrow \mathcal{B}({\rm SL}_{\V},k)$, which is continuous and $\G(k)$-equivariant. By \cite[2.2.9]{LandvogtCrelle}, $\rho_\ast$ is injective and isometrical, if the metric on $\mathcal{B}(\G,k)$ is normalized correctly.

We want to show that $\rho_\ast$ can be extended to a map $\rho_\ast: \overline{\mathcal{B}}(\G,k)_{\rho} \cong \overline{\mathcal{B}}_{\tau}(\G,k) \rightarrow \overline{\mathcal{B}}_\pi({\rm SL}_\V,k)$. Besides, we prove that this map of compactified buildings identifies $\overline{\mathcal{B}}_{\tau}(\G,k)$ as a topological $\G(k)$-space with the closure of $\rho_\ast(\mathcal{B}(\G,k))$ in $\overline{\mathcal{B}}_{\pi}({\rm SL}_\V,k)$.

\vspace{0.2cm} \noindent \textbf{(5.2)} Let us first look at compactified apartments in $\overline{\mathcal{B}}_{\tau}(\G,k)$ and $\overline{\mathcal{B}}_\pi({\rm SL}_\V,k)$.

Let $(e_0,\ldots, e_d)$ be a basis of $\V$ consisting of eigenvectors of $\T$ and denote by $\chi_0, \ldots, \chi_d$ the corresponding characters of $\T$. The map $$\Lambda(\T) \rightarrow (\mathbb{R}_{>0})^{d+1}, \ \ u \mapsto (\langle u, \chi_i \rangle)_{0\leqslant i \leqslant d}$$ identifies $\Lambda(\T)$ with the subset of $(\mathbb{R}_{>0})^{d+1}$ consisting of vectors $(r_0,\ldots, r_d)$ satisfying $r_0 \ldots r_d =~1$. The fan on $\Lambda(\T)$ defining the compactification $\overline{\A}_{\pi}(\T,k)$ of $\A(\T,k)$ consists of all faces of the cones $\C_0, \ldots, \C_d$, where $$\C_i = \{(r_0,\ldots, r_d) \in (\mathbb{R}_{>0})^{d+1} \ ; \ r_0 \cdot \ldots \cdot r_d = 1 \ \textrm{ and } \ r_i \geqslant r_j, \ \textrm{ for all j}\}.$$

The weights of the representation $\rho$ with respect to the torus $\Sr$ are the images of $\chi_0, \ldots, \chi_d$ under the projection $\X^*(\T) \rightarrow \X^{*}(\Sr)$ deduced from the morphism $\rho : \Sr \rightarrow \T$, i.e., the restrictions of $\chi_0,\ldots, \chi_d$ to $\Sr$. Setting $\lambda_i = (\chi_i)_{|\Sr}$ for all $i \in \{0,\ldots, d\}$ and identifying as above $\Lambda(\T)={\rm Hom}_{\mathbf{Ab}}(\X^*(\T),\mathbb{R}_{>0})$ with a subset of $(\mathbb{R}_{>0})^{d+1}$, the dual map $$\iota : \Lambda(\Sr) = {\rm Hom}_{\mathbf{Ab}}(\X^{*}(\Sr),\mathbb{R}_{>0}) \rightarrow (\mathbb{R}_{>0})^{d+1}$$ is simply defined by $$u \mapsto \left(\langle \lambda_i, u \rangle \right)_{0 \leqslant i \leqslant d}.$$ This is an embedding since the representation $\rho$ is faithful.

\vspace{0.1cm}
\begin{Lemma} \label{lemma.comparison.fans2} The preimage under $\iota$ of the fan $\mathcal{F}$ generated by $\{\C_0,\ldots, \C_d\}$ is the fan $\mathcal{F}_{\tau}$ on $\Lambda(\Sr)$.
\end{Lemma}

\vspace{0.1cm}
\noindent \emph{\textbf{Proof}}. By definition, $$\iota^{-1}(\C_i) = \{u \in \Lambda(\Sr) \ ; \ \langle \lambda_i, u \rangle \geqslant \langle \lambda_j, u\rangle, \ \textrm{ for all j}\} = \{u \in \Lambda(\Sr) \ ; \ \langle \lambda_i - \lambda_j, u \rangle \geqslant 1, \ \textrm{ for all j}\}.$$ Given a basis $\Delta$ of $\Phi(\G, \Sr) \subset \X^*(\Sr)$, we denote by $\rP^{\Delta}_{\varnothing}$ the corresponding minimal parabolic subgroup of $\G$ containing $\Sr$ and by $\lambda_0(\Delta)$ the highest $k$-weight of $\rho$ with respect to $\rP^{\Delta}_{\varnothing}$; we also recall that the Weyl cone $\mathfrak{C}(\rP_{\varnothing}^{\Delta})$ is defined by the conditions $\alpha \geqslant 1$ for all $\alpha \in \Delta$. If $\lambda_0(\Delta) = \lambda_i$, then $\lambda_i - \lambda_j$ is a linear combination with non-negative coefficients of elements of $\Delta$ and thus $\iota^{-1}(\C_i)$ contains $\mathfrak{C}(\rP_{\varnothing}^{\Delta})$. Therefore, it follows from the proof of Lemma 2.5 that $\iota^{-1}(\C_i)$ contains the cone $$\C_{\tau}(\rP_{\varnothing}^{\Delta}) = \bigcup_{{\tiny \begin{array}{c} \Delta' \\ \lambda_0(\Delta') = \lambda_0(\Delta) \end{array}}}\mathfrak{C}(\rP^{\Delta'}_{\varnothing})$$ if $\lambda_i$ is the highest weight of $\rho$ with respect to $\rP^{\Delta}_{\varnothing}$.

The inclusion $\C_{\tau}(\rP_{\varnothing}^{\Delta}) \subset \iota^{-1}(\C_i)$ is in fact an equality. If it were not, then $\iota^{-1}(\C_i)$ would meet the interior of some Weyl cone $\mathfrak{C}(\rP_{\varnothing}^{\Delta'})$ with $\lambda_0(\Delta') \neq \lambda_i$. Setting $\lambda_0(\Delta') = \lambda_j$, it would follow that $\iota^{-1}(\C_i \cap \C_j)$ contains a point $x$ of $\mathfrak{C}(\rP_{\varnothing}^{\Delta'})^{\circ}$. Such a situation cannot happen: on the one hand, $\iota(x) \in \C_i \cap \C_j$ implies $\lambda_i(x) = \lambda_j(x)$; on the other hand, $\lambda_j - \lambda_i$ is a non-zero linear combination with non-negative coefficients of elements of $\Delta'$, hence $\lambda_j - \lambda_i >1$ on $\mathfrak{C}(\rP_{\varnothing}^{\Delta'})^{\circ}$ and thus $\lambda_j(x) > \lambda_i(x)$.

We have therefore $\iota^{-1}(\C_i) = \C_{\tau}(\rP^{\Delta}_\varnothing)$ if $\lambda_i$ is the highest $k$-weight of $\rho$ with respect to $\rP_{\varnothing}^{\Delta}$, whereas $\iota^{-1}(\C_i^\circ)$ is empty if $\lambda_i$ doesn't occur among the highest $k$-weights of $\rho$. We have checked that the fans $\mathcal{F}_\tau$ and $\iota^{-1}(\mathcal{F})$ have the same cones of maximal dimension; since each face is the intersection of suitable cones of maximal dimension, it follows that $\mathcal{F}_{\tau} = \iota^{-1}(\mathcal{F})$. \hfill $\Box$

\vspace{0.1cm}
By the preceding lemma, the affine map $i : \A(\Sr,k) \rightarrow \A(\T,k)$ can be extended to a continuous injective map
$$i : \overline{\A}_{\tau}(\Sr,k) \longrightarrow \overline{\A}_{\pi}(\T,k)$$
which is a homeomorphism of $\overline{\A}_{\tau}(\Sr,k)$ onto the closure of $i(\A(\Sr,k))$ in $\overline{\A}_{\pi}(\T,k)$.

\vspace{0.2cm}
\noindent \textbf{(5.3)} As recalled in Section 2, $\overline{\mathcal{B}}_{\tau}(\G,k) \cong \overline{\mathcal{B}}(\G,k)_{\rho}$ can be described as the quotient of $\G(k) \times \overline{\A}_{\tau}(\Sr,k)$ by the following equivalence relation:
\begin{eqnarray*}
(g,x) \sim (h,y) & \mbox{if and only if  there exists an element } n \in \N(k) \\
~ & \mbox{such that } nx= y \mbox{ and } g \inv h n \in \rP_x.
\end{eqnarray*}

Here $\rP_x$ is defined as $\rP_x= \N(k)_x \U_x$, where $\N(k)_x$ is the subgroup of $\N(k)$ fixing $x$, and where $\U_x$ is generated by all filtration steps $\U_{\alpha}(k)_{-\log \widetilde{\alpha}(x)}$ in the root group $\U_\alpha(k)$, with
$$\widetilde{\alpha}(x) = \sup\{c \in \mathbb{R}_{>0}: x \in \overline{\{\alpha(\cdot-o) \geqslant c\}}.$$
Similarly, $\overline{\mathcal{B}}_\pi({\rm SL}_\V,k)$ can be described as the quotient of ${\rm SL}_\V (k) \times \overline{\A}_\pi(\T,k)$ with respect to the analogous equivalence relation involving the stabilizer groups $\rP'_{x}$ for $x \in \overline{\A}_\pi(\T,k)$.

Composing $i: \overline{\A}_{\tau}(\Sr,k) \rightarrow \overline{\A}_\pi(\T,k)$ with the embedding of $\overline{\A}_\pi(\T,k)$ in $\overline{\mathcal{B}}_\pi({\rm SL}_\V,k)$, we obtain a continuous and therefore $\N(k)$-equivariant map $\overline{\rho}_\ast: \overline{\A}_{\tau}(\Sr,k) \rightarrow \overline{\mathcal{B}}_\pi({\rm SL}_\V,k)$.

\vspace{0.1cm}
Now we want to continue this map to the compactified building $\overline{\mathcal{B}}_{\tau}(\G,k)$.
\begin{Lemma}
For every $x \in \overline{\A}_{\tau}(\Sr,k)$ we have $\rho(\rP_x) \subset \rP'_{i(x)}$, where $\rP_x$ denotes the stabilizer of $x$ in $\G(k)$ and $\rP'_{i(x)}$ denotes the stabilizer of $i(x)$ in ${\rm SL}_\V (k)$.
\end{Lemma}

\textbf{\emph{Proof}}. If $x \in \A(\Sr,k)$, the claim holds by (5.1), property 2. In general, we have $\rP_x = \U(k)_x \N(k)_x$ where $\N(k)_x$ is the stabilizer of $x$ in $\N(k)$. Since $\rho_\ast: \overline{\A}_{\tau}(\Sr,k) \rightarrow \overline{\mathcal{B}}_{\pi}({\rm SL}_\V,k)$ is $\N(k)$-equivariant, we find $\rho(\N(k)_x) \subset \rP'_{i(x)}$. The group $\U(k)_x$ is generated by all $\U_{\alpha}(k)_x = \U_{\alpha}(k)_{-\log \widetilde{\alpha}(x)}$ for $\alpha \in \Phired$. Hence it suffices to show $\rho(\U_{\alpha}(k)_x) \subset \rP'_{i(x)}$ for all $\alpha \in \Phired$.

If $0 < \widetilde{\alpha}(x) < \infty$, then there exists a sequence $(x_n)$ of points in $\A(\Sr,k)$ converging towards $x$ and such that $\widetilde{\alpha}(x) = \alpha(x_n)$ for all $n$, hence $\U_\alpha(k)_x = \U_\alpha(k)_{x_n}$ for all $n$. By (5.1), property 2, it follows that $\rho(\U_{\alpha}(k)_x) \subset  \rho(\rP_{x_n})$ is contained in $\rP'_{i(x_n)}$. Since $i(x_n)$ converges towards $i(x)$ and ${\rm SL}_\V (k)$ acts continuously on $\overline{\mathcal{B}}_\pi({\rm SL}_\V, k)$, this implies $\rho(\U_{\alpha}(k)_x) \subset \rP'_{i(x)}$.

If $\widetilde{\alpha}(x) = 0$, then $\U_{\alpha}(k)_x = \{1\}$ and there is nothing to prove.

It remains to address the case where $\widetilde{\alpha}(x) = \infty$, hence $\U_{\alpha}(k)_x = \U_\alpha(k)$. There exists a sequence $(x_n)$ of points in $\A(\Sr,k)$ converging to $x$ and such that $\lim \alpha(x_n) = \infty$ (observe that $x$ belongs to the closure of each half-space $\{\alpha(\cdot - o) \geq c\}$, with $c \in \mathbb{R}_{\geq 0}$). Any element $u$ of $\U_\alpha(k)$ lies in one of the filtration steps $\U_{\alpha}(k)_r$; since this filtration is decreasing, $u$ belongs to $\U_{\alpha}(k)_{x_n}$, hence to the stabilizer $\rP_{x_n}$, if $n$ is big enough. By Landvogt's results, this implies that $\rho(u)$ is contained in $\rP'_{i(x_n)}$ for $n$ big enough. Since ${\rm SL}_\V(k)$ acts continuously on $\overline{\mathcal{B}}_\pi({\rm SL}_\V,k)$, it follows that $\rho(u)$ is indeed contained in $\rP'_{i(x)}$ and the proof is complete.\hfill$\Box$

\vspace{0.1cm}
It follows immediately from the lemma above that the natural $\G(k)$-equivariant map $$\G(k) \times \overline{\A}_{\tau}(\Sr,k) \rightarrow \overline{\mathcal{B}}_\pi({\rm SL}_\V,k), \ \ (g,x) \mapsto \rho(g) \cdot \overline{\rho}_\ast(x)$$ factors through the equivalence relation defining $\overline{\mathcal{B}}_{\tau}(\G,k)$ and thus induces a $\G(k)$-equivariant and continuous map $$\overline{\rho}_* : \overline{\mathcal{B}}_{\tau}(\G,k) \rightarrow \overline{\mathcal{B}}_\pi({\rm SL}_\V,k)$$ extending Landvogt's map $\rho_\ast$.

\begin{Thm}
\label{th-Satake_via_embeddings}
The map $\overline{\rho}_\ast: \overline{\mathcal{B}}_{\tau}(\G,k) \rightarrow \overline{\mathcal{B}}_\pi({\rm SL}_\V,k)$ is a $\G(k)$-equivariant homeomorphism of $\overline{\mathcal{B}}_{\tau}(\G,k)$ onto the closure of $\rho_\ast(\mathcal{B}(\G,k))$ in $\overline{\mathcal{B}}_\pi({\rm SL}_\V,k)$.
\end{Thm}

\noindent \emph{\textbf{Proof}}. The image of the compact space $\overline{\mathcal{B}}_{\tau}(\G,k)$ under $\overline{\rho}_\ast$ is closed, hence it contains the closure of $\rho_\ast(\mathcal{B}(\G,k))$. On the other hand, any point $z$ in $\rho_\ast(\overline{\mathcal{B}}(\G,k)_\rho)$ is of the form $z = \rho(g) \cdot \overline{\rho}_\ast(x)$ for some $g \in \G$ and some $x \in \overline{\A}_{\tau}(\Sr,k)$. If $(x_n)$ is a sequence of points in $\A(\Sr,k)$ converging towards $x$, then $\left(\rho(g) \cdot \overline{\rho}_\ast(x_n)\right)$ is a sequence of points in $\rho_\ast(\mathcal{B}(\G,k))$ converging towards $z$, hence $z$ is contained in the closure of $\rho_\ast(\mathcal{B}(\G,k))$. Injectivity follows from the fact that any two points of $\overline{\mathcal{B}}_{\tau}(\G,k)$ are contained in one compactified apartment by Theorem \ref{thm.bruhat.decomposition} (i).

Therefore, the map $\overline{\rho}_\ast$ is a continuous bijection between $\overline{\mathcal{B}}_{\tau}(\G,k)$ and the closure of $\rho_\ast(\mathcal{B}(\G,k))$ in $\overline{\mathcal{B}}_\pi({\rm SL}_\V,k)$. Since both spaces are compact, this is a homeomorphism.\hfill$\Box$

\vspace{0.2cm} \noindent \textbf{(5.4)} We complete this work by identifying the $k$-rational type $\tau$ appearing in Theorems \ref{thm.comparison} and \ref{th-Satake_via_embeddings}.

\begin{Prop} \label{prop.types} The type $\tau$ is the unique $k$-rational type defining the Berkovich compactification $\overline{\mathcal{B}}_{t(\rho)}(\G,k)$. Equivalently, we have $$\overline{\mathcal{B}}(\G,k)_{\rho} \cong \overline{\mathcal{B}}_{t(\rho)}(\G,k)$$ and any Landvogt map $\rho_* : \mathcal{B}(\G,k) \rightarrow \mathcal{B}({\rm SL}_{\V},k)$ extends to a $\G(k)$-equivariant homeomorphism between $\overline{\mathcal{B}}_{t(\rho)}(\G,k)$ and a closed subspace of $\overline{\mathcal{B}}_{\pi}({\rm SL}_\V,k)$.
\end{Prop}

\vspace{0.1cm}
\noindent \textbf{Proof}. Applying Theorem \ref{thm.Berkovich.Satake} to the contragredient representation $\check{\rho}$, the Berkovich map $\underline{\check{\rho}}$ provides us with a $\G(k)$-homeomorphism between $\overline{\mathcal{B}}_{t(\rho)}(\G,k)$ and a closed subspace of $\overline{\mathcal{B}}_{\delta}({\rm SL}_{\V^{\vee}},k) \cong \overline{\mathcal{B}}_{\pi}({\rm SL}_\V,k)$. Since this map is toral (Proposition \ref{prop.Berkovich.toral}), it satisfies conditions 1 to 3 of (5.1) and we deduce from Theorem \ref{th-Satake_via_embeddings} that the compactifications $\overline{\mathcal{B}}_{t(\rho)}(\G,k)$ and $\overline{\mathcal{B}}_\tau(\G,k)$ are $\G(k)$-homeomorphic. Thus, $\tau$ is the unique $k$-rational type defining the same Berkovich compactification as the type $t(\rho)$ naturally attached to the absolutely irreducible representation $\rho$ (see \cite[Appendix C]{RTW1}). \hfill $\Box$

\bibliographystyle{amsalpha}
\bibliography{BerkovichSatake4}

\vspace{1cm}

\vspace{0.5cm}
\begin{flushleft} \textit{Bertrand R\'emy and Amaury Thuillier} \\
Universit\'e de Lyon \\
Universit\'e Lyon 1 \\
CNRS - UMR 5208 \\
Institut Camille Jordan \\
43 boulevard du 11 novembre 1918 \\
F-69622 Villeurbanne cedex \\
\vspace{1pt}
$\{$remy; thuillier$\}$@math.univ-lyon1.fr
\end{flushleft}

\vspace{0.1cm}
\begin{flushleft}
\textit{Annette Werner} \\
Institut f\"ur Mathematik \\
Goethe-Universit\"at Frankfurt \\
Robert-Maier-Str., 6-8 \\
D-60325 Frankfurt-am-Main \\

\vspace{1pt}
werner@mathematik.uni-frankfurt.de
\end{flushleft}

\end{document}